\documentclass[11pt]{article}
\usepackage{geometry}                % See geometry.pdf to learn the layout options. There are lots.
\usepackage{graphicx}
\usepackage{amssymb}
\usepackage{epstopdf}
\usepackage{amscd}
\newtheorem{lemma}{Lemma}[section] 
\newtheorem{propos}[lemma]{Proposition}
\newtheorem{example}[lemma]{Example}
\newtheorem{theorem}[lemma]{Theorem}

\newtheorem{defin}[lemma]{Definition}
\newtheorem{remark}[lemma]{Remark}

\newcommand{\C}{\mathbb{C}}

\newcommand{\Hom}{\mathrm{Hom}}
\newcommand{\Vec}{\mathrm{Vec}}

\newcommand{\ev}{\mathrm{ev}}
\newcommand{\coev}{\mathrm{coev}}

\newcommand{\extd}{\mathrm{d}}

\newcommand{\tens}{\mathop{\otimes}}
\newcommand{\la}{{\triangleright}}

\newcommand{\id}{\mathrm{id}}

\newcommand{\<}{\langle}
\renewcommand{\>}{\rangle}

\newcounter{piccie}
\title{Noncommutative differential operators, Sobolev spaces and the centre of a category}
\author{E.J.\ Beggs \& Tomasz Brzezi\'nski\\
Department of Mathematics, Swansea University, Wales.}
\begin{document}

\maketitle
%\section{}
%\subsection{}

\abstract{
We consider differential operators over a noncommutative algebra $A$ generated by vector fields. These are shown to form a unital associative algebra of differential operators, and act on $A$-modules $E$ with covariant derivative. We use the repeated differentials given in the paper to give a definition of noncommutative Sobolev space for modules with connection and Hermitian inner product. The tensor algebra of vector fields, with a modified bimodule structure and a bimodule connection, is shown to lie in the centre of the bimodule connection category ${}_A\mathcal{E}_A$, and in fact to be an algebra in the centre. The crossing natural transformation in the definition of the centre of the category is related to the action of the differential operators on bimodules with connection. 
}

\section{Introduction}
The reader will be familiar with the local description of differential operators involving partial derivatives along coordinate directions on a manifold. This can be modified to be given in terms of derivatives along globally defined vector fields. In this paper we will consider the noncommutative analogue of such differential operators. We will not consider pseudo-differential or fractional differential operators. 

Take a possibly noncommutative algebra $A$, with differential structure given by a differential graded algebra $\Omega^n A$ (with $\Omega^0A=A$), differential
$\extd:\Omega^n A\to \Omega^{n+1} A$ (with $\extd^2=0$) and product
$\wedge$ (see \cite{connesbook} for a more detailed description). We shall assume that $\Omega^1 A$ is finitely generated projective as a right $A$ module, and set the vector fields $\Vec A=\Hom_A(\Omega^1 A,A)$, the right module maps from $\Omega^1 A$ to $A$. These vector fields will not act on the algebra $A$ as derivations in general. However this dual idea is just what corresponds to vector fields in commutative geometry. 

We take the tensor algebra of vector fields,
\begin{eqnarray*}
\mathcal{T}\,\Vec A\ =\ A\bigoplus \Vec A \bigoplus \Vec A\tens_A\Vec A
\bigoplus \Vec A\tens_A\Vec A\tens_A\Vec A \bigoplus \dots
\end{eqnarray*}
and give it a new associative product $\bullet$, involving differentiating the vector fields. This new algebra $\mathcal{T}\,\Vec A_\bullet$ then acts on modules with covariant derivative by repeated differentiation. This is the content of one of the main results, Theorem \ref{kuykjvcvu}. 

One feature of this paper is the number of covariant derivatives used, but this should not be too surprising. Given a left module $E$ with a left covariant derivative $\nabla_E:E\to \Omega^1 A\tens_A E$, if we apply the covariant derivative twice
(without using the wedge product, which would give curvature) we have to differentiate elements of $\Omega^1 A$, and the most convenient way to do that is by a right 
covariant derivative $\square:\Omega^1 A\to \Omega^1 A\tens_A\Omega^1 A$. The vector fields $\Vec A$ are dually paired with $\Omega^1 A$, and we obtain an induced left covariant derivative $\square:\Vec A\to \Omega^1 A\tens_A \Vec A$.
(It will be convenient to overload the symbol $\square$ and distinguish these just by their domain.) We also have to assume that $\square$ is a bimodule covariant derivative. This idea had its origins in \cite{DVMic,DVMass} and \cite{Mourad}, and was later 
used in \cite{Madore,FioMad}. In \cite{BMHDSfinitegrps} it was shown that this idea allowed tensoring
of bimodules with connections. 

We use the category ${}_A\mathcal{E}_A$ of $A$-bimodules with bimodule covariant derivative (see Definition \ref{cvgjfxtzhc}). 
The $A$-bimodule $\mathcal{T}\Vec A_\bullet$ with the bimodule 
covariant derivative given in Proposition \ref{bchudsabvuyas} is in the centre
$\mathcal{Z}({}_A\mathcal{E}_A)$ of ${}_A\mathcal{E}_A$, using the 
natural transformation 
$\vartheta:\mathcal{T}\Vec A_\bullet\tens_A-\Rightarrow -\tens_A \mathcal{T}\Vec A_\bullet$ defined in Proposition \ref{jkytfryjhzzdt}. In addition $\mathcal{T}\Vec A_\bullet$ is a unital associative algebra in $\mathcal{Z}({}_A\mathcal{E}_A)$, using the previously mentioned product $\bullet$. 
The natural transformation $\vartheta$ is related to the action of differential fields on modules with connection by Fig.\ \ref{newsig56svvuyd}. This is the content of the other of the main results, Theorem \ref{cbadiosiuvcuycf}. 

 We can view classical vector fields as the Lie algebra of the diffeomorphism group, or `infinitesimal' diffeomorphisms. Then we can give a coproduct on the tensor algebra over the Lie algebra, and make a bialgebra. It might be expected that a similar construction in the noncommutative case would give rise to a Hopf algebroid \cite {Boh:Hop}, or something similar. However those familiar with the theory of Hopf algebroids will be expecting complications. We have given the crossing natural transformation $\vartheta$ and the centre of the category construction as an alternative method to describe the action of differential operators on tensor products of bimodules with bimodule connection. It is not at all obvious whether a coproduct can be made to work, and we would be happy to hear from other authors on the matter. One interesting coproduct, though not quite what is needed in this case, is the braided shuffle coproduct \cite{Ros:gro}. 
 
In Section \ref{lkjsbvkcyxth} we define noncommutative Sobolev spaces for modules with connection, using Hermitian inner products. [Here make the assumption, not required elsewhere, that the algebra $A$ is a dense subalgebra of a $C^*$-algebra.] The reader who is familiar with Sobolev spaces will probably (and correctly) think that this is the easy part. However it is our hope that this definition will prompt people to look at the difficult part, the Sobolev embedding results. The interesting part here is that the dimension of the manifold explicitly appears in the classical Sobolev embedding results. Presumably a noncommutative version would have to make reference to a dimension of the algebra $A$, and there are quite a few different definitions of dimension to choose from.

It would be rather dishonest if we were to imply that many of the results were originally proved in the line by line form used for proofs here, it was much more usual to originally work with a diagrammatic form, and we give several figures to illustrate this. This form is basically the usual diagrams for the monoidal category of $A$-bimodules with $\tens_A$ as the tensor operation. 
The reader is reminded that the tensor product $E\tens_A F$ is the tensor product over the field $E\tens_{\mathbb{C}} F$, quotiented by the additional relations
\begin{eqnarray*}
e.a\tens f-e\tens a.f\ =\ 0\ ,
\end{eqnarray*}
for all $e\in E$, $f\in F$ and $a\in A$. It is the appropriate tensor product to use to replace the pointwise tensor product of sections of vector bundles in the commutative case. If $A=C(X)$ (the continuous functions on $X$) and $E$ and $F$ are modules of sections on vector bundles over $X$, then $E\tens_A F$ is the sections of the tensor product bundle. Now, many of the diagrams or equations we will use involve operations which are not bimodule maps, and therefore do not sit neatly within this bimodule tensor category formulation. The best example of this is in Proposition~\ref{tens45}, where we have an operation $\nabla_{E\tens F}$ on a tensor product
$E\tens_A F$, which is a sum of two operations, but neither of these two operations is well defined on $E\tens_A F$ (i.e.\ they may differ when applied to $e.a\tens f$ and $e\tens a.f$). The reader should think of these operations as being applied to a fixed element of the tensor product over the field $E\tens_{\mathbb{C}} F$, with operations being lifted to operations on  the tensor product over the field also. In this manner each diagram or equation can be given a meaning, and the problem of well definition over $\tens_A$ checked later. 

The reader should note that we continue with the bimodule covariant derivative and finitely generated projective assumption for $\Omega^1 A$ stated in Section~\ref{mghchgxa}, and the consequences for $\Vec A$, throughout the rest of the paper. For notation, the reader should note that we use the `generalised braidings' $\sigma$ and $\sigma^{-1}$ in a manner consistent with the notation on braided categories, so we may define $\sigma^{-1}$ as the basic object, and
denote its inverse (if it exists) by $\sigma$. Similarly we often use the symbols $\sigma$ and $\sigma^{-1}$ to stand for different maps, which are distinguished by their domains. This overloading of symbols is quite well defined, and was found to be preferable to adding many subscripts to the formulae. 

The authors would like to thank N.C.\ Phillips for useful advice on operator algebras, and S.\ Majid and R.\ Street for very helpful comments regarding category theory.

\section{Preliminaries}  \label{mghchgxa}

\subsection{Finitely generated projective modules}\label{vcxgaksujhgc}
Let $A$ be a possibly noncommutative unital algebra. 
 The \emph{(right) dual} $F'$ of a right
$A$-module $F$ is defined to be $\Hom_A(F,A)$, the right module maps from
$F$ to $A$.
Then $F'$ has a left module structure given by
$(a.\alpha)(f)=a.\alpha(f)$ for all
 $\alpha\in F'$, $a\in A$ and $f\in F$.

\begin{defin}\label{canonyy} A right $A$-module $F$ is said to be
\emph{finitely 
    generated projective} if
there are $f^i\in F$ and $f_i\in F'$ (for integer $1\le i\le n$)
(the `dual basis') so that for all $f\in F$, $f=\sum f^i.f_i(f)$. 
From this it follows directly that $\alpha=\sum \alpha(f^i).f_i$
for all $\alpha\in F'$. The $A$ valued matrix $P_{qj}=f_q(f^j)$
is an idempotent associated to the module.
\end{defin}

If $F$ is a bimodule, then $F'$ is also a bimodule, with left module
structure given previously, and right module structure
$(\alpha.a)(f)=\alpha(a.f)$. If in addition $F$ is
finitely generated projective as a right $A$-module, we have bimodule maps
\begin{eqnarray*}
\ev:F'\tens_A F\to A\ ,\quad \coev:A\to F\tens_A F'\ ,
\end{eqnarray*}
given by $\ev(\alpha\tens f)=\alpha(f)$ and $\coev(1_A)=\sum_i f^i\tens
f_i$. These obey the identities
\begin{eqnarray} \label{xerakjcgh}
(\ev\tens\id)(\id\tens\coev(1_A))\ =\ \id:F' &\to & F'\ ,\cr
(\id\tens\ev)(\coev(1_A)\tens\id)\ =\ \id:F &\to & F\ .
\end{eqnarray}
 
 In this paper we shall use the usual diagrammatic notation for tensor or monoidal categories. Here vertical lines represent the identity maps between objects (modules in our case), and tensor product is given by writing the vertical lines next to each other. The diagrams are meant to be read from top to bottom.

\unitlength 0.5 mm
\begin{picture}(220,65)(-20,30)
\linethickness{0.3mm}
\multiput(109.99,60.5)(0.01,-0.5){1}{\line(0,-1){0.5}}
\multiput(109.95,61)(0.04,-0.5){1}{\line(0,-1){0.5}}
\multiput(109.89,61.49)(0.06,-0.49){1}{\line(0,-1){0.49}}
\multiput(109.8,61.98)(0.09,-0.49){1}{\line(0,-1){0.49}}
\multiput(109.69,62.47)(0.11,-0.49){1}{\line(0,-1){0.49}}
\multiput(109.56,62.95)(0.14,-0.48){1}{\line(0,-1){0.48}}
\multiput(109.4,63.42)(0.16,-0.47){1}{\line(0,-1){0.47}}
\multiput(109.21,63.88)(0.09,-0.23){2}{\line(0,-1){0.23}}
\multiput(109.01,64.34)(0.1,-0.23){2}{\line(0,-1){0.23}}
\multiput(108.78,64.78)(0.11,-0.22){2}{\line(0,-1){0.22}}
\multiput(108.53,65.21)(0.12,-0.22){2}{\line(0,-1){0.22}}
\multiput(108.26,65.63)(0.14,-0.21){2}{\line(0,-1){0.21}}
\multiput(107.97,66.04)(0.15,-0.2){2}{\line(0,-1){0.2}}
\multiput(107.66,66.43)(0.1,-0.13){3}{\line(0,-1){0.13}}
\multiput(107.33,66.8)(0.11,-0.12){3}{\line(0,-1){0.12}}
\multiput(106.98,67.16)(0.12,-0.12){3}{\line(0,-1){0.12}}
\multiput(106.62,67.5)(0.12,-0.11){3}{\line(1,0){0.12}}
\multiput(106.23,67.82)(0.13,-0.11){3}{\line(1,0){0.13}}
\multiput(105.84,68.12)(0.13,-0.1){3}{\line(1,0){0.13}}
\multiput(105.43,68.4)(0.21,-0.14){2}{\line(1,0){0.21}}
\multiput(105,68.66)(0.21,-0.13){2}{\line(1,0){0.21}}
\multiput(104.56,68.9)(0.22,-0.12){2}{\line(1,0){0.22}}
\multiput(104.11,69.12)(0.22,-0.11){2}{\line(1,0){0.22}}
\multiput(103.65,69.31)(0.23,-0.1){2}{\line(1,0){0.23}}
\multiput(103.18,69.48)(0.47,-0.17){1}{\line(1,0){0.47}}
\multiput(102.71,69.63)(0.48,-0.15){1}{\line(1,0){0.48}}
\multiput(102.23,69.75)(0.48,-0.12){1}{\line(1,0){0.48}}
\multiput(101.74,69.85)(0.49,-0.1){1}{\line(1,0){0.49}}
\multiput(101.24,69.92)(0.49,-0.07){1}{\line(1,0){0.49}}
\multiput(100.75,69.97)(0.5,-0.05){1}{\line(1,0){0.5}}
\multiput(100.25,70)(0.5,-0.02){1}{\line(1,0){0.5}}
\put(99.75,70){\line(1,0){0.5}}
\multiput(99.25,69.97)(0.5,0.02){1}{\line(1,0){0.5}}
\multiput(98.76,69.92)(0.5,0.05){1}{\line(1,0){0.5}}
\multiput(98.26,69.85)(0.49,0.07){1}{\line(1,0){0.49}}
\multiput(97.77,69.75)(0.49,0.1){1}{\line(1,0){0.49}}
\multiput(97.29,69.63)(0.48,0.12){1}{\line(1,0){0.48}}
\multiput(96.82,69.48)(0.48,0.15){1}{\line(1,0){0.48}}
\multiput(96.35,69.31)(0.47,0.17){1}{\line(1,0){0.47}}
\multiput(95.89,69.12)(0.23,0.1){2}{\line(1,0){0.23}}
\multiput(95.44,68.9)(0.22,0.11){2}{\line(1,0){0.22}}
\multiput(95,68.66)(0.22,0.12){2}{\line(1,0){0.22}}
\multiput(94.57,68.4)(0.21,0.13){2}{\line(1,0){0.21}}
\multiput(94.16,68.12)(0.21,0.14){2}{\line(1,0){0.21}}
\multiput(93.77,67.82)(0.13,0.1){3}{\line(1,0){0.13}}
\multiput(93.38,67.5)(0.13,0.11){3}{\line(1,0){0.13}}
\multiput(93.02,67.16)(0.12,0.11){3}{\line(1,0){0.12}}
\multiput(92.67,66.8)(0.12,0.12){3}{\line(0,1){0.12}}
\multiput(92.34,66.43)(0.11,0.12){3}{\line(0,1){0.12}}
\multiput(92.03,66.04)(0.1,0.13){3}{\line(0,1){0.13}}
\multiput(91.74,65.63)(0.15,0.2){2}{\line(0,1){0.2}}
\multiput(91.47,65.21)(0.14,0.21){2}{\line(0,1){0.21}}
\multiput(91.22,64.78)(0.12,0.22){2}{\line(0,1){0.22}}
\multiput(90.99,64.34)(0.11,0.22){2}{\line(0,1){0.22}}
\multiput(90.79,63.88)(0.1,0.23){2}{\line(0,1){0.23}}
\multiput(90.6,63.42)(0.09,0.23){2}{\line(0,1){0.23}}
\multiput(90.44,62.95)(0.16,0.47){1}{\line(0,1){0.47}}
\multiput(90.31,62.47)(0.14,0.48){1}{\line(0,1){0.48}}
\multiput(90.2,61.98)(0.11,0.49){1}{\line(0,1){0.49}}
\multiput(90.11,61.49)(0.09,0.49){1}{\line(0,1){0.49}}
\multiput(90.05,61)(0.06,0.49){1}{\line(0,1){0.49}}
\multiput(90.01,60.5)(0.04,0.5){1}{\line(0,1){0.5}}
\multiput(90,60)(0.01,0.5){1}{\line(0,1){0.5}}

\linethickness{0.3mm}
\multiput(110,60)(0.01,-0.5){1}{\line(0,-1){0.5}}
\multiput(110.01,59.5)(0.04,-0.5){1}{\line(0,-1){0.5}}
\multiput(110.05,59)(0.06,-0.49){1}{\line(0,-1){0.49}}
\multiput(110.11,58.51)(0.09,-0.49){1}{\line(0,-1){0.49}}
\multiput(110.2,58.02)(0.11,-0.49){1}{\line(0,-1){0.49}}
\multiput(110.31,57.53)(0.14,-0.48){1}{\line(0,-1){0.48}}
\multiput(110.44,57.05)(0.16,-0.47){1}{\line(0,-1){0.47}}
\multiput(110.6,56.58)(0.09,-0.23){2}{\line(0,-1){0.23}}
\multiput(110.79,56.12)(0.1,-0.23){2}{\line(0,-1){0.23}}
\multiput(110.99,55.66)(0.11,-0.22){2}{\line(0,-1){0.22}}
\multiput(111.22,55.22)(0.12,-0.22){2}{\line(0,-1){0.22}}
\multiput(111.47,54.79)(0.14,-0.21){2}{\line(0,-1){0.21}}
\multiput(111.74,54.37)(0.15,-0.2){2}{\line(0,-1){0.2}}
\multiput(112.03,53.96)(0.1,-0.13){3}{\line(0,-1){0.13}}
\multiput(112.34,53.57)(0.11,-0.12){3}{\line(0,-1){0.12}}
\multiput(112.67,53.2)(0.12,-0.12){3}{\line(0,-1){0.12}}
\multiput(113.02,52.84)(0.12,-0.11){3}{\line(1,0){0.12}}
\multiput(113.38,52.5)(0.13,-0.11){3}{\line(1,0){0.13}}
\multiput(113.77,52.18)(0.13,-0.1){3}{\line(1,0){0.13}}
\multiput(114.16,51.88)(0.21,-0.14){2}{\line(1,0){0.21}}
\multiput(114.57,51.6)(0.21,-0.13){2}{\line(1,0){0.21}}
\multiput(115,51.34)(0.22,-0.12){2}{\line(1,0){0.22}}
\multiput(115.44,51.1)(0.22,-0.11){2}{\line(1,0){0.22}}
\multiput(115.89,50.88)(0.23,-0.1){2}{\line(1,0){0.23}}
\multiput(116.35,50.69)(0.47,-0.17){1}{\line(1,0){0.47}}
\multiput(116.82,50.52)(0.48,-0.15){1}{\line(1,0){0.48}}
\multiput(117.29,50.37)(0.48,-0.12){1}{\line(1,0){0.48}}
\multiput(117.77,50.25)(0.49,-0.1){1}{\line(1,0){0.49}}
\multiput(118.26,50.15)(0.49,-0.07){1}{\line(1,0){0.49}}
\multiput(118.76,50.08)(0.5,-0.05){1}{\line(1,0){0.5}}
\multiput(119.25,50.03)(0.5,-0.02){1}{\line(1,0){0.5}}
\put(119.75,50){\line(1,0){0.5}}
\multiput(120.25,50)(0.5,0.02){1}{\line(1,0){0.5}}
\multiput(120.75,50.03)(0.5,0.05){1}{\line(1,0){0.5}}
\multiput(121.24,50.08)(0.49,0.07){1}{\line(1,0){0.49}}
\multiput(121.74,50.15)(0.49,0.1){1}{\line(1,0){0.49}}
\multiput(122.23,50.25)(0.48,0.12){1}{\line(1,0){0.48}}
\multiput(122.71,50.37)(0.48,0.15){1}{\line(1,0){0.48}}
\multiput(123.18,50.52)(0.47,0.17){1}{\line(1,0){0.47}}
\multiput(123.65,50.69)(0.23,0.1){2}{\line(1,0){0.23}}
\multiput(124.11,50.88)(0.22,0.11){2}{\line(1,0){0.22}}
\multiput(124.56,51.1)(0.22,0.12){2}{\line(1,0){0.22}}
\multiput(125,51.34)(0.21,0.13){2}{\line(1,0){0.21}}
\multiput(125.43,51.6)(0.21,0.14){2}{\line(1,0){0.21}}
\multiput(125.84,51.88)(0.13,0.1){3}{\line(1,0){0.13}}
\multiput(126.23,52.18)(0.13,0.11){3}{\line(1,0){0.13}}
\multiput(126.62,52.5)(0.12,0.11){3}{\line(1,0){0.12}}
\multiput(126.98,52.84)(0.12,0.12){3}{\line(0,1){0.12}}
\multiput(127.33,53.2)(0.11,0.12){3}{\line(0,1){0.12}}
\multiput(127.66,53.57)(0.1,0.13){3}{\line(0,1){0.13}}
\multiput(127.97,53.96)(0.15,0.2){2}{\line(0,1){0.2}}
\multiput(128.26,54.37)(0.14,0.21){2}{\line(0,1){0.21}}
\multiput(128.53,54.79)(0.12,0.22){2}{\line(0,1){0.22}}
\multiput(128.78,55.22)(0.11,0.22){2}{\line(0,1){0.22}}
\multiput(129.01,55.66)(0.1,0.23){2}{\line(0,1){0.23}}
\multiput(129.21,56.12)(0.09,0.23){2}{\line(0,1){0.23}}
\multiput(129.4,56.58)(0.16,0.47){1}{\line(0,1){0.47}}
\multiput(129.56,57.05)(0.14,0.48){1}{\line(0,1){0.48}}
\multiput(129.69,57.53)(0.11,0.49){1}{\line(0,1){0.49}}
\multiput(129.8,58.02)(0.09,0.49){1}{\line(0,1){0.49}}
\multiput(129.89,58.51)(0.06,0.49){1}{\line(0,1){0.49}}
\multiput(129.95,59)(0.04,0.5){1}{\line(0,1){0.5}}
\multiput(129.99,59.5)(0.01,0.5){1}{\line(0,1){0.5}}

\linethickness{0.3mm}
\put(90,40){\line(0,1){20}}
\linethickness{0.3mm}
\put(130,60){\line(0,1){20}}
\linethickness{0.3mm}
\put(150,40){\line(0,1){40}}
\linethickness{0.3mm}
\multiput(39.99,60.5)(0.01,-0.5){1}{\line(0,-1){0.5}}
\multiput(39.95,61)(0.04,-0.5){1}{\line(0,-1){0.5}}
\multiput(39.89,61.49)(0.06,-0.49){1}{\line(0,-1){0.49}}
\multiput(39.8,61.98)(0.09,-0.49){1}{\line(0,-1){0.49}}
\multiput(39.69,62.47)(0.11,-0.49){1}{\line(0,-1){0.49}}
\multiput(39.56,62.95)(0.14,-0.48){1}{\line(0,-1){0.48}}
\multiput(39.4,63.42)(0.16,-0.47){1}{\line(0,-1){0.47}}
\multiput(39.21,63.88)(0.09,-0.23){2}{\line(0,-1){0.23}}
\multiput(39.01,64.34)(0.1,-0.23){2}{\line(0,-1){0.23}}
\multiput(38.78,64.78)(0.11,-0.22){2}{\line(0,-1){0.22}}
\multiput(38.53,65.21)(0.12,-0.22){2}{\line(0,-1){0.22}}
\multiput(38.26,65.63)(0.14,-0.21){2}{\line(0,-1){0.21}}
\multiput(37.97,66.04)(0.15,-0.2){2}{\line(0,-1){0.2}}
\multiput(37.66,66.43)(0.1,-0.13){3}{\line(0,-1){0.13}}
\multiput(37.33,66.8)(0.11,-0.12){3}{\line(0,-1){0.12}}
\multiput(36.98,67.16)(0.12,-0.12){3}{\line(0,-1){0.12}}
\multiput(36.62,67.5)(0.12,-0.11){3}{\line(1,0){0.12}}
\multiput(36.23,67.82)(0.13,-0.11){3}{\line(1,0){0.13}}
\multiput(35.84,68.12)(0.13,-0.1){3}{\line(1,0){0.13}}
\multiput(35.43,68.4)(0.21,-0.14){2}{\line(1,0){0.21}}
\multiput(35,68.66)(0.21,-0.13){2}{\line(1,0){0.21}}
\multiput(34.56,68.9)(0.22,-0.12){2}{\line(1,0){0.22}}
\multiput(34.11,69.12)(0.22,-0.11){2}{\line(1,0){0.22}}
\multiput(33.65,69.31)(0.23,-0.1){2}{\line(1,0){0.23}}
\multiput(33.18,69.48)(0.47,-0.17){1}{\line(1,0){0.47}}
\multiput(32.71,69.63)(0.48,-0.15){1}{\line(1,0){0.48}}
\multiput(32.23,69.75)(0.48,-0.12){1}{\line(1,0){0.48}}
\multiput(31.74,69.85)(0.49,-0.1){1}{\line(1,0){0.49}}
\multiput(31.24,69.92)(0.49,-0.07){1}{\line(1,0){0.49}}
\multiput(30.75,69.97)(0.5,-0.05){1}{\line(1,0){0.5}}
\multiput(30.25,70)(0.5,-0.02){1}{\line(1,0){0.5}}
\put(29.75,70){\line(1,0){0.5}}
\multiput(29.25,69.97)(0.5,0.02){1}{\line(1,0){0.5}}
\multiput(28.76,69.92)(0.5,0.05){1}{\line(1,0){0.5}}
\multiput(28.26,69.85)(0.49,0.07){1}{\line(1,0){0.49}}
\multiput(27.77,69.75)(0.49,0.1){1}{\line(1,0){0.49}}
\multiput(27.29,69.63)(0.48,0.12){1}{\line(1,0){0.48}}
\multiput(26.82,69.48)(0.48,0.15){1}{\line(1,0){0.48}}
\multiput(26.35,69.31)(0.47,0.17){1}{\line(1,0){0.47}}
\multiput(25.89,69.12)(0.23,0.1){2}{\line(1,0){0.23}}
\multiput(25.44,68.9)(0.22,0.11){2}{\line(1,0){0.22}}
\multiput(25,68.66)(0.22,0.12){2}{\line(1,0){0.22}}
\multiput(24.57,68.4)(0.21,0.13){2}{\line(1,0){0.21}}
\multiput(24.16,68.12)(0.21,0.14){2}{\line(1,0){0.21}}
\multiput(23.77,67.82)(0.13,0.1){3}{\line(1,0){0.13}}
\multiput(23.38,67.5)(0.13,0.11){3}{\line(1,0){0.13}}
\multiput(23.02,67.16)(0.12,0.11){3}{\line(1,0){0.12}}
\multiput(22.67,66.8)(0.12,0.12){3}{\line(0,1){0.12}}
\multiput(22.34,66.43)(0.11,0.12){3}{\line(0,1){0.12}}
\multiput(22.03,66.04)(0.1,0.13){3}{\line(0,1){0.13}}
\multiput(21.74,65.63)(0.15,0.2){2}{\line(0,1){0.2}}
\multiput(21.47,65.21)(0.14,0.21){2}{\line(0,1){0.21}}
\multiput(21.22,64.78)(0.12,0.22){2}{\line(0,1){0.22}}
\multiput(20.99,64.34)(0.11,0.22){2}{\line(0,1){0.22}}
\multiput(20.79,63.88)(0.1,0.23){2}{\line(0,1){0.23}}
\multiput(20.6,63.42)(0.09,0.23){2}{\line(0,1){0.23}}
\multiput(20.44,62.95)(0.16,0.47){1}{\line(0,1){0.47}}
\multiput(20.31,62.47)(0.14,0.48){1}{\line(0,1){0.48}}
\multiput(20.2,61.98)(0.11,0.49){1}{\line(0,1){0.49}}
\multiput(20.11,61.49)(0.09,0.49){1}{\line(0,1){0.49}}
\multiput(20.05,61)(0.06,0.49){1}{\line(0,1){0.49}}
\multiput(20.01,60.5)(0.04,0.5){1}{\line(0,1){0.5}}
\multiput(20,60)(0.01,0.5){1}{\line(0,1){0.5}}

\linethickness{0.3mm}
\multiput(0,60)(0.01,-0.5){1}{\line(0,-1){0.5}}
\multiput(0.01,59.5)(0.04,-0.5){1}{\line(0,-1){0.5}}
\multiput(0.05,59)(0.06,-0.49){1}{\line(0,-1){0.49}}
\multiput(0.11,58.51)(0.09,-0.49){1}{\line(0,-1){0.49}}
\multiput(0.2,58.02)(0.11,-0.49){1}{\line(0,-1){0.49}}
\multiput(0.31,57.53)(0.14,-0.48){1}{\line(0,-1){0.48}}
\multiput(0.44,57.05)(0.16,-0.47){1}{\line(0,-1){0.47}}
\multiput(0.6,56.58)(0.09,-0.23){2}{\line(0,-1){0.23}}
\multiput(0.79,56.12)(0.1,-0.23){2}{\line(0,-1){0.23}}
\multiput(0.99,55.66)(0.11,-0.22){2}{\line(0,-1){0.22}}
\multiput(1.22,55.22)(0.12,-0.22){2}{\line(0,-1){0.22}}
\multiput(1.47,54.79)(0.14,-0.21){2}{\line(0,-1){0.21}}
\multiput(1.74,54.37)(0.15,-0.2){2}{\line(0,-1){0.2}}
\multiput(2.03,53.96)(0.1,-0.13){3}{\line(0,-1){0.13}}
\multiput(2.34,53.57)(0.11,-0.12){3}{\line(0,-1){0.12}}
\multiput(2.67,53.2)(0.12,-0.12){3}{\line(0,-1){0.12}}
\multiput(3.02,52.84)(0.12,-0.11){3}{\line(1,0){0.12}}
\multiput(3.38,52.5)(0.13,-0.11){3}{\line(1,0){0.13}}
\multiput(3.77,52.18)(0.13,-0.1){3}{\line(1,0){0.13}}
\multiput(4.16,51.88)(0.21,-0.14){2}{\line(1,0){0.21}}
\multiput(4.57,51.6)(0.21,-0.13){2}{\line(1,0){0.21}}
\multiput(5,51.34)(0.22,-0.12){2}{\line(1,0){0.22}}
\multiput(5.44,51.1)(0.22,-0.11){2}{\line(1,0){0.22}}
\multiput(5.89,50.88)(0.23,-0.1){2}{\line(1,0){0.23}}
\multiput(6.35,50.69)(0.47,-0.17){1}{\line(1,0){0.47}}
\multiput(6.82,50.52)(0.48,-0.15){1}{\line(1,0){0.48}}
\multiput(7.29,50.37)(0.48,-0.12){1}{\line(1,0){0.48}}
\multiput(7.77,50.25)(0.49,-0.1){1}{\line(1,0){0.49}}
\multiput(8.26,50.15)(0.49,-0.07){1}{\line(1,0){0.49}}
\multiput(8.76,50.08)(0.5,-0.05){1}{\line(1,0){0.5}}
\multiput(9.25,50.03)(0.5,-0.02){1}{\line(1,0){0.5}}
\put(9.75,50){\line(1,0){0.5}}
\multiput(10.25,50)(0.5,0.02){1}{\line(1,0){0.5}}
\multiput(10.75,50.03)(0.5,0.05){1}{\line(1,0){0.5}}
\multiput(11.24,50.08)(0.49,0.07){1}{\line(1,0){0.49}}
\multiput(11.74,50.15)(0.49,0.1){1}{\line(1,0){0.49}}
\multiput(12.23,50.25)(0.48,0.12){1}{\line(1,0){0.48}}
\multiput(12.71,50.37)(0.48,0.15){1}{\line(1,0){0.48}}
\multiput(13.18,50.52)(0.47,0.17){1}{\line(1,0){0.47}}
\multiput(13.65,50.69)(0.23,0.1){2}{\line(1,0){0.23}}
\multiput(14.11,50.88)(0.22,0.11){2}{\line(1,0){0.22}}
\multiput(14.56,51.1)(0.22,0.12){2}{\line(1,0){0.22}}
\multiput(15,51.34)(0.21,0.13){2}{\line(1,0){0.21}}
\multiput(15.43,51.6)(0.21,0.14){2}{\line(1,0){0.21}}
\multiput(15.84,51.88)(0.13,0.1){3}{\line(1,0){0.13}}
\multiput(16.23,52.18)(0.13,0.11){3}{\line(1,0){0.13}}
\multiput(16.62,52.5)(0.12,0.11){3}{\line(1,0){0.12}}
\multiput(16.98,52.84)(0.12,0.12){3}{\line(0,1){0.12}}
\multiput(17.33,53.2)(0.11,0.12){3}{\line(0,1){0.12}}
\multiput(17.66,53.57)(0.1,0.13){3}{\line(0,1){0.13}}
\multiput(17.97,53.96)(0.15,0.2){2}{\line(0,1){0.2}}
\multiput(18.26,54.37)(0.14,0.21){2}{\line(0,1){0.21}}
\multiput(18.53,54.79)(0.12,0.22){2}{\line(0,1){0.22}}
\multiput(18.78,55.22)(0.11,0.22){2}{\line(0,1){0.22}}
\multiput(19.01,55.66)(0.1,0.23){2}{\line(0,1){0.23}}
\multiput(19.21,56.12)(0.09,0.23){2}{\line(0,1){0.23}}
\multiput(19.4,56.58)(0.16,0.47){1}{\line(0,1){0.47}}
\multiput(19.56,57.05)(0.14,0.48){1}{\line(0,1){0.48}}
\multiput(19.69,57.53)(0.11,0.49){1}{\line(0,1){0.49}}
\multiput(19.8,58.02)(0.09,0.49){1}{\line(0,1){0.49}}
\multiput(19.89,58.51)(0.06,0.49){1}{\line(0,1){0.49}}
\multiput(19.95,59)(0.04,0.5){1}{\line(0,1){0.5}}
\multiput(19.99,59.5)(0.01,0.5){1}{\line(0,1){0.5}}

\linethickness{0.3mm}
\put(40,40){\line(0,1){20}}
\linethickness{0.3mm}
\put(0,60){\line(0,1){20}}
\put(50,60){\makebox(0,0)[cc]{$=$}}

\linethickness{0.3mm}
\put(60,40){\line(0,1){40}}
\put(140,60){\makebox(0,0)[cc]{$=$}}

\put(25,75){\makebox(0,0)[cc]{$\coev$}}

\put(25,75){\makebox(0,0)[cc]{}}

\put(95,75){\makebox(0,0)[cc]{$\coev$}}

\put(115,45){\makebox(0,0)[cc]{$\ev$}}

\put(5,45){\makebox(0,0)[cc]{$\ev$}}

\put(0,85){\makebox(0,0)[cc]{$F'$}}

\put(58,85){\makebox(0,0)[cc]{$F'$}}

\put(38,35){\makebox(0,0)[cc]{$F'$}}

\put(58,35){\makebox(0,0)[cc]{$F'$}}

\put(128,85){\makebox(0,0)[cc]{$F$}}

\put(148,85){\makebox(0,0)[cc]{$F$}}

\put(88,35){\makebox(0,0)[cc]{$F$}}

%%%% Edwins Piccie Counter Stuff %%%%%
\refstepcounter{piccie} \label{someLabel}
\put(185,55){\makebox(0,0)[cc]{Fig.\ \arabic{piccie}}}
%%%%%%%%%%%%%%%%%%%%%%%

\put(148,35){\makebox(0,0)[cc]{$F$}}

\end{picture}

 We give the usual relation between the evaluation and coevaluation maps (already given in (\ref{xerakjcgh})) in pictorial form in Fig.\ \ref{someLabel}. 
Note that the algebra $A$ is represented by an `invisible line' in accordance with the usual diagrammatic notation of the unit element in monoidal categories.

\subsection{Differential calculi and covariant derivatives}

Let $A$ be a unital algebra over $\C$. Suppose that the algebra $A$ has a
differential structure $(\Omega A,\extd)$ in the sense of a differential
graded exterior algebra $\Omega A=\oplus_n\Omega^n A$ with $\extd$ increasing
degree by 1 and obeying the graded Leibniz rule, and $\extd^2=0$. We suppose that
$\Omega^1A$ generates the exterior algebra over $A$, and that $\Omega^1 A=A.\extd A$. The notion of a
covariant derivative in this context is standard \cite{connesIHES62}:

\begin{defin}\label{bcaiskyf} Given a left $A$-module $E$, a left $A$-covariant derivative
is a map $\nabla:E\to\Omega^{1}A\tens_{A}E$ which obeys the condition
$\nabla(a.e)=\extd a\tens e+a.\nabla e$
for all $e\in E$ and $a\in A$.
\end{defin}

In classical differential geometry there is no difference in whether we
multiply a section of a vector bundle by a function on the left or right. In the
noncommutative case there is a difference, and
for a bimodule we could require the Leibniz rule for both
left and right multiplication, but this would turn out to be too
restrictive. Instead,
following \cite{DVMic,DVMass} and \cite{Mourad}, we make the following
definition:

\begin{defin} \label{ppll}  A bimodule covariant derivative on an
 $A$-bimodule $E$ is a triple $(E,\nabla,\sigma)$,
where $\nabla:E\to\Omega^{1}A\tens_{A}E$
is a left $A$-covariant derivative, and $\sigma:E\tens_A\Omega^1 A\to
\Omega^1 A \tens_A E$ is a bimodule map obeying
\[
\nabla(e.a)\,=\,\nabla(e).a\,+\,\sigma(e\tens da)\ ,\qquad \forall\, e\in
E,\ a\in A
\]
\end{defin}

Now we consider one of the most immediate reasons to define the
bimodule covariant derivative, that is to have a covariant derivative on
tensor
products of bimodules. As mentioned in the introduction, this occurred in 
\cite{BMHDSfinitegrps}.

\begin{propos} \label{tens45}  Given $(E,\nabla_E,\sigma_E)$ a
bimodule covariant derivative on the bimodule $E$ and $\nabla_F$ a left
covariant derivative on the left module $F$, there is a left
$A$-covariant derivative on $E\tens_{A}F$ given by
\begin{eqnarray*}
\nabla_{E\tens F} &=& \nabla_E\tens
\id_{F}+(\sigma_E\tens\id_{F})(\id_{E}\tens\nabla_F)
\end{eqnarray*}
Further if $F$ is also an $A$-bimodule with a bimodule covariant
derivative $(\nabla_F,\sigma_F)$, then there is a
bimodule covariant
derivative $(\nabla_{E\tens_A F},\sigma_{E\tens_A F})$ on $E\tens_A F$ with
\begin{eqnarray*}
\sigma_{E\tens F} \,=\, (\sigma_{E}\tens
\id)(\id\tens \sigma_{F})\ .
\end{eqnarray*}
\end{propos}

\medskip To illustrate Proposition~\ref{tens45} it will be useful to introduce a picture, Fig.\ \ref{someLabel2}. This is a formula for $\nabla_{E\tens F}$ which illustrates Proposition~\ref{tens45} in a pictorial fashion, and should be read in the manner of a braid diagram, going down the page. Note that $\sigma_E:E\tens_A \Omega^1 A\to \Omega^1 A\tens_A E$ has been written in Fig.\ \ref{someLabel2} after the manner of a braid, so that we may distinguish it from its inverse (if it exists) $\sigma_E^{-1}:\Omega^1 A\tens_A E\to E\tens_A \Omega^1 A$, which would be written with the crossing the other way.

\unitlength 0.5 mm
\begin{picture}(110,70)(-20,24)\label{bdhwvvggcv}
\linethickness{0.3mm}
\qbezier(30,60)(29.89,64.25)(27.07,67.07)
\qbezier(27.07,67.07)(24.25,69.89)(20,70)
\qbezier(20,70)(15.75,69.89)(12.93,67.07)
\qbezier(12.93,67.07)(10.11,64.25)(10,60)
\linethickness{0.3mm}
\multiput(109.99,60.5)(0.01,-0.5){1}{\line(0,-1){0.5}}
\multiput(109.95,61)(0.04,-0.5){1}{\line(0,-1){0.5}}
\multiput(109.89,61.49)(0.06,-0.49){1}{\line(0,-1){0.49}}
\multiput(109.8,61.98)(0.09,-0.49){1}{\line(0,-1){0.49}}
\multiput(109.69,62.47)(0.11,-0.49){1}{\line(0,-1){0.49}}
\multiput(109.56,62.95)(0.14,-0.48){1}{\line(0,-1){0.48}}
\multiput(109.4,63.42)(0.16,-0.47){1}{\line(0,-1){0.47}}
\multiput(109.21,63.88)(0.09,-0.23){2}{\line(0,-1){0.23}}
\multiput(109.01,64.34)(0.1,-0.23){2}{\line(0,-1){0.23}}
\multiput(108.78,64.78)(0.11,-0.22){2}{\line(0,-1){0.22}}
\multiput(108.53,65.21)(0.12,-0.22){2}{\line(0,-1){0.22}}
\multiput(108.26,65.63)(0.14,-0.21){2}{\line(0,-1){0.21}}
\multiput(107.97,66.04)(0.15,-0.2){2}{\line(0,-1){0.2}}
\multiput(107.66,66.43)(0.1,-0.13){3}{\line(0,-1){0.13}}
\multiput(107.33,66.8)(0.11,-0.12){3}{\line(0,-1){0.12}}
\multiput(106.98,67.16)(0.12,-0.12){3}{\line(0,-1){0.12}}
\multiput(106.62,67.5)(0.12,-0.11){3}{\line(1,0){0.12}}
\multiput(106.23,67.82)(0.13,-0.11){3}{\line(1,0){0.13}}
\multiput(105.84,68.12)(0.13,-0.1){3}{\line(1,0){0.13}}
\multiput(105.43,68.4)(0.21,-0.14){2}{\line(1,0){0.21}}
\multiput(105,68.66)(0.21,-0.13){2}{\line(1,0){0.21}}
\multiput(104.56,68.9)(0.22,-0.12){2}{\line(1,0){0.22}}
\multiput(104.11,69.12)(0.22,-0.11){2}{\line(1,0){0.22}}
\multiput(103.65,69.31)(0.23,-0.1){2}{\line(1,0){0.23}}
\multiput(103.18,69.48)(0.47,-0.17){1}{\line(1,0){0.47}}
\multiput(102.71,69.63)(0.48,-0.15){1}{\line(1,0){0.48}}
\multiput(102.23,69.75)(0.48,-0.12){1}{\line(1,0){0.48}}
\multiput(101.74,69.85)(0.49,-0.1){1}{\line(1,0){0.49}}
\multiput(101.24,69.92)(0.49,-0.07){1}{\line(1,0){0.49}}
\multiput(100.75,69.97)(0.5,-0.05){1}{\line(1,0){0.5}}
\multiput(100.25,70)(0.5,-0.02){1}{\line(1,0){0.5}}
\put(99.75,70){\line(1,0){0.5}}
\multiput(99.25,69.97)(0.5,0.02){1}{\line(1,0){0.5}}
\multiput(98.76,69.92)(0.5,0.05){1}{\line(1,0){0.5}}
\multiput(98.26,69.85)(0.49,0.07){1}{\line(1,0){0.49}}
\multiput(97.77,69.75)(0.49,0.1){1}{\line(1,0){0.49}}
\multiput(97.29,69.63)(0.48,0.12){1}{\line(1,0){0.48}}
\multiput(96.82,69.48)(0.48,0.15){1}{\line(1,0){0.48}}
\multiput(96.35,69.31)(0.47,0.17){1}{\line(1,0){0.47}}
\multiput(95.89,69.12)(0.23,0.1){2}{\line(1,0){0.23}}
\multiput(95.44,68.9)(0.22,0.11){2}{\line(1,0){0.22}}
\multiput(95,68.66)(0.22,0.12){2}{\line(1,0){0.22}}
\multiput(94.57,68.4)(0.21,0.13){2}{\line(1,0){0.21}}
\multiput(94.16,68.12)(0.21,0.14){2}{\line(1,0){0.21}}
\multiput(93.77,67.82)(0.13,0.1){3}{\line(1,0){0.13}}
\multiput(93.38,67.5)(0.13,0.11){3}{\line(1,0){0.13}}
\multiput(93.02,67.16)(0.12,0.11){3}{\line(1,0){0.12}}
\multiput(92.67,66.8)(0.12,0.12){3}{\line(0,1){0.12}}
\multiput(92.34,66.43)(0.11,0.12){3}{\line(0,1){0.12}}
\multiput(92.03,66.04)(0.1,0.13){3}{\line(0,1){0.13}}
\multiput(91.74,65.63)(0.15,0.2){2}{\line(0,1){0.2}}
\multiput(91.47,65.21)(0.14,0.21){2}{\line(0,1){0.21}}
\multiput(91.22,64.78)(0.12,0.22){2}{\line(0,1){0.22}}
\multiput(90.99,64.34)(0.11,0.22){2}{\line(0,1){0.22}}
\multiput(90.79,63.88)(0.1,0.23){2}{\line(0,1){0.23}}
\multiput(90.6,63.42)(0.09,0.23){2}{\line(0,1){0.23}}
\multiput(90.44,62.95)(0.16,0.47){1}{\line(0,1){0.47}}
\multiput(90.31,62.47)(0.14,0.48){1}{\line(0,1){0.48}}
\multiput(90.2,61.98)(0.11,0.49){1}{\line(0,1){0.49}}
\multiput(90.11,61.49)(0.09,0.49){1}{\line(0,1){0.49}}
\multiput(90.05,61)(0.06,0.49){1}{\line(0,1){0.49}}
\multiput(90.01,60.5)(0.04,0.5){1}{\line(0,1){0.5}}
\multiput(90,60)(0.01,0.5){1}{\line(0,1){0.5}}

\linethickness{0.3mm}
\multiput(70,40)(0.12,0.12){167}{\line(1,0){0.12}}
\linethickness{0.3mm}
\put(100,70){\line(0,1){10}}
\linethickness{0.3mm}
\put(70,60){\line(0,1){20}}
\linethickness{0.3mm}
\multiput(70,60)(0.12,-0.12){67}{\line(1,0){0.12}}
\linethickness{0.3mm}
\multiput(82,48)(0.12,-0.12){67}{\line(1,0){0.12}}
\linethickness{0.3mm}
\put(110,40){\line(0,1){20}}
\linethickness{0.3mm}
\put(20,70){\line(0,1){10}}
\linethickness{0.3mm}
\put(50,40){\line(0,1){40}}
\linethickness{0.3mm}
\put(30,40){\line(0,1){20}}
\linethickness{0.3mm}
\put(10,40){\line(0,1){20}}
\put(60,60){\makebox(0,0)[cc]{$+$}}

\put(28,72){\makebox(0,0)[cc]{$\nabla_E$}}

\put(108,72){\makebox(0,0)[cc]{$\nabla_F$}}

\put(110,72){\makebox(0,0)[cc]{}}

\put(20,84){\makebox(0,0)[cc]{$E$}}

\put(50,84){\makebox(0,0)[cc]{$F$}}

\put(70,84){\makebox(0,0)[cc]{$E$}}

\put(98,84){\makebox(0,0)[cc]{$F$}}

\put(10,36){\makebox(0,0)[cc]{$\Omega^1A$}}

\put(30,36){\makebox(0,0)[cc]{$E$}}

\put(50,36){\makebox(0,0)[cc]{$F$}}

\put(70,36){\makebox(0,0)[cc]{$\Omega^1 A$}}

\put(90,36){\makebox(0,0)[cc]{$E$}}

\put(110,36){\makebox(0,0)[cc]{$F$}}

\put(88,50){\makebox(0,0)[cc]{$\sigma_E$}}
%%%%%%%%%%%%%%%
\refstepcounter{piccie} \label{someLabel2}
\put(140,36){\makebox(0,0)[cc]{Fig.\ \arabic{piccie}}}
%%%%%%%%%%%%%%%
\end{picture}

\noindent The reader should note that neither term in Fig.\ \ref{someLabel2} in isolation is well defined on $E\tens_A F$, only the sum is. It will be a general principle of our formulae or diagrams that care must be taken with operators which are not module maps. It will be convenient to think of the operators as defined on $\tens_{\mathbb{C}}$ rather than on $\tens_A$, with a separate check that the required sum of operations is well defined on the given domain. 
Similarly, to illustrate the formula for $\sigma_{E\tens F}$ in Proposition~\ref{tens45} we use Fig.\ \ref{someLabel3}:

\unitlength 0.5 mm
\begin{picture}(50,44)(-20,30)
\linethickness{0.3mm}
\multiput(10,40)(0.24,0.12){167}{\line(1,0){0.24}}
\linethickness{0.3mm}
\multiput(10,60)(0.12,-0.12){83}{\line(1,0){0.12}}
\linethickness{0.3mm}
\multiput(40,50)(0.12,-0.12){83}{\line(1,0){0.12}}
\linethickness{0.3mm}
\multiput(26,44)(0.12,-0.12){33}{\line(1,0){0.12}}
\linethickness{0.3mm}
\multiput(30,60)(0.12,-0.12){33}{\line(1,0){0.12}}
\put(10,65){\makebox(0,0)[cc]{$E$}}

\put(32,65){\makebox(0,0)[cc]{$F$}}

\put(48,65){\makebox(0,0)[cc]{$\Omega^1 A$}}

\put(8,36){\makebox(0,0)[cc]{$\Omega^1 A$}}

\put(30,35){\makebox(0,0)[cc]{$E$}}

\put(50,35){\makebox(0,0)[cc]{$F$}}
%%%%%%%%%%%%%%%
\refstepcounter{piccie} \label{someLabel3}
\put(97,35){\makebox(0,0)[cc]{Fig.\ \arabic{piccie}}}
%%%%%%%%%%%%%%%

\put(12,48){\makebox(0,0)[cc]{$\sigma_E$}}

\put(46,52){\makebox(0,0)[cc]{$\sigma_F$}}

\end{picture}

\noindent Note that we are not implying that we have any form of braid relation by using this notation. However we must comment on notation in the paper. From hard experience of the confusion caused by trying it the other way, we shall shall describe products of various $\sigma$ and $\sigma^{-1}$ as braids (we continue to use this term even in the absence of any braid relation, and may use the term `generalised braids'), and frequently not give any explicit dependence on domain. The domains can easily be established by a study of the formula or diagram. 

\begin{defin}    \label{cvgjfxtzhc}
The category ${}_A\mathcal{E}$ consists of objects $(E,\nabla_E)$, where
$E$ is a left  $A$-module, and $\nabla_E$ is a left covariant derivative on $E$. The morphisms $T:(E,\nabla_E)\to (F,\nabla_F)$ consist of left module maps
$T:E\to F$ for which $(\id\tens T)\nabla_E=\nabla_F\,T:E\to \Omega^1 A\tens_A F$.

The category ${}_A\mathcal{E}_A$ consists of objects $(E,\nabla_E,\sigma_E)$, where
$E$ is an  $A$-bimodule, and $(\nabla_E,\sigma_E)$ is a bimodule covariant derivative on $E$,
and where $\sigma_E:E\tens_A\Omega^1 A\to \Omega^1 A \tens_A E$ is invertible. The morphisms $T:(E,\nabla_E,\sigma_E)\to (F,\nabla_F,\sigma_F)$ consist of bimodule maps
$T:E\to F$ for which $(\id\tens T)\nabla_E=\nabla_F\,T:E\to \Omega^1 A\tens_A F$. It is then automatically true that $\sigma_F(T\tens\id)=(\id\tens T)\sigma_E$. Taking the identity object as 
the algebra $A$ itself, with $\nabla=\extd:A\to \Omega^1 A\tens_A A\cong\Omega^1 A$ and tensor product as in Proposition \ref{tens45}, makes ${}_A\mathcal{E}_A$ into a monoidal category. 
\end{defin}

\subsection{Vector fields and $n$-tuples of vector fields}
We now come to a controversial question: what is a noncommutative vector field on an algebra $A$? We may as well be honest and say that we shall not consider derivations on an algebra, though it is possible to use a generalised idea of derivation, see \cite{JarLle,Borowiec}, and in certain circumstances the idea of braided derivation would apply, see \cite{Mabicdifcalc,LychBraid}. 
We shall take the approach that the vector fields $\Vec A$ are the dual of the 1-forms $\Omega^1 A$ (see\cite{BegNCvec}). 

With quite weak conditions on the topology, the space of sections of a locally trivial finite dimensional vector bundle on 
a finite dimensional space $X$ obeys the finitely generated projective property as a module over $C(X)$, the continuous functions on $X$. This is part of the content of the Serre-Swan theorem, and is the foundation of $K$-theory. From this, it is reasonable to suggest that the module of 1-forms on a  `finite dimensional noncommutative manifold', whatever that is, should be finitely generated projective. As we have seen in Subsection~\ref{vcxgaksujhgc}, there is then a nice idea of the dual, and we define, as in the classical case, the vector fields $\Vec A$ as the dual of the 1-forms $\Omega^1 A$. As we are in the noncommutative case, we need to choose a side, so we choose $\Omega^1 A$ to be finitely generated projective as a right module, and then we take $\Vec A$ to be the
right module maps from $\Omega^1 A$ to $A$. 
There is of course a certain symmetry here, we could have defined $\Omega^1 A$ as the left dual of $\Vec A$. However as $A$ is considered to be an algebra, the idea of a differential calculus as a differential graded algebra extending $A$ has proven to be more popular, and thus we take $\Omega^1 A$ as the fundamental object. 

We suppose that $\Omega^1 A$ is finitely generated projective as a right $A$ module, and set $\Vec A=\Hom_A(\Omega^1 A,A)$. We denote the evaluation and coevaluation maps by
\begin{eqnarray*}
\ev:\Vec A\tens_A \Omega^1 A\to A\ ,\quad \coev:A\to \Omega^1 A\tens_A \Vec A\ .
\end{eqnarray*}
Copying from Fig.\ \ref{someLabel}, the relation between the evaluation and coevaluation maps is given by Fig.\ \ref{dshbsh}.

\unitlength 0.5 mm
\begin{picture}(220,65)(-20,30)
\linethickness{0.3mm}
\multiput(109.99,60.5)(0.01,-0.5){1}{\line(0,-1){0.5}}
\multiput(109.95,61)(0.04,-0.5){1}{\line(0,-1){0.5}}
\multiput(109.89,61.49)(0.06,-0.49){1}{\line(0,-1){0.49}}
\multiput(109.8,61.98)(0.09,-0.49){1}{\line(0,-1){0.49}}
\multiput(109.69,62.47)(0.11,-0.49){1}{\line(0,-1){0.49}}
\multiput(109.56,62.95)(0.14,-0.48){1}{\line(0,-1){0.48}}
\multiput(109.4,63.42)(0.16,-0.47){1}{\line(0,-1){0.47}}
\multiput(109.21,63.88)(0.09,-0.23){2}{\line(0,-1){0.23}}
\multiput(109.01,64.34)(0.1,-0.23){2}{\line(0,-1){0.23}}
\multiput(108.78,64.78)(0.11,-0.22){2}{\line(0,-1){0.22}}
\multiput(108.53,65.21)(0.12,-0.22){2}{\line(0,-1){0.22}}
\multiput(108.26,65.63)(0.14,-0.21){2}{\line(0,-1){0.21}}
\multiput(107.97,66.04)(0.15,-0.2){2}{\line(0,-1){0.2}}
\multiput(107.66,66.43)(0.1,-0.13){3}{\line(0,-1){0.13}}
\multiput(107.33,66.8)(0.11,-0.12){3}{\line(0,-1){0.12}}
\multiput(106.98,67.16)(0.12,-0.12){3}{\line(0,-1){0.12}}
\multiput(106.62,67.5)(0.12,-0.11){3}{\line(1,0){0.12}}
\multiput(106.23,67.82)(0.13,-0.11){3}{\line(1,0){0.13}}
\multiput(105.84,68.12)(0.13,-0.1){3}{\line(1,0){0.13}}
\multiput(105.43,68.4)(0.21,-0.14){2}{\line(1,0){0.21}}
\multiput(105,68.66)(0.21,-0.13){2}{\line(1,0){0.21}}
\multiput(104.56,68.9)(0.22,-0.12){2}{\line(1,0){0.22}}
\multiput(104.11,69.12)(0.22,-0.11){2}{\line(1,0){0.22}}
\multiput(103.65,69.31)(0.23,-0.1){2}{\line(1,0){0.23}}
\multiput(103.18,69.48)(0.47,-0.17){1}{\line(1,0){0.47}}
\multiput(102.71,69.63)(0.48,-0.15){1}{\line(1,0){0.48}}
\multiput(102.23,69.75)(0.48,-0.12){1}{\line(1,0){0.48}}
\multiput(101.74,69.85)(0.49,-0.1){1}{\line(1,0){0.49}}
\multiput(101.24,69.92)(0.49,-0.07){1}{\line(1,0){0.49}}
\multiput(100.75,69.97)(0.5,-0.05){1}{\line(1,0){0.5}}
\multiput(100.25,70)(0.5,-0.02){1}{\line(1,0){0.5}}
\put(99.75,70){\line(1,0){0.5}}
\multiput(99.25,69.97)(0.5,0.02){1}{\line(1,0){0.5}}
\multiput(98.76,69.92)(0.5,0.05){1}{\line(1,0){0.5}}
\multiput(98.26,69.85)(0.49,0.07){1}{\line(1,0){0.49}}
\multiput(97.77,69.75)(0.49,0.1){1}{\line(1,0){0.49}}
\multiput(97.29,69.63)(0.48,0.12){1}{\line(1,0){0.48}}
\multiput(96.82,69.48)(0.48,0.15){1}{\line(1,0){0.48}}
\multiput(96.35,69.31)(0.47,0.17){1}{\line(1,0){0.47}}
\multiput(95.89,69.12)(0.23,0.1){2}{\line(1,0){0.23}}
\multiput(95.44,68.9)(0.22,0.11){2}{\line(1,0){0.22}}
\multiput(95,68.66)(0.22,0.12){2}{\line(1,0){0.22}}
\multiput(94.57,68.4)(0.21,0.13){2}{\line(1,0){0.21}}
\multiput(94.16,68.12)(0.21,0.14){2}{\line(1,0){0.21}}
\multiput(93.77,67.82)(0.13,0.1){3}{\line(1,0){0.13}}
\multiput(93.38,67.5)(0.13,0.11){3}{\line(1,0){0.13}}
\multiput(93.02,67.16)(0.12,0.11){3}{\line(1,0){0.12}}
\multiput(92.67,66.8)(0.12,0.12){3}{\line(0,1){0.12}}
\multiput(92.34,66.43)(0.11,0.12){3}{\line(0,1){0.12}}
\multiput(92.03,66.04)(0.1,0.13){3}{\line(0,1){0.13}}
\multiput(91.74,65.63)(0.15,0.2){2}{\line(0,1){0.2}}
\multiput(91.47,65.21)(0.14,0.21){2}{\line(0,1){0.21}}
\multiput(91.22,64.78)(0.12,0.22){2}{\line(0,1){0.22}}
\multiput(90.99,64.34)(0.11,0.22){2}{\line(0,1){0.22}}
\multiput(90.79,63.88)(0.1,0.23){2}{\line(0,1){0.23}}
\multiput(90.6,63.42)(0.09,0.23){2}{\line(0,1){0.23}}
\multiput(90.44,62.95)(0.16,0.47){1}{\line(0,1){0.47}}
\multiput(90.31,62.47)(0.14,0.48){1}{\line(0,1){0.48}}
\multiput(90.2,61.98)(0.11,0.49){1}{\line(0,1){0.49}}
\multiput(90.11,61.49)(0.09,0.49){1}{\line(0,1){0.49}}
\multiput(90.05,61)(0.06,0.49){1}{\line(0,1){0.49}}
\multiput(90.01,60.5)(0.04,0.5){1}{\line(0,1){0.5}}
\multiput(90,60)(0.01,0.5){1}{\line(0,1){0.5}}

\linethickness{0.3mm}
\multiput(110,60)(0.01,-0.5){1}{\line(0,-1){0.5}}
\multiput(110.01,59.5)(0.04,-0.5){1}{\line(0,-1){0.5}}
\multiput(110.05,59)(0.06,-0.49){1}{\line(0,-1){0.49}}
\multiput(110.11,58.51)(0.09,-0.49){1}{\line(0,-1){0.49}}
\multiput(110.2,58.02)(0.11,-0.49){1}{\line(0,-1){0.49}}
\multiput(110.31,57.53)(0.14,-0.48){1}{\line(0,-1){0.48}}
\multiput(110.44,57.05)(0.16,-0.47){1}{\line(0,-1){0.47}}
\multiput(110.6,56.58)(0.09,-0.23){2}{\line(0,-1){0.23}}
\multiput(110.79,56.12)(0.1,-0.23){2}{\line(0,-1){0.23}}
\multiput(110.99,55.66)(0.11,-0.22){2}{\line(0,-1){0.22}}
\multiput(111.22,55.22)(0.12,-0.22){2}{\line(0,-1){0.22}}
\multiput(111.47,54.79)(0.14,-0.21){2}{\line(0,-1){0.21}}
\multiput(111.74,54.37)(0.15,-0.2){2}{\line(0,-1){0.2}}
\multiput(112.03,53.96)(0.1,-0.13){3}{\line(0,-1){0.13}}
\multiput(112.34,53.57)(0.11,-0.12){3}{\line(0,-1){0.12}}
\multiput(112.67,53.2)(0.12,-0.12){3}{\line(0,-1){0.12}}
\multiput(113.02,52.84)(0.12,-0.11){3}{\line(1,0){0.12}}
\multiput(113.38,52.5)(0.13,-0.11){3}{\line(1,0){0.13}}
\multiput(113.77,52.18)(0.13,-0.1){3}{\line(1,0){0.13}}
\multiput(114.16,51.88)(0.21,-0.14){2}{\line(1,0){0.21}}
\multiput(114.57,51.6)(0.21,-0.13){2}{\line(1,0){0.21}}
\multiput(115,51.34)(0.22,-0.12){2}{\line(1,0){0.22}}
\multiput(115.44,51.1)(0.22,-0.11){2}{\line(1,0){0.22}}
\multiput(115.89,50.88)(0.23,-0.1){2}{\line(1,0){0.23}}
\multiput(116.35,50.69)(0.47,-0.17){1}{\line(1,0){0.47}}
\multiput(116.82,50.52)(0.48,-0.15){1}{\line(1,0){0.48}}
\multiput(117.29,50.37)(0.48,-0.12){1}{\line(1,0){0.48}}
\multiput(117.77,50.25)(0.49,-0.1){1}{\line(1,0){0.49}}
\multiput(118.26,50.15)(0.49,-0.07){1}{\line(1,0){0.49}}
\multiput(118.76,50.08)(0.5,-0.05){1}{\line(1,0){0.5}}
\multiput(119.25,50.03)(0.5,-0.02){1}{\line(1,0){0.5}}
\put(119.75,50){\line(1,0){0.5}}
\multiput(120.25,50)(0.5,0.02){1}{\line(1,0){0.5}}
\multiput(120.75,50.03)(0.5,0.05){1}{\line(1,0){0.5}}
\multiput(121.24,50.08)(0.49,0.07){1}{\line(1,0){0.49}}
\multiput(121.74,50.15)(0.49,0.1){1}{\line(1,0){0.49}}
\multiput(122.23,50.25)(0.48,0.12){1}{\line(1,0){0.48}}
\multiput(122.71,50.37)(0.48,0.15){1}{\line(1,0){0.48}}
\multiput(123.18,50.52)(0.47,0.17){1}{\line(1,0){0.47}}
\multiput(123.65,50.69)(0.23,0.1){2}{\line(1,0){0.23}}
\multiput(124.11,50.88)(0.22,0.11){2}{\line(1,0){0.22}}
\multiput(124.56,51.1)(0.22,0.12){2}{\line(1,0){0.22}}
\multiput(125,51.34)(0.21,0.13){2}{\line(1,0){0.21}}
\multiput(125.43,51.6)(0.21,0.14){2}{\line(1,0){0.21}}
\multiput(125.84,51.88)(0.13,0.1){3}{\line(1,0){0.13}}
\multiput(126.23,52.18)(0.13,0.11){3}{\line(1,0){0.13}}
\multiput(126.62,52.5)(0.12,0.11){3}{\line(1,0){0.12}}
\multiput(126.98,52.84)(0.12,0.12){3}{\line(0,1){0.12}}
\multiput(127.33,53.2)(0.11,0.12){3}{\line(0,1){0.12}}
\multiput(127.66,53.57)(0.1,0.13){3}{\line(0,1){0.13}}
\multiput(127.97,53.96)(0.15,0.2){2}{\line(0,1){0.2}}
\multiput(128.26,54.37)(0.14,0.21){2}{\line(0,1){0.21}}
\multiput(128.53,54.79)(0.12,0.22){2}{\line(0,1){0.22}}
\multiput(128.78,55.22)(0.11,0.22){2}{\line(0,1){0.22}}
\multiput(129.01,55.66)(0.1,0.23){2}{\line(0,1){0.23}}
\multiput(129.21,56.12)(0.09,0.23){2}{\line(0,1){0.23}}
\multiput(129.4,56.58)(0.16,0.47){1}{\line(0,1){0.47}}
\multiput(129.56,57.05)(0.14,0.48){1}{\line(0,1){0.48}}
\multiput(129.69,57.53)(0.11,0.49){1}{\line(0,1){0.49}}
\multiput(129.8,58.02)(0.09,0.49){1}{\line(0,1){0.49}}
\multiput(129.89,58.51)(0.06,0.49){1}{\line(0,1){0.49}}
\multiput(129.95,59)(0.04,0.5){1}{\line(0,1){0.5}}
\multiput(129.99,59.5)(0.01,0.5){1}{\line(0,1){0.5}}

\linethickness{0.3mm}
\put(90,40){\line(0,1){20}}
\linethickness{0.3mm}
\put(130,60){\line(0,1){20}}
\linethickness{0.3mm}
\put(150,40){\line(0,1){40}}
\linethickness{0.3mm}
\multiput(39.99,60.5)(0.01,-0.5){1}{\line(0,-1){0.5}}
\multiput(39.95,61)(0.04,-0.5){1}{\line(0,-1){0.5}}
\multiput(39.89,61.49)(0.06,-0.49){1}{\line(0,-1){0.49}}
\multiput(39.8,61.98)(0.09,-0.49){1}{\line(0,-1){0.49}}
\multiput(39.69,62.47)(0.11,-0.49){1}{\line(0,-1){0.49}}
\multiput(39.56,62.95)(0.14,-0.48){1}{\line(0,-1){0.48}}
\multiput(39.4,63.42)(0.16,-0.47){1}{\line(0,-1){0.47}}
\multiput(39.21,63.88)(0.09,-0.23){2}{\line(0,-1){0.23}}
\multiput(39.01,64.34)(0.1,-0.23){2}{\line(0,-1){0.23}}
\multiput(38.78,64.78)(0.11,-0.22){2}{\line(0,-1){0.22}}
\multiput(38.53,65.21)(0.12,-0.22){2}{\line(0,-1){0.22}}
\multiput(38.26,65.63)(0.14,-0.21){2}{\line(0,-1){0.21}}
\multiput(37.97,66.04)(0.15,-0.2){2}{\line(0,-1){0.2}}
\multiput(37.66,66.43)(0.1,-0.13){3}{\line(0,-1){0.13}}
\multiput(37.33,66.8)(0.11,-0.12){3}{\line(0,-1){0.12}}
\multiput(36.98,67.16)(0.12,-0.12){3}{\line(0,-1){0.12}}
\multiput(36.62,67.5)(0.12,-0.11){3}{\line(1,0){0.12}}
\multiput(36.23,67.82)(0.13,-0.11){3}{\line(1,0){0.13}}
\multiput(35.84,68.12)(0.13,-0.1){3}{\line(1,0){0.13}}
\multiput(35.43,68.4)(0.21,-0.14){2}{\line(1,0){0.21}}
\multiput(35,68.66)(0.21,-0.13){2}{\line(1,0){0.21}}
\multiput(34.56,68.9)(0.22,-0.12){2}{\line(1,0){0.22}}
\multiput(34.11,69.12)(0.22,-0.11){2}{\line(1,0){0.22}}
\multiput(33.65,69.31)(0.23,-0.1){2}{\line(1,0){0.23}}
\multiput(33.18,69.48)(0.47,-0.17){1}{\line(1,0){0.47}}
\multiput(32.71,69.63)(0.48,-0.15){1}{\line(1,0){0.48}}
\multiput(32.23,69.75)(0.48,-0.12){1}{\line(1,0){0.48}}
\multiput(31.74,69.85)(0.49,-0.1){1}{\line(1,0){0.49}}
\multiput(31.24,69.92)(0.49,-0.07){1}{\line(1,0){0.49}}
\multiput(30.75,69.97)(0.5,-0.05){1}{\line(1,0){0.5}}
\multiput(30.25,70)(0.5,-0.02){1}{\line(1,0){0.5}}
\put(29.75,70){\line(1,0){0.5}}
\multiput(29.25,69.97)(0.5,0.02){1}{\line(1,0){0.5}}
\multiput(28.76,69.92)(0.5,0.05){1}{\line(1,0){0.5}}
\multiput(28.26,69.85)(0.49,0.07){1}{\line(1,0){0.49}}
\multiput(27.77,69.75)(0.49,0.1){1}{\line(1,0){0.49}}
\multiput(27.29,69.63)(0.48,0.12){1}{\line(1,0){0.48}}
\multiput(26.82,69.48)(0.48,0.15){1}{\line(1,0){0.48}}
\multiput(26.35,69.31)(0.47,0.17){1}{\line(1,0){0.47}}
\multiput(25.89,69.12)(0.23,0.1){2}{\line(1,0){0.23}}
\multiput(25.44,68.9)(0.22,0.11){2}{\line(1,0){0.22}}
\multiput(25,68.66)(0.22,0.12){2}{\line(1,0){0.22}}
\multiput(24.57,68.4)(0.21,0.13){2}{\line(1,0){0.21}}
\multiput(24.16,68.12)(0.21,0.14){2}{\line(1,0){0.21}}
\multiput(23.77,67.82)(0.13,0.1){3}{\line(1,0){0.13}}
\multiput(23.38,67.5)(0.13,0.11){3}{\line(1,0){0.13}}
\multiput(23.02,67.16)(0.12,0.11){3}{\line(1,0){0.12}}
\multiput(22.67,66.8)(0.12,0.12){3}{\line(0,1){0.12}}
\multiput(22.34,66.43)(0.11,0.12){3}{\line(0,1){0.12}}
\multiput(22.03,66.04)(0.1,0.13){3}{\line(0,1){0.13}}
\multiput(21.74,65.63)(0.15,0.2){2}{\line(0,1){0.2}}
\multiput(21.47,65.21)(0.14,0.21){2}{\line(0,1){0.21}}
\multiput(21.22,64.78)(0.12,0.22){2}{\line(0,1){0.22}}
\multiput(20.99,64.34)(0.11,0.22){2}{\line(0,1){0.22}}
\multiput(20.79,63.88)(0.1,0.23){2}{\line(0,1){0.23}}
\multiput(20.6,63.42)(0.09,0.23){2}{\line(0,1){0.23}}
\multiput(20.44,62.95)(0.16,0.47){1}{\line(0,1){0.47}}
\multiput(20.31,62.47)(0.14,0.48){1}{\line(0,1){0.48}}
\multiput(20.2,61.98)(0.11,0.49){1}{\line(0,1){0.49}}
\multiput(20.11,61.49)(0.09,0.49){1}{\line(0,1){0.49}}
\multiput(20.05,61)(0.06,0.49){1}{\line(0,1){0.49}}
\multiput(20.01,60.5)(0.04,0.5){1}{\line(0,1){0.5}}
\multiput(20,60)(0.01,0.5){1}{\line(0,1){0.5}}

\linethickness{0.3mm}
\multiput(0,60)(0.01,-0.5){1}{\line(0,-1){0.5}}
\multiput(0.01,59.5)(0.04,-0.5){1}{\line(0,-1){0.5}}
\multiput(0.05,59)(0.06,-0.49){1}{\line(0,-1){0.49}}
\multiput(0.11,58.51)(0.09,-0.49){1}{\line(0,-1){0.49}}
\multiput(0.2,58.02)(0.11,-0.49){1}{\line(0,-1){0.49}}
\multiput(0.31,57.53)(0.14,-0.48){1}{\line(0,-1){0.48}}
\multiput(0.44,57.05)(0.16,-0.47){1}{\line(0,-1){0.47}}
\multiput(0.6,56.58)(0.09,-0.23){2}{\line(0,-1){0.23}}
\multiput(0.79,56.12)(0.1,-0.23){2}{\line(0,-1){0.23}}
\multiput(0.99,55.66)(0.11,-0.22){2}{\line(0,-1){0.22}}
\multiput(1.22,55.22)(0.12,-0.22){2}{\line(0,-1){0.22}}
\multiput(1.47,54.79)(0.14,-0.21){2}{\line(0,-1){0.21}}
\multiput(1.74,54.37)(0.15,-0.2){2}{\line(0,-1){0.2}}
\multiput(2.03,53.96)(0.1,-0.13){3}{\line(0,-1){0.13}}
\multiput(2.34,53.57)(0.11,-0.12){3}{\line(0,-1){0.12}}
\multiput(2.67,53.2)(0.12,-0.12){3}{\line(0,-1){0.12}}
\multiput(3.02,52.84)(0.12,-0.11){3}{\line(1,0){0.12}}
\multiput(3.38,52.5)(0.13,-0.11){3}{\line(1,0){0.13}}
\multiput(3.77,52.18)(0.13,-0.1){3}{\line(1,0){0.13}}
\multiput(4.16,51.88)(0.21,-0.14){2}{\line(1,0){0.21}}
\multiput(4.57,51.6)(0.21,-0.13){2}{\line(1,0){0.21}}
\multiput(5,51.34)(0.22,-0.12){2}{\line(1,0){0.22}}
\multiput(5.44,51.1)(0.22,-0.11){2}{\line(1,0){0.22}}
\multiput(5.89,50.88)(0.23,-0.1){2}{\line(1,0){0.23}}
\multiput(6.35,50.69)(0.47,-0.17){1}{\line(1,0){0.47}}
\multiput(6.82,50.52)(0.48,-0.15){1}{\line(1,0){0.48}}
\multiput(7.29,50.37)(0.48,-0.12){1}{\line(1,0){0.48}}
\multiput(7.77,50.25)(0.49,-0.1){1}{\line(1,0){0.49}}
\multiput(8.26,50.15)(0.49,-0.07){1}{\line(1,0){0.49}}
\multiput(8.76,50.08)(0.5,-0.05){1}{\line(1,0){0.5}}
\multiput(9.25,50.03)(0.5,-0.02){1}{\line(1,0){0.5}}
\put(9.75,50){\line(1,0){0.5}}
\multiput(10.25,50)(0.5,0.02){1}{\line(1,0){0.5}}
\multiput(10.75,50.03)(0.5,0.05){1}{\line(1,0){0.5}}
\multiput(11.24,50.08)(0.49,0.07){1}{\line(1,0){0.49}}
\multiput(11.74,50.15)(0.49,0.1){1}{\line(1,0){0.49}}
\multiput(12.23,50.25)(0.48,0.12){1}{\line(1,0){0.48}}
\multiput(12.71,50.37)(0.48,0.15){1}{\line(1,0){0.48}}
\multiput(13.18,50.52)(0.47,0.17){1}{\line(1,0){0.47}}
\multiput(13.65,50.69)(0.23,0.1){2}{\line(1,0){0.23}}
\multiput(14.11,50.88)(0.22,0.11){2}{\line(1,0){0.22}}
\multiput(14.56,51.1)(0.22,0.12){2}{\line(1,0){0.22}}
\multiput(15,51.34)(0.21,0.13){2}{\line(1,0){0.21}}
\multiput(15.43,51.6)(0.21,0.14){2}{\line(1,0){0.21}}
\multiput(15.84,51.88)(0.13,0.1){3}{\line(1,0){0.13}}
\multiput(16.23,52.18)(0.13,0.11){3}{\line(1,0){0.13}}
\multiput(16.62,52.5)(0.12,0.11){3}{\line(1,0){0.12}}
\multiput(16.98,52.84)(0.12,0.12){3}{\line(0,1){0.12}}
\multiput(17.33,53.2)(0.11,0.12){3}{\line(0,1){0.12}}
\multiput(17.66,53.57)(0.1,0.13){3}{\line(0,1){0.13}}
\multiput(17.97,53.96)(0.15,0.2){2}{\line(0,1){0.2}}
\multiput(18.26,54.37)(0.14,0.21){2}{\line(0,1){0.21}}
\multiput(18.53,54.79)(0.12,0.22){2}{\line(0,1){0.22}}
\multiput(18.78,55.22)(0.11,0.22){2}{\line(0,1){0.22}}
\multiput(19.01,55.66)(0.1,0.23){2}{\line(0,1){0.23}}
\multiput(19.21,56.12)(0.09,0.23){2}{\line(0,1){0.23}}
\multiput(19.4,56.58)(0.16,0.47){1}{\line(0,1){0.47}}
\multiput(19.56,57.05)(0.14,0.48){1}{\line(0,1){0.48}}
\multiput(19.69,57.53)(0.11,0.49){1}{\line(0,1){0.49}}
\multiput(19.8,58.02)(0.09,0.49){1}{\line(0,1){0.49}}
\multiput(19.89,58.51)(0.06,0.49){1}{\line(0,1){0.49}}
\multiput(19.95,59)(0.04,0.5){1}{\line(0,1){0.5}}
\multiput(19.99,59.5)(0.01,0.5){1}{\line(0,1){0.5}}

\linethickness{0.3mm}
\put(40,40){\line(0,1){20}}
\linethickness{0.3mm}
\put(0,60){\line(0,1){20}}
\put(50,60){\makebox(0,0)[cc]{$=$}}

\linethickness{0.3mm}
\put(60,40){\line(0,1){40}}
\put(140,60){\makebox(0,0)[cc]{$=$}}

\put(25,75){\makebox(0,0)[cc]{$\coev$}}

\put(25,75){\makebox(0,0)[cc]{}}

\put(95,75){\makebox(0,0)[cc]{$\coev$}}

\put(115,45){\makebox(0,0)[cc]{$\ev$}}

\put(5,45){\makebox(0,0)[cc]{$\ev$}}

\put(0,85){\makebox(0,0)[cc]{$\Vec A$}}

\put(55,85){\makebox(0,0)[cc]{$\Vec A$}}

\put(35,35){\makebox(0,0)[cc]{$\Vec A$}}

\put(55,35){\makebox(0,0)[cc]{$\Vec A$}}

\put(125,85){\makebox(0,0)[cc]{$\Omega^1 A$}}

\put(145,85){\makebox(0,0)[cc]{$\Omega^1 A$}}

\put(85,35){\makebox(0,0)[cc]{$\Omega^1 A$}}

\put(145,35){\makebox(0,0)[cc]{$\Omega^1 A$}}

%%%%%%%%%%%%%%%
\refstepcounter{piccie} \label{dshbsh}
\put(185,55){\makebox(0,0)[cc]{Fig.\ \arabic{piccie}}}
%%%%%%%%%%%%%%%

\end{picture}

\noindent We shall require multiple copies of 1-forms and vector fields, so we define
\begin{eqnarray*}
&& \Vec^{\tens 0}A\ =\ A\ ,\quad
\Vec^{\tens n}A\ =\ \Vec A\tens_A \Vec A\tens_A\dots \tens_A \Vec A\ ,\cr
&& \Omega^{\tens 0}A\ =\ A\ ,\quad
\Omega^{\tens n}A\ =\ \Omega^1 A\tens_A \Omega^1 A\tens_A\dots \tens_A \Omega^1 A\ ,
\end{eqnarray*}
where we have $n$ copies of $\Vec A$ and $\Omega^1 A$. It is important to note that the definition of 
$\Vec^{\tens n}A$ and $\Omega^{\tens n}A$ uses $\tens_A$, the tensor product over the algebra.
For clarity, we will often use $\id^{\tens n}$ as the identity on $\Vec^{\tens n}A$
or $\Omega^{\tens n}A$.

\begin{defin}\label{bbieob}
Define the $n$-fold evaluation map $\ev^{\<n\>}:\Vec^{\tens n}A\tens_A \Omega^{\tens n}A\to A$
recursively by
\begin{eqnarray*}
\ev^{\<1\>}\ =\ \ev\ ,\quad \ev^{\<n+1\>}\ =\ \ev\,(\id\tens \ev^{\<n\>}  \tens \id)\ ,
\end{eqnarray*}
and the coevaluation map $\coev^{\<n\>}:A\to \Omega^{\tens n}A \tens_A \Vec^{\tens n}A$  by
\begin{eqnarray*}
\coev^{\<1\>}\ =\ \coev\ ,\quad \coev^{\<n+1\>}\ =\ (\id\tens \coev^{\<n\>}  \tens \id)\,\coev\ ,
\end{eqnarray*}
\end{defin}

To illustrate Definition~\ref{bbieob} we use Fig.\ \ref{dshbshdxh}, the coevaluation diagram looks similar but upside down!

\unitlength 0.5 mm
\begin{picture}(130,60)(-20,45)
\put(55,70){\makebox(0,0)[cc]{$=$}}

\put(0,85){\makebox(0,0)[cc]{$\Vec^{\tens n+1}A$}}

\put(40,85){\makebox(0,0)[cc]{$\Omega^{\tens n+1}A$}}

%\put(130,85){\makebox(0,0)[cc]{}}

\linethickness{0.3mm}
\put(90,70){\line(0,1){10}}
\linethickness{0.3mm}
\put(110,70){\line(0,1){10}}
\linethickness{0.3mm}
\multiput(90,70)(0.01,-0.5){1}{\line(0,-1){0.5}}
\multiput(90.01,69.5)(0.04,-0.5){1}{\line(0,-1){0.5}}
\multiput(90.05,69)(0.06,-0.49){1}{\line(0,-1){0.49}}
\multiput(90.11,68.51)(0.09,-0.49){1}{\line(0,-1){0.49}}
\multiput(90.2,68.02)(0.11,-0.49){1}{\line(0,-1){0.49}}
\multiput(90.31,67.53)(0.14,-0.48){1}{\line(0,-1){0.48}}
\multiput(90.44,67.05)(0.16,-0.47){1}{\line(0,-1){0.47}}
\multiput(90.6,66.58)(0.09,-0.23){2}{\line(0,-1){0.23}}
\multiput(90.79,66.12)(0.1,-0.23){2}{\line(0,-1){0.23}}
\multiput(90.99,65.66)(0.11,-0.22){2}{\line(0,-1){0.22}}
\multiput(91.22,65.22)(0.12,-0.22){2}{\line(0,-1){0.22}}
\multiput(91.47,64.79)(0.14,-0.21){2}{\line(0,-1){0.21}}
\multiput(91.74,64.37)(0.15,-0.2){2}{\line(0,-1){0.2}}
\multiput(92.03,63.96)(0.1,-0.13){3}{\line(0,-1){0.13}}
\multiput(92.34,63.57)(0.11,-0.12){3}{\line(0,-1){0.12}}
\multiput(92.67,63.2)(0.12,-0.12){3}{\line(0,-1){0.12}}
\multiput(93.02,62.84)(0.12,-0.11){3}{\line(1,0){0.12}}
\multiput(93.38,62.5)(0.13,-0.11){3}{\line(1,0){0.13}}
\multiput(93.77,62.18)(0.13,-0.1){3}{\line(1,0){0.13}}
\multiput(94.16,61.88)(0.21,-0.14){2}{\line(1,0){0.21}}
\multiput(94.57,61.6)(0.21,-0.13){2}{\line(1,0){0.21}}
\multiput(95,61.34)(0.22,-0.12){2}{\line(1,0){0.22}}
\multiput(95.44,61.1)(0.22,-0.11){2}{\line(1,0){0.22}}
\multiput(95.89,60.88)(0.23,-0.1){2}{\line(1,0){0.23}}
\multiput(96.35,60.69)(0.47,-0.17){1}{\line(1,0){0.47}}
\multiput(96.82,60.52)(0.48,-0.15){1}{\line(1,0){0.48}}
\multiput(97.29,60.37)(0.48,-0.12){1}{\line(1,0){0.48}}
\multiput(97.77,60.25)(0.49,-0.1){1}{\line(1,0){0.49}}
\multiput(98.26,60.15)(0.49,-0.07){1}{\line(1,0){0.49}}
\multiput(98.76,60.08)(0.5,-0.05){1}{\line(1,0){0.5}}
\multiput(99.25,60.03)(0.5,-0.02){1}{\line(1,0){0.5}}
\put(99.75,60){\line(1,0){0.5}}
\multiput(100.25,60)(0.5,0.02){1}{\line(1,0){0.5}}
\multiput(100.75,60.03)(0.5,0.05){1}{\line(1,0){0.5}}
\multiput(101.24,60.08)(0.49,0.07){1}{\line(1,0){0.49}}
\multiput(101.74,60.15)(0.49,0.1){1}{\line(1,0){0.49}}
\multiput(102.23,60.25)(0.48,0.12){1}{\line(1,0){0.48}}
\multiput(102.71,60.37)(0.48,0.15){1}{\line(1,0){0.48}}
\multiput(103.18,60.52)(0.47,0.17){1}{\line(1,0){0.47}}
\multiput(103.65,60.69)(0.23,0.1){2}{\line(1,0){0.23}}
\multiput(104.11,60.88)(0.22,0.11){2}{\line(1,0){0.22}}
\multiput(104.56,61.1)(0.22,0.12){2}{\line(1,0){0.22}}
\multiput(105,61.34)(0.21,0.13){2}{\line(1,0){0.21}}
\multiput(105.43,61.6)(0.21,0.14){2}{\line(1,0){0.21}}
\multiput(105.84,61.88)(0.13,0.1){3}{\line(1,0){0.13}}
\multiput(106.23,62.18)(0.13,0.11){3}{\line(1,0){0.13}}
\multiput(106.62,62.5)(0.12,0.11){3}{\line(1,0){0.12}}
\multiput(106.98,62.84)(0.12,0.12){3}{\line(0,1){0.12}}
\multiput(107.33,63.2)(0.11,0.12){3}{\line(0,1){0.12}}
\multiput(107.66,63.57)(0.1,0.13){3}{\line(0,1){0.13}}
\multiput(107.97,63.96)(0.15,0.2){2}{\line(0,1){0.2}}
\multiput(108.26,64.37)(0.14,0.21){2}{\line(0,1){0.21}}
\multiput(108.53,64.79)(0.12,0.22){2}{\line(0,1){0.22}}
\multiput(108.78,65.22)(0.11,0.22){2}{\line(0,1){0.22}}
\multiput(109.01,65.66)(0.1,0.23){2}{\line(0,1){0.23}}
\multiput(109.21,66.12)(0.09,0.23){2}{\line(0,1){0.23}}
\multiput(109.4,66.58)(0.16,0.47){1}{\line(0,1){0.47}}
\multiput(109.56,67.05)(0.14,0.48){1}{\line(0,1){0.48}}
\multiput(109.69,67.53)(0.11,0.49){1}{\line(0,1){0.49}}
\multiput(109.8,68.02)(0.09,0.49){1}{\line(0,1){0.49}}
\multiput(109.89,68.51)(0.06,0.49){1}{\line(0,1){0.49}}
\multiput(109.95,69)(0.04,0.5){1}{\line(0,1){0.5}}
\multiput(109.99,69.5)(0.01,0.5){1}{\line(0,1){0.5}}

\linethickness{0.3mm}
\put(80,70){\line(0,1){10}}
\linethickness{0.3mm}
\put(120,70){\line(0,1){10}}
\linethickness{0.3mm}
\multiput(80,70)(0.01,-0.5){1}{\line(0,-1){0.5}}
\multiput(80.01,69.5)(0.02,-0.5){1}{\line(0,-1){0.5}}
\multiput(80.02,69)(0.03,-0.5){1}{\line(0,-1){0.5}}
\multiput(80.06,68.51)(0.04,-0.5){1}{\line(0,-1){0.5}}
\multiput(80.1,68.01)(0.06,-0.5){1}{\line(0,-1){0.5}}
\multiput(80.16,67.51)(0.07,-0.49){1}{\line(0,-1){0.49}}
\multiput(80.22,67.02)(0.08,-0.49){1}{\line(0,-1){0.49}}
\multiput(80.3,66.53)(0.09,-0.49){1}{\line(0,-1){0.49}}
\multiput(80.4,66.04)(0.1,-0.49){1}{\line(0,-1){0.49}}
\multiput(80.5,65.55)(0.12,-0.48){1}{\line(0,-1){0.48}}
\multiput(80.62,65.06)(0.13,-0.48){1}{\line(0,-1){0.48}}
\multiput(80.75,64.58)(0.14,-0.48){1}{\line(0,-1){0.48}}
\multiput(80.89,64.1)(0.15,-0.47){1}{\line(0,-1){0.47}}
\multiput(81.04,63.63)(0.16,-0.47){1}{\line(0,-1){0.47}}
\multiput(81.21,63.16)(0.18,-0.47){1}{\line(0,-1){0.47}}
\multiput(81.38,62.69)(0.09,-0.23){2}{\line(0,-1){0.23}}
\multiput(81.57,62.23)(0.1,-0.23){2}{\line(0,-1){0.23}}
\multiput(81.77,61.77)(0.11,-0.23){2}{\line(0,-1){0.23}}
\multiput(81.98,61.32)(0.11,-0.22){2}{\line(0,-1){0.22}}
\multiput(82.2,60.88)(0.12,-0.22){2}{\line(0,-1){0.22}}
\multiput(82.44,60.43)(0.12,-0.22){2}{\line(0,-1){0.22}}
\multiput(82.68,60)(0.13,-0.21){2}{\line(0,-1){0.21}}
\multiput(82.93,59.57)(0.13,-0.21){2}{\line(0,-1){0.21}}
\multiput(83.2,59.15)(0.14,-0.21){2}{\line(0,-1){0.21}}
\multiput(83.48,58.73)(0.14,-0.2){2}{\line(0,-1){0.2}}
\multiput(83.76,58.33)(0.15,-0.2){2}{\line(0,-1){0.2}}
\multiput(84.06,57.92)(0.1,-0.13){3}{\line(0,-1){0.13}}
\multiput(84.36,57.53)(0.11,-0.13){3}{\line(0,-1){0.13}}
\multiput(84.68,57.14)(0.11,-0.13){3}{\line(0,-1){0.13}}
\multiput(85,56.77)(0.11,-0.12){3}{\line(0,-1){0.12}}
\multiput(85.34,56.4)(0.11,-0.12){3}{\line(0,-1){0.12}}
\multiput(85.68,56.04)(0.12,-0.12){3}{\line(0,-1){0.12}}
\multiput(86.04,55.68)(0.12,-0.11){3}{\line(1,0){0.12}}
\multiput(86.4,55.34)(0.12,-0.11){3}{\line(1,0){0.12}}
\multiput(86.77,55)(0.13,-0.11){3}{\line(1,0){0.13}}
\multiput(87.14,54.68)(0.13,-0.11){3}{\line(1,0){0.13}}
\multiput(87.53,54.36)(0.13,-0.1){3}{\line(1,0){0.13}}
\multiput(87.92,54.06)(0.2,-0.15){2}{\line(1,0){0.2}}
\multiput(88.33,53.76)(0.2,-0.14){2}{\line(1,0){0.2}}
\multiput(88.73,53.48)(0.21,-0.14){2}{\line(1,0){0.21}}
\multiput(89.15,53.2)(0.21,-0.13){2}{\line(1,0){0.21}}
\multiput(89.57,52.93)(0.21,-0.13){2}{\line(1,0){0.21}}
\multiput(90,52.68)(0.22,-0.12){2}{\line(1,0){0.22}}
\multiput(90.43,52.44)(0.22,-0.12){2}{\line(1,0){0.22}}
\multiput(90.88,52.2)(0.22,-0.11){2}{\line(1,0){0.22}}
\multiput(91.32,51.98)(0.23,-0.11){2}{\line(1,0){0.23}}
\multiput(91.77,51.77)(0.23,-0.1){2}{\line(1,0){0.23}}
\multiput(92.23,51.57)(0.23,-0.09){2}{\line(1,0){0.23}}
\multiput(92.69,51.38)(0.47,-0.18){1}{\line(1,0){0.47}}
\multiput(93.16,51.21)(0.47,-0.16){1}{\line(1,0){0.47}}
\multiput(93.63,51.04)(0.47,-0.15){1}{\line(1,0){0.47}}
\multiput(94.1,50.89)(0.48,-0.14){1}{\line(1,0){0.48}}
\multiput(94.58,50.75)(0.48,-0.13){1}{\line(1,0){0.48}}
\multiput(95.06,50.62)(0.48,-0.12){1}{\line(1,0){0.48}}
\multiput(95.55,50.5)(0.49,-0.1){1}{\line(1,0){0.49}}
\multiput(96.04,50.4)(0.49,-0.09){1}{\line(1,0){0.49}}
\multiput(96.53,50.3)(0.49,-0.08){1}{\line(1,0){0.49}}
\multiput(97.02,50.22)(0.49,-0.07){1}{\line(1,0){0.49}}
\multiput(97.51,50.16)(0.5,-0.06){1}{\line(1,0){0.5}}
\multiput(98.01,50.1)(0.5,-0.04){1}{\line(1,0){0.5}}
\multiput(98.51,50.06)(0.5,-0.03){1}{\line(1,0){0.5}}
\multiput(99,50.02)(0.5,-0.02){1}{\line(1,0){0.5}}
\multiput(99.5,50.01)(0.5,-0.01){1}{\line(1,0){0.5}}
\multiput(100,50)(0.5,0.01){1}{\line(1,0){0.5}}
\multiput(100.5,50.01)(0.5,0.02){1}{\line(1,0){0.5}}
\multiput(101,50.02)(0.5,0.03){1}{\line(1,0){0.5}}
\multiput(101.49,50.06)(0.5,0.04){1}{\line(1,0){0.5}}
\multiput(101.99,50.1)(0.5,0.06){1}{\line(1,0){0.5}}
\multiput(102.49,50.16)(0.49,0.07){1}{\line(1,0){0.49}}
\multiput(102.98,50.22)(0.49,0.08){1}{\line(1,0){0.49}}
\multiput(103.47,50.3)(0.49,0.09){1}{\line(1,0){0.49}}
\multiput(103.96,50.4)(0.49,0.1){1}{\line(1,0){0.49}}
\multiput(104.45,50.5)(0.48,0.12){1}{\line(1,0){0.48}}
\multiput(104.94,50.62)(0.48,0.13){1}{\line(1,0){0.48}}
\multiput(105.42,50.75)(0.48,0.14){1}{\line(1,0){0.48}}
\multiput(105.9,50.89)(0.47,0.15){1}{\line(1,0){0.47}}
\multiput(106.37,51.04)(0.47,0.16){1}{\line(1,0){0.47}}
\multiput(106.84,51.21)(0.47,0.18){1}{\line(1,0){0.47}}
\multiput(107.31,51.38)(0.23,0.09){2}{\line(1,0){0.23}}
\multiput(107.77,51.57)(0.23,0.1){2}{\line(1,0){0.23}}
\multiput(108.23,51.77)(0.23,0.11){2}{\line(1,0){0.23}}
\multiput(108.68,51.98)(0.22,0.11){2}{\line(1,0){0.22}}
\multiput(109.12,52.2)(0.22,0.12){2}{\line(1,0){0.22}}
\multiput(109.57,52.44)(0.22,0.12){2}{\line(1,0){0.22}}
\multiput(110,52.68)(0.21,0.13){2}{\line(1,0){0.21}}
\multiput(110.43,52.93)(0.21,0.13){2}{\line(1,0){0.21}}
\multiput(110.85,53.2)(0.21,0.14){2}{\line(1,0){0.21}}
\multiput(111.27,53.48)(0.2,0.14){2}{\line(1,0){0.2}}
\multiput(111.67,53.76)(0.2,0.15){2}{\line(1,0){0.2}}
\multiput(112.08,54.06)(0.13,0.1){3}{\line(1,0){0.13}}
\multiput(112.47,54.36)(0.13,0.11){3}{\line(1,0){0.13}}
\multiput(112.86,54.68)(0.13,0.11){3}{\line(1,0){0.13}}
\multiput(113.23,55)(0.12,0.11){3}{\line(1,0){0.12}}
\multiput(113.6,55.34)(0.12,0.11){3}{\line(1,0){0.12}}
\multiput(113.96,55.68)(0.12,0.12){3}{\line(0,1){0.12}}
\multiput(114.32,56.04)(0.11,0.12){3}{\line(0,1){0.12}}
\multiput(114.66,56.4)(0.11,0.12){3}{\line(0,1){0.12}}
\multiput(115,56.77)(0.11,0.13){3}{\line(0,1){0.13}}
\multiput(115.32,57.14)(0.11,0.13){3}{\line(0,1){0.13}}
\multiput(115.64,57.53)(0.1,0.13){3}{\line(0,1){0.13}}
\multiput(115.94,57.92)(0.15,0.2){2}{\line(0,1){0.2}}
\multiput(116.24,58.33)(0.14,0.2){2}{\line(0,1){0.2}}
\multiput(116.52,58.73)(0.14,0.21){2}{\line(0,1){0.21}}
\multiput(116.8,59.15)(0.13,0.21){2}{\line(0,1){0.21}}
\multiput(117.07,59.57)(0.13,0.21){2}{\line(0,1){0.21}}
\multiput(117.32,60)(0.12,0.22){2}{\line(0,1){0.22}}
\multiput(117.56,60.43)(0.12,0.22){2}{\line(0,1){0.22}}
\multiput(117.8,60.88)(0.11,0.22){2}{\line(0,1){0.22}}
\multiput(118.02,61.32)(0.11,0.23){2}{\line(0,1){0.23}}
\multiput(118.23,61.77)(0.1,0.23){2}{\line(0,1){0.23}}
\multiput(118.43,62.23)(0.09,0.23){2}{\line(0,1){0.23}}
\multiput(118.62,62.69)(0.18,0.47){1}{\line(0,1){0.47}}
\multiput(118.79,63.16)(0.16,0.47){1}{\line(0,1){0.47}}
\multiput(118.96,63.63)(0.15,0.47){1}{\line(0,1){0.47}}
\multiput(119.11,64.1)(0.14,0.48){1}{\line(0,1){0.48}}
\multiput(119.25,64.58)(0.13,0.48){1}{\line(0,1){0.48}}
\multiput(119.38,65.06)(0.12,0.48){1}{\line(0,1){0.48}}
\multiput(119.5,65.55)(0.1,0.49){1}{\line(0,1){0.49}}
\multiput(119.6,66.04)(0.09,0.49){1}{\line(0,1){0.49}}
\multiput(119.7,66.53)(0.08,0.49){1}{\line(0,1){0.49}}
\multiput(119.78,67.02)(0.07,0.49){1}{\line(0,1){0.49}}
\multiput(119.84,67.51)(0.06,0.5){1}{\line(0,1){0.5}}
\multiput(119.9,68.01)(0.04,0.5){1}{\line(0,1){0.5}}
\multiput(119.94,68.51)(0.03,0.5){1}{\line(0,1){0.5}}
\multiput(119.98,69)(0.02,0.5){1}{\line(0,1){0.5}}
\multiput(119.99,69.5)(0.01,0.5){1}{\line(0,1){0.5}}

\put(75,85){\makebox(0,0)[cc]{$\Vec A$}}

\put(88,93){\makebox(0,0)[cc]{$\Vec^{\tens n}A$}}

\put(112,93){\makebox(0,0)[cc]{$\Omega^{\tens n}A$}}

\put(122,85){\makebox(0,0)[cc]{$\Omega^1 A$}}

%%%%%%%%%%%%%%%
\refstepcounter{piccie} \label{dshbshdxh}
\put(182,55){\makebox(0,0)[cc]{Fig.\ \arabic{piccie}}}
%%%%%%%%%%%%%%%

\linethickness{0.3mm}
\put(0,70){\line(0,1){10}}
\linethickness{0.3mm}
\put(40,70){\line(0,1){10}}
\linethickness{0.3mm}
\multiput(0,70)(0.01,-0.5){1}{\line(0,-1){0.5}}
\multiput(0.01,69.5)(0.02,-0.5){1}{\line(0,-1){0.5}}
\multiput(0.02,69)(0.03,-0.5){1}{\line(0,-1){0.5}}
\multiput(0.06,68.51)(0.04,-0.5){1}{\line(0,-1){0.5}}
\multiput(0.1,68.01)(0.06,-0.5){1}{\line(0,-1){0.5}}
\multiput(0.16,67.51)(0.07,-0.49){1}{\line(0,-1){0.49}}
\multiput(0.22,67.02)(0.08,-0.49){1}{\line(0,-1){0.49}}
\multiput(0.3,66.53)(0.09,-0.49){1}{\line(0,-1){0.49}}
\multiput(0.4,66.04)(0.1,-0.49){1}{\line(0,-1){0.49}}
\multiput(0.5,65.55)(0.12,-0.48){1}{\line(0,-1){0.48}}
\multiput(0.62,65.06)(0.13,-0.48){1}{\line(0,-1){0.48}}
\multiput(0.75,64.58)(0.14,-0.48){1}{\line(0,-1){0.48}}
\multiput(0.89,64.1)(0.15,-0.47){1}{\line(0,-1){0.47}}
\multiput(1.04,63.63)(0.16,-0.47){1}{\line(0,-1){0.47}}
\multiput(1.21,63.16)(0.18,-0.47){1}{\line(0,-1){0.47}}
\multiput(1.38,62.69)(0.09,-0.23){2}{\line(0,-1){0.23}}
\multiput(1.57,62.23)(0.1,-0.23){2}{\line(0,-1){0.23}}
\multiput(1.77,61.77)(0.11,-0.23){2}{\line(0,-1){0.23}}
\multiput(1.98,61.32)(0.11,-0.22){2}{\line(0,-1){0.22}}
\multiput(2.2,60.88)(0.12,-0.22){2}{\line(0,-1){0.22}}
\multiput(2.44,60.43)(0.12,-0.22){2}{\line(0,-1){0.22}}
\multiput(2.68,60)(0.13,-0.21){2}{\line(0,-1){0.21}}
\multiput(2.93,59.57)(0.13,-0.21){2}{\line(0,-1){0.21}}
\multiput(3.2,59.15)(0.14,-0.21){2}{\line(0,-1){0.21}}
\multiput(3.48,58.73)(0.14,-0.2){2}{\line(0,-1){0.2}}
\multiput(3.76,58.33)(0.15,-0.2){2}{\line(0,-1){0.2}}
\multiput(4.06,57.92)(0.1,-0.13){3}{\line(0,-1){0.13}}
\multiput(4.36,57.53)(0.11,-0.13){3}{\line(0,-1){0.13}}
\multiput(4.68,57.14)(0.11,-0.13){3}{\line(0,-1){0.13}}
\multiput(5,56.77)(0.11,-0.12){3}{\line(0,-1){0.12}}
\multiput(5.34,56.4)(0.11,-0.12){3}{\line(0,-1){0.12}}
\multiput(5.68,56.04)(0.12,-0.12){3}{\line(0,-1){0.12}}
\multiput(6.04,55.68)(0.12,-0.11){3}{\line(1,0){0.12}}
\multiput(6.4,55.34)(0.12,-0.11){3}{\line(1,0){0.12}}
\multiput(6.77,55)(0.13,-0.11){3}{\line(1,0){0.13}}
\multiput(7.14,54.68)(0.13,-0.11){3}{\line(1,0){0.13}}
\multiput(7.53,54.36)(0.13,-0.1){3}{\line(1,0){0.13}}
\multiput(7.92,54.06)(0.2,-0.15){2}{\line(1,0){0.2}}
\multiput(8.33,53.76)(0.2,-0.14){2}{\line(1,0){0.2}}
\multiput(8.73,53.48)(0.21,-0.14){2}{\line(1,0){0.21}}
\multiput(9.15,53.2)(0.21,-0.13){2}{\line(1,0){0.21}}
\multiput(9.57,52.93)(0.21,-0.13){2}{\line(1,0){0.21}}
\multiput(10,52.68)(0.22,-0.12){2}{\line(1,0){0.22}}
\multiput(10.43,52.44)(0.22,-0.12){2}{\line(1,0){0.22}}
\multiput(10.88,52.2)(0.22,-0.11){2}{\line(1,0){0.22}}
\multiput(11.32,51.98)(0.23,-0.11){2}{\line(1,0){0.23}}
\multiput(11.77,51.77)(0.23,-0.1){2}{\line(1,0){0.23}}
\multiput(12.23,51.57)(0.23,-0.09){2}{\line(1,0){0.23}}
\multiput(12.69,51.38)(0.47,-0.18){1}{\line(1,0){0.47}}
\multiput(13.16,51.21)(0.47,-0.16){1}{\line(1,0){0.47}}
\multiput(13.63,51.04)(0.47,-0.15){1}{\line(1,0){0.47}}
\multiput(14.1,50.89)(0.48,-0.14){1}{\line(1,0){0.48}}
\multiput(14.58,50.75)(0.48,-0.13){1}{\line(1,0){0.48}}
\multiput(15.06,50.62)(0.48,-0.12){1}{\line(1,0){0.48}}
\multiput(15.55,50.5)(0.49,-0.1){1}{\line(1,0){0.49}}
\multiput(16.04,50.4)(0.49,-0.09){1}{\line(1,0){0.49}}
\multiput(16.53,50.3)(0.49,-0.08){1}{\line(1,0){0.49}}
\multiput(17.02,50.22)(0.49,-0.07){1}{\line(1,0){0.49}}
\multiput(17.51,50.16)(0.5,-0.06){1}{\line(1,0){0.5}}
\multiput(18.01,50.1)(0.5,-0.04){1}{\line(1,0){0.5}}
\multiput(18.51,50.06)(0.5,-0.03){1}{\line(1,0){0.5}}
\multiput(19,50.02)(0.5,-0.02){1}{\line(1,0){0.5}}
\multiput(19.5,50.01)(0.5,-0.01){1}{\line(1,0){0.5}}
\multiput(20,50)(0.5,0.01){1}{\line(1,0){0.5}}
\multiput(20.5,50.01)(0.5,0.02){1}{\line(1,0){0.5}}
\multiput(21,50.02)(0.5,0.03){1}{\line(1,0){0.5}}
\multiput(21.49,50.06)(0.5,0.04){1}{\line(1,0){0.5}}
\multiput(21.99,50.1)(0.5,0.06){1}{\line(1,0){0.5}}
\multiput(22.49,50.16)(0.49,0.07){1}{\line(1,0){0.49}}
\multiput(22.98,50.22)(0.49,0.08){1}{\line(1,0){0.49}}
\multiput(23.47,50.3)(0.49,0.09){1}{\line(1,0){0.49}}
\multiput(23.96,50.4)(0.49,0.1){1}{\line(1,0){0.49}}
\multiput(24.45,50.5)(0.48,0.12){1}{\line(1,0){0.48}}
\multiput(24.94,50.62)(0.48,0.13){1}{\line(1,0){0.48}}
\multiput(25.42,50.75)(0.48,0.14){1}{\line(1,0){0.48}}
\multiput(25.9,50.89)(0.47,0.15){1}{\line(1,0){0.47}}
\multiput(26.37,51.04)(0.47,0.16){1}{\line(1,0){0.47}}
\multiput(26.84,51.21)(0.47,0.18){1}{\line(1,0){0.47}}
\multiput(27.31,51.38)(0.23,0.09){2}{\line(1,0){0.23}}
\multiput(27.77,51.57)(0.23,0.1){2}{\line(1,0){0.23}}
\multiput(28.23,51.77)(0.23,0.11){2}{\line(1,0){0.23}}
\multiput(28.68,51.98)(0.22,0.11){2}{\line(1,0){0.22}}
\multiput(29.12,52.2)(0.22,0.12){2}{\line(1,0){0.22}}
\multiput(29.57,52.44)(0.22,0.12){2}{\line(1,0){0.22}}
\multiput(30,52.68)(0.21,0.13){2}{\line(1,0){0.21}}
\multiput(30.43,52.93)(0.21,0.13){2}{\line(1,0){0.21}}
\multiput(30.85,53.2)(0.21,0.14){2}{\line(1,0){0.21}}
\multiput(31.27,53.48)(0.2,0.14){2}{\line(1,0){0.2}}
\multiput(31.67,53.76)(0.2,0.15){2}{\line(1,0){0.2}}
\multiput(32.08,54.06)(0.13,0.1){3}{\line(1,0){0.13}}
\multiput(32.47,54.36)(0.13,0.11){3}{\line(1,0){0.13}}
\multiput(32.86,54.68)(0.13,0.11){3}{\line(1,0){0.13}}
\multiput(33.23,55)(0.12,0.11){3}{\line(1,0){0.12}}
\multiput(33.6,55.34)(0.12,0.11){3}{\line(1,0){0.12}}
\multiput(33.96,55.68)(0.12,0.12){3}{\line(0,1){0.12}}
\multiput(34.32,56.04)(0.11,0.12){3}{\line(0,1){0.12}}
\multiput(34.66,56.4)(0.11,0.12){3}{\line(0,1){0.12}}
\multiput(35,56.77)(0.11,0.13){3}{\line(0,1){0.13}}
\multiput(35.32,57.14)(0.11,0.13){3}{\line(0,1){0.13}}
\multiput(35.64,57.53)(0.1,0.13){3}{\line(0,1){0.13}}
\multiput(35.94,57.92)(0.15,0.2){2}{\line(0,1){0.2}}
\multiput(36.24,58.33)(0.14,0.2){2}{\line(0,1){0.2}}
\multiput(36.52,58.73)(0.14,0.21){2}{\line(0,1){0.21}}
\multiput(36.8,59.15)(0.13,0.21){2}{\line(0,1){0.21}}
\multiput(37.07,59.57)(0.13,0.21){2}{\line(0,1){0.21}}
\multiput(37.32,60)(0.12,0.22){2}{\line(0,1){0.22}}
\multiput(37.56,60.43)(0.12,0.22){2}{\line(0,1){0.22}}
\multiput(37.8,60.88)(0.11,0.22){2}{\line(0,1){0.22}}
\multiput(38.02,61.32)(0.11,0.23){2}{\line(0,1){0.23}}
\multiput(38.23,61.77)(0.1,0.23){2}{\line(0,1){0.23}}
\multiput(38.43,62.23)(0.09,0.23){2}{\line(0,1){0.23}}
\multiput(38.62,62.69)(0.18,0.47){1}{\line(0,1){0.47}}
\multiput(38.79,63.16)(0.16,0.47){1}{\line(0,1){0.47}}
\multiput(38.96,63.63)(0.15,0.47){1}{\line(0,1){0.47}}
\multiput(39.11,64.1)(0.14,0.48){1}{\line(0,1){0.48}}
\multiput(39.25,64.58)(0.13,0.48){1}{\line(0,1){0.48}}
\multiput(39.38,65.06)(0.12,0.48){1}{\line(0,1){0.48}}
\multiput(39.5,65.55)(0.1,0.49){1}{\line(0,1){0.49}}
\multiput(39.6,66.04)(0.09,0.49){1}{\line(0,1){0.49}}
\multiput(39.7,66.53)(0.08,0.49){1}{\line(0,1){0.49}}
\multiput(39.78,67.02)(0.07,0.49){1}{\line(0,1){0.49}}
\multiput(39.84,67.51)(0.06,0.5){1}{\line(0,1){0.5}}
\multiput(39.9,68.01)(0.04,0.5){1}{\line(0,1){0.5}}
\multiput(39.94,68.51)(0.03,0.5){1}{\line(0,1){0.5}}
\multiput(39.98,69)(0.02,0.5){1}{\line(0,1){0.5}}
\multiput(39.99,69.5)(0.01,0.5){1}{\line(0,1){0.5}}

\put(20,45){\makebox(0,0)[cc]{$\ev^{\<n+1\>}$}}

\put(100,45){\makebox(0,0)[cc]{$\ev$}}

\put(99,56){\makebox(0,0)[cc]{$\ev^{\<n\>}$}}

\end{picture}

\section{Covariant derivatives on fields and forms}
Consider a right covariant derivative on $\Omega^1 A$ and the corresponding dual 
left covariant derivative on $\Vec A$. We then extend these to $n$-tuples of forms and fields, and consider the generalised braiding given by bimodule covariant derivatives.

\subsection{A right connection on tensor products of 1-forms} \label{vcgsajhdx}
We begin with the right handed analogue of Definitions~\ref{bcaiskyf} and \ref{ppll}.
Suppose that $\square:\Omega^1 A\to \Omega^1 A\tens_A \Omega^1 A$ is a right bimodule covariant derivative, so it satisfies, for
$\xi\in \Omega^1 A$ and $a\in A$,
\begin{eqnarray}\label{sjgacvaghjchjk}
\square(\xi.a) &=& \square(\xi).a+\xi\tens\extd a\ ,\cr
\square(a.\xi) &=& a.\square(\xi)+\sigma^{-1}(\extd a\tens\xi)\ .
\end{eqnarray}
Here $\sigma^{-1}:\Omega^1 A\tens_A \Omega^1 A\to \Omega^1 A\tens_A \Omega^1 A$ is a  bimodule map, and we shall assume that it is invertible with inverse $\sigma$. [The reader should note that the inverse
in $\sigma^{-1}$ is a convention, set to be compatible with the use of $\sigma$ for left bimodule covariant derivatives, and to be compatible with the notation for braided categories as previously described.]  Now $\square$ extends to a right bimodule covariant derivative $\square^{\<n\>}:\Omega^{\tens n}A\to \Omega^{\tens n+1}A$, defined recursively by
\begin{eqnarray*}
\square^{\<0\>} &=& \extd\ ,\quad \square^{\<1\>}\ =\ \square\ ,\cr
\square^{\<n+1\>} &=& \id^{\tens n}\tens\square+(\id^{\tens n}\tens\sigma^{-1})(\square^{\<n\>}\tens\id^{\tens 1})\ .
\end{eqnarray*}

\subsection{A left connection on tensor products of vector fields}\label{gcysiuaf}
Suppose that $\square:\Omega^1 A\to \Omega^1 A\tens_A \Omega^1 A$ is a right bimodule covariant derivative, with generalised braiding $\sigma^{-1}:\Omega^1 A\tens_A \Omega^1 A\to \Omega^1 A\tens_A \Omega^1 A$ invertible (the inverse being denoted $\sigma$).
Then there is a left covariant derivative $\square:\Vec A\to \Omega^1 A\tens_A \Vec A$
so that
\begin{eqnarray} \label{ghjksvdgwv}
\extd\circ\ev\ =\ (\id\tens\ev)(\square\tens\id)+(\ev\tens\id)(\id\tens\square)\ :\ 
\Vec A\tens_A \Omega^1 A\to\Omega^1  A\ .
\end{eqnarray}
[It will be convenient to use $\square$ for both this covariant derivative
and the one on $\Omega^1 A$, and distinguish them by their domains.] 
We give a pictorial interpretation of (\ref{ghjksvdgwv}) in Fig.\ \ref{dshtttt}. Note that in the diagram to the left of the equality in Fig.\ \ref{dshtttt} we assign a visible line to $A$, so that we can see what the derivative $\extd:A\to \Omega^1 A$ applies to. 

\unitlength 0.5 mm
\begin{picture}(160,70)(-20,8)
\linethickness{0.3mm}
\multiput(10,60)(0.01,-0.5){1}{\line(0,-1){0.5}}
\multiput(10.01,59.5)(0.04,-0.5){1}{\line(0,-1){0.5}}
\multiput(10.05,59)(0.06,-0.49){1}{\line(0,-1){0.49}}
\multiput(10.11,58.51)(0.09,-0.49){1}{\line(0,-1){0.49}}
\multiput(10.2,58.02)(0.11,-0.49){1}{\line(0,-1){0.49}}
\multiput(10.31,57.53)(0.14,-0.48){1}{\line(0,-1){0.48}}
\multiput(10.44,57.05)(0.16,-0.47){1}{\line(0,-1){0.47}}
\multiput(10.6,56.58)(0.09,-0.23){2}{\line(0,-1){0.23}}
\multiput(10.79,56.12)(0.1,-0.23){2}{\line(0,-1){0.23}}
\multiput(10.99,55.66)(0.11,-0.22){2}{\line(0,-1){0.22}}
\multiput(11.22,55.22)(0.12,-0.22){2}{\line(0,-1){0.22}}
\multiput(11.47,54.79)(0.14,-0.21){2}{\line(0,-1){0.21}}
\multiput(11.74,54.37)(0.15,-0.2){2}{\line(0,-1){0.2}}
\multiput(12.03,53.96)(0.1,-0.13){3}{\line(0,-1){0.13}}
\multiput(12.34,53.57)(0.11,-0.12){3}{\line(0,-1){0.12}}
\multiput(12.67,53.2)(0.12,-0.12){3}{\line(0,-1){0.12}}
\multiput(13.02,52.84)(0.12,-0.11){3}{\line(1,0){0.12}}
\multiput(13.38,52.5)(0.13,-0.11){3}{\line(1,0){0.13}}
\multiput(13.77,52.18)(0.13,-0.1){3}{\line(1,0){0.13}}
\multiput(14.16,51.88)(0.21,-0.14){2}{\line(1,0){0.21}}
\multiput(14.57,51.6)(0.21,-0.13){2}{\line(1,0){0.21}}
\multiput(15,51.34)(0.22,-0.12){2}{\line(1,0){0.22}}
\multiput(15.44,51.1)(0.22,-0.11){2}{\line(1,0){0.22}}
\multiput(15.89,50.88)(0.23,-0.1){2}{\line(1,0){0.23}}
\multiput(16.35,50.69)(0.47,-0.17){1}{\line(1,0){0.47}}
\multiput(16.82,50.52)(0.48,-0.15){1}{\line(1,0){0.48}}
\multiput(17.29,50.37)(0.48,-0.12){1}{\line(1,0){0.48}}
\multiput(17.77,50.25)(0.49,-0.1){1}{\line(1,0){0.49}}
\multiput(18.26,50.15)(0.49,-0.07){1}{\line(1,0){0.49}}
\multiput(18.76,50.08)(0.5,-0.05){1}{\line(1,0){0.5}}
\multiput(19.25,50.03)(0.5,-0.02){1}{\line(1,0){0.5}}
\put(19.75,50){\line(1,0){0.5}}
\multiput(20.25,50)(0.5,0.02){1}{\line(1,0){0.5}}
\multiput(20.75,50.03)(0.5,0.05){1}{\line(1,0){0.5}}
\multiput(21.24,50.08)(0.49,0.07){1}{\line(1,0){0.49}}
\multiput(21.74,50.15)(0.49,0.1){1}{\line(1,0){0.49}}
\multiput(22.23,50.25)(0.48,0.12){1}{\line(1,0){0.48}}
\multiput(22.71,50.37)(0.48,0.15){1}{\line(1,0){0.48}}
\multiput(23.18,50.52)(0.47,0.17){1}{\line(1,0){0.47}}
\multiput(23.65,50.69)(0.23,0.1){2}{\line(1,0){0.23}}
\multiput(24.11,50.88)(0.22,0.11){2}{\line(1,0){0.22}}
\multiput(24.56,51.1)(0.22,0.12){2}{\line(1,0){0.22}}
\multiput(25,51.34)(0.21,0.13){2}{\line(1,0){0.21}}
\multiput(25.43,51.6)(0.21,0.14){2}{\line(1,0){0.21}}
\multiput(25.84,51.88)(0.13,0.1){3}{\line(1,0){0.13}}
\multiput(26.23,52.18)(0.13,0.11){3}{\line(1,0){0.13}}
\multiput(26.62,52.5)(0.12,0.11){3}{\line(1,0){0.12}}
\multiput(26.98,52.84)(0.12,0.12){3}{\line(0,1){0.12}}
\multiput(27.33,53.2)(0.11,0.12){3}{\line(0,1){0.12}}
\multiput(27.66,53.57)(0.1,0.13){3}{\line(0,1){0.13}}
\multiput(27.97,53.96)(0.15,0.2){2}{\line(0,1){0.2}}
\multiput(28.26,54.37)(0.14,0.21){2}{\line(0,1){0.21}}
\multiput(28.53,54.79)(0.12,0.22){2}{\line(0,1){0.22}}
\multiput(28.78,55.22)(0.11,0.22){2}{\line(0,1){0.22}}
\multiput(29.01,55.66)(0.1,0.23){2}{\line(0,1){0.23}}
\multiput(29.21,56.12)(0.09,0.23){2}{\line(0,1){0.23}}
\multiput(29.4,56.58)(0.16,0.47){1}{\line(0,1){0.47}}
\multiput(29.56,57.05)(0.14,0.48){1}{\line(0,1){0.48}}
\multiput(29.69,57.53)(0.11,0.49){1}{\line(0,1){0.49}}
\multiput(29.8,58.02)(0.09,0.49){1}{\line(0,1){0.49}}
\multiput(29.89,58.51)(0.06,0.49){1}{\line(0,1){0.49}}
\multiput(29.95,59)(0.04,0.5){1}{\line(0,1){0.5}}
\multiput(29.99,59.5)(0.01,0.5){1}{\line(0,1){0.5}}

\linethickness{0.3mm}
\put(20,40){\line(0,1){10}}
\put(20,35){\makebox(0,0)[cc]{$\extd$}}

\put(15,35){\makebox(0,0)[cc]{}}

\put(15,35){\makebox(0,0)[cc]{}}

\linethickness{0.3mm}
\put(20,20){\line(0,1){10}}
\put(45,40){\makebox(0,0)[cc]{=}}

%%%%%%%%%%%%%%%
\refstepcounter{piccie} \label{dshtttt}
\put(200,40){\makebox(0,0)[cc]{Fig.\ \arabic{piccie}}}
%%%%%%%%%%%%%%%

\linethickness{0.3mm}
\put(70,50){\line(0,1){10}}
\linethickness{0.3mm}
\multiput(79.99,40.5)(0.01,-0.5){1}{\line(0,-1){0.5}}
\multiput(79.95,41)(0.04,-0.5){1}{\line(0,-1){0.5}}
\multiput(79.89,41.49)(0.06,-0.49){1}{\line(0,-1){0.49}}
\multiput(79.8,41.98)(0.09,-0.49){1}{\line(0,-1){0.49}}
\multiput(79.69,42.47)(0.11,-0.49){1}{\line(0,-1){0.49}}
\multiput(79.56,42.95)(0.14,-0.48){1}{\line(0,-1){0.48}}
\multiput(79.4,43.42)(0.16,-0.47){1}{\line(0,-1){0.47}}
\multiput(79.21,43.88)(0.09,-0.23){2}{\line(0,-1){0.23}}
\multiput(79.01,44.34)(0.1,-0.23){2}{\line(0,-1){0.23}}
\multiput(78.78,44.78)(0.11,-0.22){2}{\line(0,-1){0.22}}
\multiput(78.53,45.21)(0.12,-0.22){2}{\line(0,-1){0.22}}
\multiput(78.26,45.63)(0.14,-0.21){2}{\line(0,-1){0.21}}
\multiput(77.97,46.04)(0.15,-0.2){2}{\line(0,-1){0.2}}
\multiput(77.66,46.43)(0.1,-0.13){3}{\line(0,-1){0.13}}
\multiput(77.33,46.8)(0.11,-0.12){3}{\line(0,-1){0.12}}
\multiput(76.98,47.16)(0.12,-0.12){3}{\line(0,-1){0.12}}
\multiput(76.62,47.5)(0.12,-0.11){3}{\line(1,0){0.12}}
\multiput(76.23,47.82)(0.13,-0.11){3}{\line(1,0){0.13}}
\multiput(75.84,48.12)(0.13,-0.1){3}{\line(1,0){0.13}}
\multiput(75.43,48.4)(0.21,-0.14){2}{\line(1,0){0.21}}
\multiput(75,48.66)(0.21,-0.13){2}{\line(1,0){0.21}}
\multiput(74.56,48.9)(0.22,-0.12){2}{\line(1,0){0.22}}
\multiput(74.11,49.12)(0.22,-0.11){2}{\line(1,0){0.22}}
\multiput(73.65,49.31)(0.23,-0.1){2}{\line(1,0){0.23}}
\multiput(73.18,49.48)(0.47,-0.17){1}{\line(1,0){0.47}}
\multiput(72.71,49.63)(0.48,-0.15){1}{\line(1,0){0.48}}
\multiput(72.23,49.75)(0.48,-0.12){1}{\line(1,0){0.48}}
\multiput(71.74,49.85)(0.49,-0.1){1}{\line(1,0){0.49}}
\multiput(71.24,49.92)(0.49,-0.07){1}{\line(1,0){0.49}}
\multiput(70.75,49.97)(0.5,-0.05){1}{\line(1,0){0.5}}
\multiput(70.25,50)(0.5,-0.02){1}{\line(1,0){0.5}}
\put(69.75,50){\line(1,0){0.5}}
\multiput(69.25,49.97)(0.5,0.02){1}{\line(1,0){0.5}}
\multiput(68.76,49.92)(0.5,0.05){1}{\line(1,0){0.5}}
\multiput(68.26,49.85)(0.49,0.07){1}{\line(1,0){0.49}}
\multiput(67.77,49.75)(0.49,0.1){1}{\line(1,0){0.49}}
\multiput(67.29,49.63)(0.48,0.12){1}{\line(1,0){0.48}}
\multiput(66.82,49.48)(0.48,0.15){1}{\line(1,0){0.48}}
\multiput(66.35,49.31)(0.47,0.17){1}{\line(1,0){0.47}}
\multiput(65.89,49.12)(0.23,0.1){2}{\line(1,0){0.23}}
\multiput(65.44,48.9)(0.22,0.11){2}{\line(1,0){0.22}}
\multiput(65,48.66)(0.22,0.12){2}{\line(1,0){0.22}}
\multiput(64.57,48.4)(0.21,0.13){2}{\line(1,0){0.21}}
\multiput(64.16,48.12)(0.21,0.14){2}{\line(1,0){0.21}}
\multiput(63.77,47.82)(0.13,0.1){3}{\line(1,0){0.13}}
\multiput(63.38,47.5)(0.13,0.11){3}{\line(1,0){0.13}}
\multiput(63.02,47.16)(0.12,0.11){3}{\line(1,0){0.12}}
\multiput(62.67,46.8)(0.12,0.12){3}{\line(0,1){0.12}}
\multiput(62.34,46.43)(0.11,0.12){3}{\line(0,1){0.12}}
\multiput(62.03,46.04)(0.1,0.13){3}{\line(0,1){0.13}}
\multiput(61.74,45.63)(0.15,0.2){2}{\line(0,1){0.2}}
\multiput(61.47,45.21)(0.14,0.21){2}{\line(0,1){0.21}}
\multiput(61.22,44.78)(0.12,0.22){2}{\line(0,1){0.22}}
\multiput(60.99,44.34)(0.11,0.22){2}{\line(0,1){0.22}}
\multiput(60.79,43.88)(0.1,0.23){2}{\line(0,1){0.23}}
\multiput(60.6,43.42)(0.09,0.23){2}{\line(0,1){0.23}}
\multiput(60.44,42.95)(0.16,0.47){1}{\line(0,1){0.47}}
\multiput(60.31,42.47)(0.14,0.48){1}{\line(0,1){0.48}}
\multiput(60.2,41.98)(0.11,0.49){1}{\line(0,1){0.49}}
\multiput(60.11,41.49)(0.09,0.49){1}{\line(0,1){0.49}}
\multiput(60.05,41)(0.06,0.49){1}{\line(0,1){0.49}}
\multiput(60.01,40.5)(0.04,0.5){1}{\line(0,1){0.5}}
\multiput(60,40)(0.01,0.5){1}{\line(0,1){0.5}}

\linethickness{0.3mm}
\multiput(80,40)(0.01,-0.5){1}{\line(0,-1){0.5}}
\multiput(80.01,39.5)(0.04,-0.5){1}{\line(0,-1){0.5}}
\multiput(80.05,39)(0.06,-0.49){1}{\line(0,-1){0.49}}
\multiput(80.11,38.51)(0.09,-0.49){1}{\line(0,-1){0.49}}
\multiput(80.2,38.02)(0.11,-0.49){1}{\line(0,-1){0.49}}
\multiput(80.31,37.53)(0.14,-0.48){1}{\line(0,-1){0.48}}
\multiput(80.44,37.05)(0.16,-0.47){1}{\line(0,-1){0.47}}
\multiput(80.6,36.58)(0.09,-0.23){2}{\line(0,-1){0.23}}
\multiput(80.79,36.12)(0.1,-0.23){2}{\line(0,-1){0.23}}
\multiput(80.99,35.66)(0.11,-0.22){2}{\line(0,-1){0.22}}
\multiput(81.22,35.22)(0.12,-0.22){2}{\line(0,-1){0.22}}
\multiput(81.47,34.79)(0.14,-0.21){2}{\line(0,-1){0.21}}
\multiput(81.74,34.37)(0.15,-0.2){2}{\line(0,-1){0.2}}
\multiput(82.03,33.96)(0.1,-0.13){3}{\line(0,-1){0.13}}
\multiput(82.34,33.57)(0.11,-0.12){3}{\line(0,-1){0.12}}
\multiput(82.67,33.2)(0.12,-0.12){3}{\line(0,-1){0.12}}
\multiput(83.02,32.84)(0.12,-0.11){3}{\line(1,0){0.12}}
\multiput(83.38,32.5)(0.13,-0.11){3}{\line(1,0){0.13}}
\multiput(83.77,32.18)(0.13,-0.1){3}{\line(1,0){0.13}}
\multiput(84.16,31.88)(0.21,-0.14){2}{\line(1,0){0.21}}
\multiput(84.57,31.6)(0.21,-0.13){2}{\line(1,0){0.21}}
\multiput(85,31.34)(0.22,-0.12){2}{\line(1,0){0.22}}
\multiput(85.44,31.1)(0.22,-0.11){2}{\line(1,0){0.22}}
\multiput(85.89,30.88)(0.23,-0.1){2}{\line(1,0){0.23}}
\multiput(86.35,30.69)(0.47,-0.17){1}{\line(1,0){0.47}}
\multiput(86.82,30.52)(0.48,-0.15){1}{\line(1,0){0.48}}
\multiput(87.29,30.37)(0.48,-0.12){1}{\line(1,0){0.48}}
\multiput(87.77,30.25)(0.49,-0.1){1}{\line(1,0){0.49}}
\multiput(88.26,30.15)(0.49,-0.07){1}{\line(1,0){0.49}}
\multiput(88.76,30.08)(0.5,-0.05){1}{\line(1,0){0.5}}
\multiput(89.25,30.03)(0.5,-0.02){1}{\line(1,0){0.5}}
\put(89.75,30){\line(1,0){0.5}}
\multiput(90.25,30)(0.5,0.02){1}{\line(1,0){0.5}}
\multiput(90.75,30.03)(0.5,0.05){1}{\line(1,0){0.5}}
\multiput(91.24,30.08)(0.49,0.07){1}{\line(1,0){0.49}}
\multiput(91.74,30.15)(0.49,0.1){1}{\line(1,0){0.49}}
\multiput(92.23,30.25)(0.48,0.12){1}{\line(1,0){0.48}}
\multiput(92.71,30.37)(0.48,0.15){1}{\line(1,0){0.48}}
\multiput(93.18,30.52)(0.47,0.17){1}{\line(1,0){0.47}}
\multiput(93.65,30.69)(0.23,0.1){2}{\line(1,0){0.23}}
\multiput(94.11,30.88)(0.22,0.11){2}{\line(1,0){0.22}}
\multiput(94.56,31.1)(0.22,0.12){2}{\line(1,0){0.22}}
\multiput(95,31.34)(0.21,0.13){2}{\line(1,0){0.21}}
\multiput(95.43,31.6)(0.21,0.14){2}{\line(1,0){0.21}}
\multiput(95.84,31.88)(0.13,0.1){3}{\line(1,0){0.13}}
\multiput(96.23,32.18)(0.13,0.11){3}{\line(1,0){0.13}}
\multiput(96.62,32.5)(0.12,0.11){3}{\line(1,0){0.12}}
\multiput(96.98,32.84)(0.12,0.12){3}{\line(0,1){0.12}}
\multiput(97.33,33.2)(0.11,0.12){3}{\line(0,1){0.12}}
\multiput(97.66,33.57)(0.1,0.13){3}{\line(0,1){0.13}}
\multiput(97.97,33.96)(0.15,0.2){2}{\line(0,1){0.2}}
\multiput(98.26,34.37)(0.14,0.21){2}{\line(0,1){0.21}}
\multiput(98.53,34.79)(0.12,0.22){2}{\line(0,1){0.22}}
\multiput(98.78,35.22)(0.11,0.22){2}{\line(0,1){0.22}}
\multiput(99.01,35.66)(0.1,0.23){2}{\line(0,1){0.23}}
\multiput(99.21,36.12)(0.09,0.23){2}{\line(0,1){0.23}}
\multiput(99.4,36.58)(0.16,0.47){1}{\line(0,1){0.47}}
\multiput(99.56,37.05)(0.14,0.48){1}{\line(0,1){0.48}}
\multiput(99.69,37.53)(0.11,0.49){1}{\line(0,1){0.49}}
\multiput(99.8,38.02)(0.09,0.49){1}{\line(0,1){0.49}}
\multiput(99.89,38.51)(0.06,0.49){1}{\line(0,1){0.49}}
\multiput(99.95,39)(0.04,0.5){1}{\line(0,1){0.5}}
\multiput(99.99,39.5)(0.01,0.5){1}{\line(0,1){0.5}}

\linethickness{0.3mm}
\put(100,40){\line(0,1){20}}
\linethickness{0.3mm}
\put(60,20){\line(0,1){20}}
\put(110,45){\makebox(0,0)[cc]{$+$}}

\linethickness{0.3mm}
\put(120,40){\line(0,1){20}}
\linethickness{0.3mm}
\multiput(120,40)(0.01,-0.5){1}{\line(0,-1){0.5}}
\multiput(120.01,39.5)(0.04,-0.5){1}{\line(0,-1){0.5}}
\multiput(120.05,39)(0.06,-0.49){1}{\line(0,-1){0.49}}
\multiput(120.11,38.51)(0.09,-0.49){1}{\line(0,-1){0.49}}
\multiput(120.2,38.02)(0.11,-0.49){1}{\line(0,-1){0.49}}
\multiput(120.31,37.53)(0.14,-0.48){1}{\line(0,-1){0.48}}
\multiput(120.44,37.05)(0.16,-0.47){1}{\line(0,-1){0.47}}
\multiput(120.6,36.58)(0.09,-0.23){2}{\line(0,-1){0.23}}
\multiput(120.79,36.12)(0.1,-0.23){2}{\line(0,-1){0.23}}
\multiput(120.99,35.66)(0.11,-0.22){2}{\line(0,-1){0.22}}
\multiput(121.22,35.22)(0.12,-0.22){2}{\line(0,-1){0.22}}
\multiput(121.47,34.79)(0.14,-0.21){2}{\line(0,-1){0.21}}
\multiput(121.74,34.37)(0.15,-0.2){2}{\line(0,-1){0.2}}
\multiput(122.03,33.96)(0.1,-0.13){3}{\line(0,-1){0.13}}
\multiput(122.34,33.57)(0.11,-0.12){3}{\line(0,-1){0.12}}
\multiput(122.67,33.2)(0.12,-0.12){3}{\line(0,-1){0.12}}
\multiput(123.02,32.84)(0.12,-0.11){3}{\line(1,0){0.12}}
\multiput(123.38,32.5)(0.13,-0.11){3}{\line(1,0){0.13}}
\multiput(123.77,32.18)(0.13,-0.1){3}{\line(1,0){0.13}}
\multiput(124.16,31.88)(0.21,-0.14){2}{\line(1,0){0.21}}
\multiput(124.57,31.6)(0.21,-0.13){2}{\line(1,0){0.21}}
\multiput(125,31.34)(0.22,-0.12){2}{\line(1,0){0.22}}
\multiput(125.44,31.1)(0.22,-0.11){2}{\line(1,0){0.22}}
\multiput(125.89,30.88)(0.23,-0.1){2}{\line(1,0){0.23}}
\multiput(126.35,30.69)(0.47,-0.17){1}{\line(1,0){0.47}}
\multiput(126.82,30.52)(0.48,-0.15){1}{\line(1,0){0.48}}
\multiput(127.29,30.37)(0.48,-0.12){1}{\line(1,0){0.48}}
\multiput(127.77,30.25)(0.49,-0.1){1}{\line(1,0){0.49}}
\multiput(128.26,30.15)(0.49,-0.07){1}{\line(1,0){0.49}}
\multiput(128.76,30.08)(0.5,-0.05){1}{\line(1,0){0.5}}
\multiput(129.25,30.03)(0.5,-0.02){1}{\line(1,0){0.5}}
\put(129.75,30){\line(1,0){0.5}}
\multiput(130.25,30)(0.5,0.02){1}{\line(1,0){0.5}}
\multiput(130.75,30.03)(0.5,0.05){1}{\line(1,0){0.5}}
\multiput(131.24,30.08)(0.49,0.07){1}{\line(1,0){0.49}}
\multiput(131.74,30.15)(0.49,0.1){1}{\line(1,0){0.49}}
\multiput(132.23,30.25)(0.48,0.12){1}{\line(1,0){0.48}}
\multiput(132.71,30.37)(0.48,0.15){1}{\line(1,0){0.48}}
\multiput(133.18,30.52)(0.47,0.17){1}{\line(1,0){0.47}}
\multiput(133.65,30.69)(0.23,0.1){2}{\line(1,0){0.23}}
\multiput(134.11,30.88)(0.22,0.11){2}{\line(1,0){0.22}}
\multiput(134.56,31.1)(0.22,0.12){2}{\line(1,0){0.22}}
\multiput(135,31.34)(0.21,0.13){2}{\line(1,0){0.21}}
\multiput(135.43,31.6)(0.21,0.14){2}{\line(1,0){0.21}}
\multiput(135.84,31.88)(0.13,0.1){3}{\line(1,0){0.13}}
\multiput(136.23,32.18)(0.13,0.11){3}{\line(1,0){0.13}}
\multiput(136.62,32.5)(0.12,0.11){3}{\line(1,0){0.12}}
\multiput(136.98,32.84)(0.12,0.12){3}{\line(0,1){0.12}}
\multiput(137.33,33.2)(0.11,0.12){3}{\line(0,1){0.12}}
\multiput(137.66,33.57)(0.1,0.13){3}{\line(0,1){0.13}}
\multiput(137.97,33.96)(0.15,0.2){2}{\line(0,1){0.2}}
\multiput(138.26,34.37)(0.14,0.21){2}{\line(0,1){0.21}}
\multiput(138.53,34.79)(0.12,0.22){2}{\line(0,1){0.22}}
\multiput(138.78,35.22)(0.11,0.22){2}{\line(0,1){0.22}}
\multiput(139.01,35.66)(0.1,0.23){2}{\line(0,1){0.23}}
\multiput(139.21,36.12)(0.09,0.23){2}{\line(0,1){0.23}}
\multiput(139.4,36.58)(0.16,0.47){1}{\line(0,1){0.47}}
\multiput(139.56,37.05)(0.14,0.48){1}{\line(0,1){0.48}}
\multiput(139.69,37.53)(0.11,0.49){1}{\line(0,1){0.49}}
\multiput(139.8,38.02)(0.09,0.49){1}{\line(0,1){0.49}}
\multiput(139.89,38.51)(0.06,0.49){1}{\line(0,1){0.49}}
\multiput(139.95,39)(0.04,0.5){1}{\line(0,1){0.5}}
\multiput(139.99,39.5)(0.01,0.5){1}{\line(0,1){0.5}}

\linethickness{0.3mm}
\multiput(159.99,40.5)(0.01,-0.5){1}{\line(0,-1){0.5}}
\multiput(159.95,41)(0.04,-0.5){1}{\line(0,-1){0.5}}
\multiput(159.89,41.49)(0.06,-0.49){1}{\line(0,-1){0.49}}
\multiput(159.8,41.98)(0.09,-0.49){1}{\line(0,-1){0.49}}
\multiput(159.69,42.47)(0.11,-0.49){1}{\line(0,-1){0.49}}
\multiput(159.56,42.95)(0.14,-0.48){1}{\line(0,-1){0.48}}
\multiput(159.4,43.42)(0.16,-0.47){1}{\line(0,-1){0.47}}
\multiput(159.21,43.88)(0.09,-0.23){2}{\line(0,-1){0.23}}
\multiput(159.01,44.34)(0.1,-0.23){2}{\line(0,-1){0.23}}
\multiput(158.78,44.78)(0.11,-0.22){2}{\line(0,-1){0.22}}
\multiput(158.53,45.21)(0.12,-0.22){2}{\line(0,-1){0.22}}
\multiput(158.26,45.63)(0.14,-0.21){2}{\line(0,-1){0.21}}
\multiput(157.97,46.04)(0.15,-0.2){2}{\line(0,-1){0.2}}
\multiput(157.66,46.43)(0.1,-0.13){3}{\line(0,-1){0.13}}
\multiput(157.33,46.8)(0.11,-0.12){3}{\line(0,-1){0.12}}
\multiput(156.98,47.16)(0.12,-0.12){3}{\line(0,-1){0.12}}
\multiput(156.62,47.5)(0.12,-0.11){3}{\line(1,0){0.12}}
\multiput(156.23,47.82)(0.13,-0.11){3}{\line(1,0){0.13}}
\multiput(155.84,48.12)(0.13,-0.1){3}{\line(1,0){0.13}}
\multiput(155.43,48.4)(0.21,-0.14){2}{\line(1,0){0.21}}
\multiput(155,48.66)(0.21,-0.13){2}{\line(1,0){0.21}}
\multiput(154.56,48.9)(0.22,-0.12){2}{\line(1,0){0.22}}
\multiput(154.11,49.12)(0.22,-0.11){2}{\line(1,0){0.22}}
\multiput(153.65,49.31)(0.23,-0.1){2}{\line(1,0){0.23}}
\multiput(153.18,49.48)(0.47,-0.17){1}{\line(1,0){0.47}}
\multiput(152.71,49.63)(0.48,-0.15){1}{\line(1,0){0.48}}
\multiput(152.23,49.75)(0.48,-0.12){1}{\line(1,0){0.48}}
\multiput(151.74,49.85)(0.49,-0.1){1}{\line(1,0){0.49}}
\multiput(151.24,49.92)(0.49,-0.07){1}{\line(1,0){0.49}}
\multiput(150.75,49.97)(0.5,-0.05){1}{\line(1,0){0.5}}
\multiput(150.25,50)(0.5,-0.02){1}{\line(1,0){0.5}}
\put(149.75,50){\line(1,0){0.5}}
\multiput(149.25,49.97)(0.5,0.02){1}{\line(1,0){0.5}}
\multiput(148.76,49.92)(0.5,0.05){1}{\line(1,0){0.5}}
\multiput(148.26,49.85)(0.49,0.07){1}{\line(1,0){0.49}}
\multiput(147.77,49.75)(0.49,0.1){1}{\line(1,0){0.49}}
\multiput(147.29,49.63)(0.48,0.12){1}{\line(1,0){0.48}}
\multiput(146.82,49.48)(0.48,0.15){1}{\line(1,0){0.48}}
\multiput(146.35,49.31)(0.47,0.17){1}{\line(1,0){0.47}}
\multiput(145.89,49.12)(0.23,0.1){2}{\line(1,0){0.23}}
\multiput(145.44,48.9)(0.22,0.11){2}{\line(1,0){0.22}}
\multiput(145,48.66)(0.22,0.12){2}{\line(1,0){0.22}}
\multiput(144.57,48.4)(0.21,0.13){2}{\line(1,0){0.21}}
\multiput(144.16,48.12)(0.21,0.14){2}{\line(1,0){0.21}}
\multiput(143.77,47.82)(0.13,0.1){3}{\line(1,0){0.13}}
\multiput(143.38,47.5)(0.13,0.11){3}{\line(1,0){0.13}}
\multiput(143.02,47.16)(0.12,0.11){3}{\line(1,0){0.12}}
\multiput(142.67,46.8)(0.12,0.12){3}{\line(0,1){0.12}}
\multiput(142.34,46.43)(0.11,0.12){3}{\line(0,1){0.12}}
\multiput(142.03,46.04)(0.1,0.13){3}{\line(0,1){0.13}}
\multiput(141.74,45.63)(0.15,0.2){2}{\line(0,1){0.2}}
\multiput(141.47,45.21)(0.14,0.21){2}{\line(0,1){0.21}}
\multiput(141.22,44.78)(0.12,0.22){2}{\line(0,1){0.22}}
\multiput(140.99,44.34)(0.11,0.22){2}{\line(0,1){0.22}}
\multiput(140.79,43.88)(0.1,0.23){2}{\line(0,1){0.23}}
\multiput(140.6,43.42)(0.09,0.23){2}{\line(0,1){0.23}}
\multiput(140.44,42.95)(0.16,0.47){1}{\line(0,1){0.47}}
\multiput(140.31,42.47)(0.14,0.48){1}{\line(0,1){0.48}}
\multiput(140.2,41.98)(0.11,0.49){1}{\line(0,1){0.49}}
\multiput(140.11,41.49)(0.09,0.49){1}{\line(0,1){0.49}}
\multiput(140.05,41)(0.06,0.49){1}{\line(0,1){0.49}}
\multiput(140.01,40.5)(0.04,0.5){1}{\line(0,1){0.5}}
\multiput(140,40)(0.01,0.5){1}{\line(0,1){0.5}}

\linethickness{0.3mm}
\put(150,50){\line(0,1){10}}
\linethickness{0.3mm}
\put(160,20){\line(0,1){20}}
\put(8,65){\makebox(0,0)[cc]{$\Vec A$}}

\put(68,65){\makebox(0,0)[cc]{$\Vec A$}}

\put(118,65){\makebox(0,0)[cc]{$\Vec A$}}

\put(28,65){\makebox(0,0)[cc]{$\Omega^1 A$}}

\put(98,65){\makebox(0,0)[cc]{$\Omega^1 A$}}

\put(150,65){\makebox(0,0)[cc]{$\Omega^1 A$}}

\put(160,15){\makebox(0,0)[cc]{$\Omega^1 A$}}

\put(60,15){\makebox(0,0)[cc]{$\Omega^1 A$}}

\put(18,15){\makebox(0,0)[cc]{$\Omega^1 A$}}

\put(128,25){\makebox(0,0)[cc]{$\ev$}}

\put(90,25){\makebox(0,0)[cc]{$\ev$}}

\put(25,45){\makebox(0,0)[cc]{$\ev$}}

\put(75,55){\makebox(0,0)[cc]{$\square$}}

\put(142,55){\makebox(0,0)[cc]{$\square$}}

\end{picture}

The covariant derivative in (\ref{ghjksvdgwv}) is defined, using $\coev(1)=\alpha\tens w\in \Omega^1 A\tens_A \Vec A$ (summation implicit), as
\begin{eqnarray*}
\square(v)\ =\ \extd(v(\alpha))\tens w-(\ev\tens \id\tens\id)(v\tens\square(\alpha)\tens w)\ .
\end{eqnarray*}
To show that this works, first we need to note that it only depends on $\alpha\tens w\in \Omega^1 A\tens_A \Vec A$ (emphasis on $\tens_A$), as $v$ is a right $A$-module map. 
Then we need to use $a.\coev(1)=\coev(1).a\in 
\Omega^1 A\tens_A \Vec A$ for all $a\in A$, which follows from $\coev$ being a bimodule map.
Then $\square$ is a left bimodule covariant derivative, as, for $v\in \Vec A$,
\begin{eqnarray*}
\square(v.a) &=& \square(v).a+\sigma(v\tens\extd a)\ ,\cr
\square(a.v) &=& a.\square(v)+\extd a\tens v\ ,
\end{eqnarray*}
where $\sigma:\Vec A\tens_A \Omega^1 A\to\Omega^1  A\tens_A \Vec A$ is defined by
\begin{eqnarray}\label{sigrefdefdual}
\sigma\ =\ (\ev\tens\id\tens\id)(\id\tens\sigma^{-1}\tens\id)(\id\tens\id\tens\coev(1))\ .
\end{eqnarray}
In (\ref{sigrefdefdual}), note that in line with the comment on notation in the introduction, the $\sigma$ and $\sigma^{-1}$ are not inverses of each other, as they have different domains. 
Then (\ref{sigrefdefdual}) is equivalent to Fig.\ \ref{dszdfgzg}, by using the usual properties of evaluation and coevaluation: 

\unitlength 0.5 mm
\begin{picture}(120,65)(-20,20)  \label{fighhh}
\linethickness{0.3mm}
\multiput(10,70)(0.12,-0.12){42}{\line(1,0){0.12}}
\linethickness{0.3mm}
\multiput(25,55)(0.12,-0.12){42}{\line(1,0){0.12}}
\linethickness{0.3mm}
\multiput(10,50)(0.12,0.12){167}{\line(1,0){0.12}}
\linethickness{0.3mm}
\multiput(30,50)(0.01,-0.5){1}{\line(0,-1){0.5}}
\multiput(30.01,49.5)(0.04,-0.5){1}{\line(0,-1){0.5}}
\multiput(30.05,49)(0.06,-0.49){1}{\line(0,-1){0.49}}
\multiput(30.11,48.51)(0.09,-0.49){1}{\line(0,-1){0.49}}
\multiput(30.2,48.02)(0.11,-0.49){1}{\line(0,-1){0.49}}
\multiput(30.31,47.53)(0.14,-0.48){1}{\line(0,-1){0.48}}
\multiput(30.44,47.05)(0.16,-0.47){1}{\line(0,-1){0.47}}
\multiput(30.6,46.58)(0.09,-0.23){2}{\line(0,-1){0.23}}
\multiput(30.79,46.12)(0.1,-0.23){2}{\line(0,-1){0.23}}
\multiput(30.99,45.66)(0.11,-0.22){2}{\line(0,-1){0.22}}
\multiput(31.22,45.22)(0.12,-0.22){2}{\line(0,-1){0.22}}
\multiput(31.47,44.79)(0.14,-0.21){2}{\line(0,-1){0.21}}
\multiput(31.74,44.37)(0.15,-0.2){2}{\line(0,-1){0.2}}
\multiput(32.03,43.96)(0.1,-0.13){3}{\line(0,-1){0.13}}
\multiput(32.34,43.57)(0.11,-0.12){3}{\line(0,-1){0.12}}
\multiput(32.67,43.2)(0.12,-0.12){3}{\line(0,-1){0.12}}
\multiput(33.02,42.84)(0.12,-0.11){3}{\line(1,0){0.12}}
\multiput(33.38,42.5)(0.13,-0.11){3}{\line(1,0){0.13}}
\multiput(33.77,42.18)(0.13,-0.1){3}{\line(1,0){0.13}}
\multiput(34.16,41.88)(0.21,-0.14){2}{\line(1,0){0.21}}
\multiput(34.57,41.6)(0.21,-0.13){2}{\line(1,0){0.21}}
\multiput(35,41.34)(0.22,-0.12){2}{\line(1,0){0.22}}
\multiput(35.44,41.1)(0.22,-0.11){2}{\line(1,0){0.22}}
\multiput(35.89,40.88)(0.23,-0.1){2}{\line(1,0){0.23}}
\multiput(36.35,40.69)(0.47,-0.17){1}{\line(1,0){0.47}}
\multiput(36.82,40.52)(0.48,-0.15){1}{\line(1,0){0.48}}
\multiput(37.29,40.37)(0.48,-0.12){1}{\line(1,0){0.48}}
\multiput(37.77,40.25)(0.49,-0.1){1}{\line(1,0){0.49}}
\multiput(38.26,40.15)(0.49,-0.07){1}{\line(1,0){0.49}}
\multiput(38.76,40.08)(0.5,-0.05){1}{\line(1,0){0.5}}
\multiput(39.25,40.03)(0.5,-0.02){1}{\line(1,0){0.5}}
\put(39.75,40){\line(1,0){0.5}}
\multiput(40.25,40)(0.5,0.02){1}{\line(1,0){0.5}}
\multiput(40.75,40.03)(0.5,0.05){1}{\line(1,0){0.5}}
\multiput(41.24,40.08)(0.49,0.07){1}{\line(1,0){0.49}}
\multiput(41.74,40.15)(0.49,0.1){1}{\line(1,0){0.49}}
\multiput(42.23,40.25)(0.48,0.12){1}{\line(1,0){0.48}}
\multiput(42.71,40.37)(0.48,0.15){1}{\line(1,0){0.48}}
\multiput(43.18,40.52)(0.47,0.17){1}{\line(1,0){0.47}}
\multiput(43.65,40.69)(0.23,0.1){2}{\line(1,0){0.23}}
\multiput(44.11,40.88)(0.22,0.11){2}{\line(1,0){0.22}}
\multiput(44.56,41.1)(0.22,0.12){2}{\line(1,0){0.22}}
\multiput(45,41.34)(0.21,0.13){2}{\line(1,0){0.21}}
\multiput(45.43,41.6)(0.21,0.14){2}{\line(1,0){0.21}}
\multiput(45.84,41.88)(0.13,0.1){3}{\line(1,0){0.13}}
\multiput(46.23,42.18)(0.13,0.11){3}{\line(1,0){0.13}}
\multiput(46.62,42.5)(0.12,0.11){3}{\line(1,0){0.12}}
\multiput(46.98,42.84)(0.12,0.12){3}{\line(0,1){0.12}}
\multiput(47.33,43.2)(0.11,0.12){3}{\line(0,1){0.12}}
\multiput(47.66,43.57)(0.1,0.13){3}{\line(0,1){0.13}}
\multiput(47.97,43.96)(0.15,0.2){2}{\line(0,1){0.2}}
\multiput(48.26,44.37)(0.14,0.21){2}{\line(0,1){0.21}}
\multiput(48.53,44.79)(0.12,0.22){2}{\line(0,1){0.22}}
\multiput(48.78,45.22)(0.11,0.22){2}{\line(0,1){0.22}}
\multiput(49.01,45.66)(0.1,0.23){2}{\line(0,1){0.23}}
\multiput(49.21,46.12)(0.09,0.23){2}{\line(0,1){0.23}}
\multiput(49.4,46.58)(0.16,0.47){1}{\line(0,1){0.47}}
\multiput(49.56,47.05)(0.14,0.48){1}{\line(0,1){0.48}}
\multiput(49.69,47.53)(0.11,0.49){1}{\line(0,1){0.49}}
\multiput(49.8,48.02)(0.09,0.49){1}{\line(0,1){0.49}}
\multiput(49.89,48.51)(0.06,0.49){1}{\line(0,1){0.49}}
\multiput(49.95,49)(0.04,0.5){1}{\line(0,1){0.5}}
\multiput(49.99,49.5)(0.01,0.5){1}{\line(0,1){0.5}}

\linethickness{0.3mm}
\put(50,50){\line(0,1){20}}
\linethickness{0.3mm}
\put(10,30){\line(0,1){20}}
\put(65,55){\makebox(0,0)[cc]{$=$}}

\linethickness{0.3mm}
\put(80,50){\line(0,1){20}}
\linethickness{0.3mm}
\multiput(100,70)(0.12,-0.12){167}{\line(1,0){0.12}}
\linethickness{0.3mm}
\multiput(115,65)(0.12,0.12){42}{\line(1,0){0.12}}
\linethickness{0.3mm}
\multiput(100,50)(0.12,0.12){42}{\line(1,0){0.12}}
\linethickness{0.3mm}
\multiput(80,50)(0.01,-0.5){1}{\line(0,-1){0.5}}
\multiput(80.01,49.5)(0.04,-0.5){1}{\line(0,-1){0.5}}
\multiput(80.05,49)(0.06,-0.49){1}{\line(0,-1){0.49}}
\multiput(80.11,48.51)(0.09,-0.49){1}{\line(0,-1){0.49}}
\multiput(80.2,48.02)(0.11,-0.49){1}{\line(0,-1){0.49}}
\multiput(80.31,47.53)(0.14,-0.48){1}{\line(0,-1){0.48}}
\multiput(80.44,47.05)(0.16,-0.47){1}{\line(0,-1){0.47}}
\multiput(80.6,46.58)(0.09,-0.23){2}{\line(0,-1){0.23}}
\multiput(80.79,46.12)(0.1,-0.23){2}{\line(0,-1){0.23}}
\multiput(80.99,45.66)(0.11,-0.22){2}{\line(0,-1){0.22}}
\multiput(81.22,45.22)(0.12,-0.22){2}{\line(0,-1){0.22}}
\multiput(81.47,44.79)(0.14,-0.21){2}{\line(0,-1){0.21}}
\multiput(81.74,44.37)(0.15,-0.2){2}{\line(0,-1){0.2}}
\multiput(82.03,43.96)(0.1,-0.13){3}{\line(0,-1){0.13}}
\multiput(82.34,43.57)(0.11,-0.12){3}{\line(0,-1){0.12}}
\multiput(82.67,43.2)(0.12,-0.12){3}{\line(0,-1){0.12}}
\multiput(83.02,42.84)(0.12,-0.11){3}{\line(1,0){0.12}}
\multiput(83.38,42.5)(0.13,-0.11){3}{\line(1,0){0.13}}
\multiput(83.77,42.18)(0.13,-0.1){3}{\line(1,0){0.13}}
\multiput(84.16,41.88)(0.21,-0.14){2}{\line(1,0){0.21}}
\multiput(84.57,41.6)(0.21,-0.13){2}{\line(1,0){0.21}}
\multiput(85,41.34)(0.22,-0.12){2}{\line(1,0){0.22}}
\multiput(85.44,41.1)(0.22,-0.11){2}{\line(1,0){0.22}}
\multiput(85.89,40.88)(0.23,-0.1){2}{\line(1,0){0.23}}
\multiput(86.35,40.69)(0.47,-0.17){1}{\line(1,0){0.47}}
\multiput(86.82,40.52)(0.48,-0.15){1}{\line(1,0){0.48}}
\multiput(87.29,40.37)(0.48,-0.12){1}{\line(1,0){0.48}}
\multiput(87.77,40.25)(0.49,-0.1){1}{\line(1,0){0.49}}
\multiput(88.26,40.15)(0.49,-0.07){1}{\line(1,0){0.49}}
\multiput(88.76,40.08)(0.5,-0.05){1}{\line(1,0){0.5}}
\multiput(89.25,40.03)(0.5,-0.02){1}{\line(1,0){0.5}}
\put(89.75,40){\line(1,0){0.5}}
\multiput(90.25,40)(0.5,0.02){1}{\line(1,0){0.5}}
\multiput(90.75,40.03)(0.5,0.05){1}{\line(1,0){0.5}}
\multiput(91.24,40.08)(0.49,0.07){1}{\line(1,0){0.49}}
\multiput(91.74,40.15)(0.49,0.1){1}{\line(1,0){0.49}}
\multiput(92.23,40.25)(0.48,0.12){1}{\line(1,0){0.48}}
\multiput(92.71,40.37)(0.48,0.15){1}{\line(1,0){0.48}}
\multiput(93.18,40.52)(0.47,0.17){1}{\line(1,0){0.47}}
\multiput(93.65,40.69)(0.23,0.1){2}{\line(1,0){0.23}}
\multiput(94.11,40.88)(0.22,0.11){2}{\line(1,0){0.22}}
\multiput(94.56,41.1)(0.22,0.12){2}{\line(1,0){0.22}}
\multiput(95,41.34)(0.21,0.13){2}{\line(1,0){0.21}}
\multiput(95.43,41.6)(0.21,0.14){2}{\line(1,0){0.21}}
\multiput(95.84,41.88)(0.13,0.1){3}{\line(1,0){0.13}}
\multiput(96.23,42.18)(0.13,0.11){3}{\line(1,0){0.13}}
\multiput(96.62,42.5)(0.12,0.11){3}{\line(1,0){0.12}}
\multiput(96.98,42.84)(0.12,0.12){3}{\line(0,1){0.12}}
\multiput(97.33,43.2)(0.11,0.12){3}{\line(0,1){0.12}}
\multiput(97.66,43.57)(0.1,0.13){3}{\line(0,1){0.13}}
\multiput(97.97,43.96)(0.15,0.2){2}{\line(0,1){0.2}}
\multiput(98.26,44.37)(0.14,0.21){2}{\line(0,1){0.21}}
\multiput(98.53,44.79)(0.12,0.22){2}{\line(0,1){0.22}}
\multiput(98.78,45.22)(0.11,0.22){2}{\line(0,1){0.22}}
\multiput(99.01,45.66)(0.1,0.23){2}{\line(0,1){0.23}}
\multiput(99.21,46.12)(0.09,0.23){2}{\line(0,1){0.23}}
\multiput(99.4,46.58)(0.16,0.47){1}{\line(0,1){0.47}}
\multiput(99.56,47.05)(0.14,0.48){1}{\line(0,1){0.48}}
\multiput(99.69,47.53)(0.11,0.49){1}{\line(0,1){0.49}}
\multiput(99.8,48.02)(0.09,0.49){1}{\line(0,1){0.49}}
\multiput(99.89,48.51)(0.06,0.49){1}{\line(0,1){0.49}}
\multiput(99.95,49)(0.04,0.5){1}{\line(0,1){0.5}}
\multiput(99.99,49.5)(0.01,0.5){1}{\line(0,1){0.5}}

\linethickness{0.3mm}
\put(120,30){\line(0,1){20}}
\put(5,75){\makebox(0,0)[cc]{$\Vec A$}}

\put(80,75){\makebox(0,0)[cc]{$\Vec A$}}

\put(28,75){\makebox(0,0)[cc]{$\Omega^1 A$}}

%%%%%%%%%%%%%%%
\refstepcounter{piccie} \label{dszdfgzg}
\put(188,45){\makebox(0,0)[cc]{Fig.\ \arabic{piccie}}}
%%%%%%%%%%%%%%%

\put(50,75){\makebox(0,0)[cc]{$\Omega^1 A$}}

\put(100,75){\makebox(0,0)[cc]{$\Omega^1 A$}}

\put(118,75){\makebox(0,0)[cc]{$\Omega^1 A$}}

\put(118,26){\makebox(0,0)[cc]{$\Omega^1 A$}}

\put(8,26){\makebox(0,0)[cc]{$\Omega^1 A$}}

\put(40,35){\makebox(0,0)[cc]{$\ev$}}

\put(90,35){\makebox(0,0)[cc]{$\ev$}}

\put(30,60){\makebox(0,0)[cc]{$\sigma$}}

\put(122,60){\makebox(0,0)[cc]{$\sigma^{-1}$}}

\end{picture}

Now $\square$ extends to a left bimodule covariant derivative $\square^{\<n\>}:\Vec^{\tens n}A\to \Omega^{1}A\tens_A \Vec^{\tens n}A$, defined recursively by
\begin{eqnarray*}
\square^{\<0\>} &=& \extd\ ,\quad \square^{\<1\>}\ =\ \square\ ,\cr
\square^{\<n+1\>} &=& \square\tens \id^{\tens n}+
(\sigma\tens\id^{\tens n})(\id^{\tens 1}\tens\square^{\<n\>})\ .
\end{eqnarray*}

\begin{propos}   \label{kjhcvjhhvjh}
For all $n\ge 1$ we have
\begin{eqnarray*}
\extd\, \ev^{\<n\>}\ =\ (\id\tens \ev^{\<n\>})(\square^{\<n\>}\tens\id)+(\ev^{\<n\>}\tens\id)(\id\tens\square^{\<n\>}):\Vec^{\tens n}A\tens_A \Omega^{\tens n}A\to \Omega^1 A\ .
\end{eqnarray*}
\end{propos}
\noindent {\bf Proof:}\quad The $n=1$ case is true by definition of $\square$ on $\Vec A$. Now suppose that the statement is true for $n$, and take
\begin{eqnarray*}
v\tens\underline w\tens\underline\beta\tens\alpha\in \Vec A\tens_A\Vec^{\tens n}A\tens_A \Omega^{\tens n}A\tens_A \Omega^1 A\ .
\end{eqnarray*}
Now calculate
\begin{eqnarray*}
&& (\id\tens \ev^{\<n+1\>})(\square^{\<n+1\>}\tens\id^{\tens n+1})(v\tens\underline w\tens\underline\beta\tens\alpha)  \cr
&=&  (\id\tens \ev^{\<n+1\>})\big(\square v\tens\underline w\tens\underline\beta\tens\alpha
+ (\sigma\tens\id^{\tens n})(v\tens \square^{\<n\>}\underline w)\tens\underline\beta\tens\alpha\big) \cr
&=&  (\id\tens \ev)\big(\square v\tens\ev^{\<n\>}(\underline w\tens\underline\beta).\alpha+
 (\sigma\tens\id)(v\tens (\id\tens  \ev^{\<n\>}  )(\square^{\<n\>}\underline w\tens\underline\beta)\tens\alpha\big)\ ,
\end{eqnarray*}
and since from  (\ref{sigrefdefdual}),
\begin{eqnarray*}
(\id\tens\ev)(\sigma\tens\id)\ =\ (\ev\tens\id)(\id\tens\sigma^{-1}):\Vec A\tens_A \Omega^1 A  \tens_A\Omega^1 A \to \Omega^1 A \ ,
\end{eqnarray*}
we have
\begin{eqnarray*}
&& (\id\tens \ev^{\<n+1\>})(\square^{\<n+1\>}\tens\id^{\tens n+1})(v\tens\underline w\tens\underline\beta\tens\alpha)  \cr
&=&  (\id\tens \ev)\big(\square v\tens\ev^{\<n\>}(\underline w\tens\underline\beta).\alpha\big) \cr
& & +\
(\ev\tens\id)(\id\tens\sigma^{-1})(v\tens (\id\tens  \ev^{\<n\>}  )(\square^{\<n\>}\underline w\tens\underline\beta)\tens\alpha\big)\ .
\end{eqnarray*}
Also calculate
\begin{eqnarray*}
&&(\ev^{\<n+1\>}\tens\id)(\id^{\tens n+1}\tens\square^{\<n+1\>})(v\tens\underline w\tens\underline\beta\tens\alpha)  \cr
&=& (\ev^{\<n+1\>}\tens\id)\big(  v\tens\underline w\tens\underline\beta\tens\square\alpha+
 v\tens\underline w\tens(\id^{\tens n}\tens\sigma^{-1})(\square^{\<n\>}\underline\beta\tens\alpha)\big)
 \cr
 &=&  (\ev\tens\id)\big(  v\tens\ev^{\<n\>}(\underline w\tens\underline\beta).\square\alpha+
(\id\tens\sigma^{-1})( v\tens(\ev^{\<n\>}\tens\id)(\underline w\tens\square^{\<n\>}\underline\beta)\tens\alpha)\big)\ .
\end{eqnarray*}
Now, using our assumption,
\begin{eqnarray*}
&&\big( (\id\tens \ev^{\<n+1\>})(\square^{\<n+1\>}\tens\id^{\tens n+1})
+(\ev^{\<n+1\>}\tens\id)(\id^{\tens n+1}\tens\square^{\<n+1\>})\big)(v\tens\underline w\tens\underline\beta\tens\alpha)  \cr
&=&  (\id\tens \ev)\big(\square v\tens\ev^{\<n\>}(\underline w\tens\underline\beta).\alpha\big) 
+ (\ev\tens\id)\big(  v\tens\ev^{\<n\>}(\underline w\tens\underline\beta).\square\alpha\big)\cr
& & +\
(\ev\tens\id)(\id\tens\sigma^{-1})(v\tens \Big((\id\tens  \ev^{\<n\>}  )(\square^{\<n\>}\underline w\tens\underline\beta)+(\ev^{\<n\>}\tens\id)(\underline w\tens\square^{\<n\>}\underline\beta)\Big)\tens\alpha\big)\cr
&=&  (\id\tens \ev)\big(\square v\tens\ev^{\<n\>}(\underline w\tens\underline\beta).\alpha\big) 
+ (\ev\tens\id)\big(  v\tens\ev^{\<n\>}(\underline w\tens\underline\beta).\square\alpha\big)\cr
& & +\
(\ev\tens\id)(\id\tens\sigma^{-1})(v\tens\extd\, \ev^{\<n\>}(\underline w\tens\underline\beta)\tens\alpha\big)\cr
&=&  (\id\tens \ev)\big(\square v\tens\ev^{\<n\>}(\underline w\tens\underline\beta).\alpha\big) 
+ (\ev\tens\id)\big(  v\tens \square(\ev^{\<n\>}(\underline w\tens\underline\beta).\alpha)\big)\cr
&=& \extd\,\ev(v\tens \ev^{\<n\>}(\underline w\tens\underline\beta).\alpha) \cr
&=& \extd\, \ev^{\<n+1\>}(v\tens\underline w\tens\underline\beta\tens\alpha)\ ,
\end{eqnarray*}
which completes the proof by induction.\quad$\square$

\section{An action of tensor products of vector fields} \label{adruauifvuh}
We briefly remind the reader of some notation. 
In Section \ref{vcgsajhdx} we introduced $\square:\Omega^1 A\to \Omega^1 A\tens_A \Omega^1 A$ as a right bimodule covariant derivative, with generalised braiding $\sigma^{-1}:\Omega^1 A\tens_A \Omega^1 A\to \Omega^1 A\tens_A \Omega^1 A$ invertible (the inverse being denoted $\sigma$), and its extension to tensor products
 $\square^{\<n\>}:\Omega^{\tens n}A\to \Omega^{\tens n+1}A$. 
In Section \ref{gcysiuaf} we introduced 
 a dual left covariant derivative $\square:\Vec A\to \Omega^1 A\tens_A \Vec A$
 and its extension
$\square^{\<n\>}:\Vec^{\tens n}A\to \Omega^{1}A\tens_A \Vec^{\tens n}A$.

Suppose that $E$ is a left $A$ module, with a left covariant derivative $\nabla:E\to\Omega^1 A\tens_A E$. We iterate this to define $\nabla^{(n)}:E\to \Omega^{\tens n}A\tens_A E$
recursively by
\begin{eqnarray}\label{jktycfxtu}
\nabla^{(1)}\ =\ \nabla\ ,\quad \nabla^{(n+1)}=(\square^{\<n\>}\tens\id_E+\id^{\tens n}\tens\nabla)\,
\nabla^{(n)}\ .
\end{eqnarray}
To do this, we need to check that $\square^{\<n\>}\tens\id_E+\id^{\tens n}\tens\nabla$ is a well defined operation on $\Omega^{\tens n}A\tens_A E$, which we do as follows, for $a\in A$, $\underline\xi\in \Omega^{\tens n}A$ and $e\in E$:
\begin{eqnarray*}
(\square^{\<n\>}\tens\id_E+\id^{\tens n}\tens\nabla)(\underline\xi.a\tens e) &=& 
\square^{\<n\>}(\underline\xi).a\tens e+ \underline\xi\tens\extd a\tens e +
\underline\xi.a\tens \nabla e\ ,\cr
(\square^{\<n\>}\tens\id_E+\id^{\tens n}\tens\nabla)(\underline\xi\tens a.e) &=& 
\square^{\<n\>}(\underline\xi)\tens a.e+ \underline\xi\tens\extd a\tens e +
\underline\xi\tens a.\nabla e\ .
\end{eqnarray*}
Now we can define an `action' of $\underline v\in \Vec^{\tens n}A$ on $e\in E$ by
\begin{eqnarray}  \label{bcahkjlclkv}
\underline v\,\la\, e \ =\ (\ev^{\<n\>}\tens\id_E)\, (\underline v\tens\nabla^{(n)}e)\ .
\end{eqnarray}
The quotation marks in `action' are due to the fact that we have not yet given the product
with respect to which it is an action. This will be given in the forthcoming Theorem~\ref{kuykjvcvu}. It would, however, be slightly dishonest to pretend that Theorem~\ref{kuykjvcvu} was obtained by a claim that the product was in any sense obvious, rather than by saying that the action is the obvious quantity. The following Lemma \ref{vcjhckvjvhjvj}
about the action will be useful in several places.

\begin{lemma}\label{vcjhckvjvhjvj}
For the operation $\la$ in (\ref{bcahkjlclkv}), with $\underline v\in \Vec^{\tens n}A$, $w\in \Vec A$
and $e\in E$,
\begin{eqnarray*}
w\,\la(\underline v\,\la\, e) \ =\ (w\tens \underline v)\la\, e + \big((\ev\tens\id^{\tens n})(w\tens 
\square^{\<n\>}\underline v)\big)\la\,e\ .
\end{eqnarray*}
\end{lemma}
\noindent {\bf Proof:}\quad Begin by using Prop.\ \ref{kjhcvjhhvjh} on the definition of $\underline v\,\la\, e$ as follows;
\begin{eqnarray*}
\nabla(\underline v\,\la\, e) &=& \nabla\, (\ev^{\<n\>}\tens\id_E)\, (\underline v\tens\nabla^{(n)}e)\cr
&=&  (\extd\,\ev^{\<n\>}\tens\id_E+\ev^{\<n\>}\tens\nabla_E)\, (\underline v\tens\nabla^{(n)}e) \cr
&=& (\id^{\tens 1}\tens \ev^{\<n\>}\tens\id_E)(\square^{\<n\>}\underline v\tens\nabla^{(n)}e)  \cr
&& +\ (\ev^{\<n\>}\tens\id^{\tens 1}\tens\id_E)(\underline v\tens(\square^{\<n\>}\tens\id_E+\id^{\tens n}\tens\nabla)\nabla^{(n)}e)  \cr
&=&  (\id^{\tens 1}\tens \ev^{\<n\>}\tens\id_E)(\square^{\<n\>}\underline v\tens\nabla^{(n)}e)  \cr
&& +\ (\ev^{\<n\>}\tens\id^{\tens 1}\tens\id_E)(\underline v\tens\nabla^{(n+1)}e)  \ .
\end{eqnarray*}
Now we have
\begin{eqnarray*}
w\,\la(\underline v\,\la\, e) &=& (\ev\tens\id_E)(w\tens \nabla(\underline v\la\, e))  \cr
&=& (\ev\tens \ev^{\<n\>}\tens\id_E)(w\tens\square^{\<n\>}\underline v\tens\nabla^{(n)}e)  \cr
&&+\ (\ev^{\<n+1\>}\tens\id_E)(w\tens\underline v\tens\nabla^{(n+1)}e) \ .\quad\square
\end{eqnarray*}

\begin{lemma} \label{hjsvhjkvjhk}
For all $k\ge 0$, the following recursive procedure
gives a well defined function  $\bullet_k:\Vec^{\tens n}A\tens \Vec^{\tens m}A\to \Vec^{\tens k}A$
satisfying $(a.\underline v)\bullet_k\underline w=a.(\underline v\bullet_k\underline w)$, for all $a\in A$. The definition is recursive in $n\ge 0$: The starting cases are (for $u\in \Vec A$
and $\underline w\in\Vec^{\tens m}A$)
\begin{eqnarray*}
n=0\ , & a\bullet_k \underline w &= \left\{\begin{array}{cc}a.\underline w & k=m \\0 & k\neq m\end{array}\right. \cr
n=1\ , & u\bullet_k \underline w &= \left\{\begin{array}{cc}u\tens w & k=m+1 \\
(\ev\tens\id^{\tens m})(u\tens\square^{\<m\>}\underline w) & k= m
\\0 & \mathrm{otherwise}\end{array}\right. \ .
\end{eqnarray*}
The definition continues with, for $\underline v\in \Vec^{\tens n}A$
(setting $\bullet_{-1}$ to be zero),
\begin{eqnarray*}
(u\tens\underline v)\bullet_k\underline w &=& u\tens (\underline v\bullet_{k-1}\underline w)
+u\bullet_k(\underline v\bullet_k\underline w)
-(u\bullet_n \underline v)\bullet_k\underline w\ .
\end{eqnarray*}
\end{lemma}
\noindent {\bf Proof:}\quad We prove that (a) $\bullet_k:\Vec^{\tens n}A\tens \Vec^{\tens m}A\to \Vec^{\tens k}A$ is well defined, and that (b) $(a.\underline v)\bullet_k\underline w=a.(\underline v\bullet_k\underline w)$ for all $a\in A$, by induction on $n\ge 0$. To begin, the explicit formulae for 
$n=0$ and $n=1$ have these properties. Now assume that (a) and (b) hold for $n$, and examine the $n+1$ case. 

For (a) we note that, \textit{a priori}, the recursive definition in
the statement only defines 
\begin{eqnarray*}
\bullet_k:(\Vec A\tens\Vec^{\tens n}A)\tens \Vec^{\tens m}A\to \Vec^{\tens k}A\ .
\end{eqnarray*}
Instead of $\Vec A\tens\Vec^{\tens n}A$ we want $\Vec A\tens_A\Vec^{\tens n}A$, and that necessitates the following check: Continuing the notation of the statement, we require the equality of the following two quantities in $\Vec^{\tens k}A$,
\begin{eqnarray}  \label{vfjkcvckjcxf}
(u.a\tens\underline v)\bullet_k\underline w &=& u.a\tens (\underline v\bullet_{k-1}\underline w)
+(u.a)\bullet_k(\underline v\bullet_k\underline w)
-((u.a)\bullet_n \underline v)\bullet_k\underline w\ ,\cr
(u\tens a.\underline v)\bullet_k\underline w &=& u\tens ((a.\underline v)\bullet_{k-1}\underline w)
+u\bullet_k((a.\underline v)\bullet_k\underline w)
-(u\bullet_n (a.\underline v))\bullet_k\underline w\ .
\end{eqnarray}
Using our assumption of (b) for $n$, we have
\begin{eqnarray}  \label{cvusavc}
(u\tens a.\underline v)\bullet_k\underline w &=& u\tens a.(\underline v\bullet_{k-1}\underline w)
+u\bullet_k(a.(\underline v\bullet_k\underline w))
-(u\bullet_n (a.\underline v))\bullet_k\underline w\ .
\end{eqnarray}
From the statement, for $\underline x\in \Vec^{\tens s}A$,
\begin{eqnarray*}
(u.a)\bullet_s\underline x \ =\ u\bullet_s(a.\underline x)-\ev(u\tens\extd a).\underline x\ ,
\end{eqnarray*}
and using this twice (for $s=k$ and $s=n$) on (\ref{cvusavc}) gives
\begin{eqnarray}
(u\tens a.\underline v)\bullet_k\underline w &=& u.a\tens (\underline v\bullet_{k-1}\underline w)
+(u.a)\bullet_k(\underline v\bullet_k\underline w)
-((u.a)\bullet_n \underline v)\bullet_k\underline w\cr
&& +\ \ev(u\tens\extd a).(\underline v\bullet_k\underline w) - (\ev(u\tens\extd a).\underline v)\bullet_k\underline w\ ,
\end{eqnarray}
and by the assumption (b)
again we have the result that both quantities in (\ref{vfjkcvckjcxf}) are identical.

For (b) we have, using the inductive hypothesis,
\begin{eqnarray*}
(a.u\tens\underline v)\bullet_k\underline w &=& a.u\tens (\underline v\bullet_{k-1}\underline w)
+(a.u)\bullet_k(\underline v\bullet_k\underline w)
-((a.u)\bullet_n \underline v)\bullet_k\underline w\cr
&=& a.u\tens (\underline v\bullet_{k-1}\underline w)
+a.(u\bullet_k(\underline v\bullet_k\underline w))
-(a.(u\bullet_n \underline v))\bullet_k\underline w\cr
&=& a.\big((u\tens\underline v)\bullet_k\underline w\big)\ .\quad\square
\end{eqnarray*}

\medskip The motivation behind Lemma \ref{hjsvhjkvjhk} is given by the following result.

\begin{propos}\label{kuyfvjgccttcd} For all $\underline v\in \Vec^{\tens n}A$ and 
$\underline w\in \Vec^{\tens m}A$,
\begin{eqnarray*}
\underline v\,\la(\underline w\,\la\,e) \ =\ \sum_{k\ge 0} ( \underline v\bullet_k\underline w)\la\,e\ .
\end{eqnarray*}
Note that the indices in the summation are finite, as $\underline v\bullet_k\underline w=0$ for 
$k> n+m$. 
\end{propos}
\noindent {\bf Proof:}\quad This is proved by induction on $n\ge 0$. The $n=0$ case is automatic from the definition, and the $n=1$ case is a combination of Lemma \ref{vcjhckvjvhjvj}
and Lemma \ref{hjsvhjkvjhk}. Now suppose that the result is true for some $n\ge 1$, and consider
\begin{eqnarray*}
(u\tens\underline v)\,\la(\underline w\,\la\,e)\ , 
\end{eqnarray*}
where $u\in \Vec A$ and $\underline v\in \Vec^{\tens n} A$. By Lemma \ref{vcjhckvjvhjvj}
we can write
\begin{eqnarray*}
(u\tens\underline v)\,\la(\underline w\,\la\,e) \ =\ u\,\la(\underline v\,\la(\underline w\,\la\,e))-
\big( (\ev\tens\id^{\tens n})(u\tens 
\square^{\<n\>}\underline v)  \big)\la(\underline w\,\la\,e)\ .
\end{eqnarray*}
Now we can use our inductive assumption to write this as
\begin{eqnarray*}
(u\tens\underline v)\,\la(\underline w\,\la\,e) &=& u\,\la\sum_{k\ge 0} 
(\underline v\bullet_k\underline w)\,\la\,e-\sum_{k\ge 0} 
\Big(\big( (\ev\tens\id^{\tens n})(u\tens 
\square^{\<n\>}\underline v)  \big)\bullet_k\underline w\Big)\,\la\,e\ ,
\end{eqnarray*}
and using Lemma \ref{vcjhckvjvhjvj} again gives the result (with a shift in summation index)
\begin{eqnarray*}
(u\tens\underline v)\,\la(\underline w\,\la\,e) &=& \sum_{k\ge 0} \Big(
u\tens(\underline v\bullet_k\underline w)
+u\bullet_k(\underline v\bullet_k\underline w)
-(u\bullet_n \underline v)\bullet_k\underline w\Big)\,\la\,e\ .\quad\square
\end{eqnarray*}

\medskip Proposition \ref{kuyfvjgccttcd} suggests that the following Proposition~\ref{bdhscvkvhj} might be true, but an explicit check is necessary.

\begin{propos}\label{bdhscvkvhj}
On the free tensor algebra
\begin{eqnarray*}
\mathcal{T}\,\Vec A\ =\ \bigoplus_{n\ge 0} \Vec^{\tens n}A\ ,
\end{eqnarray*}
the operation $\bullet:\mathcal{T}\,\Vec A\tens \mathcal{T}\,\Vec A\to \mathcal{T}\,\Vec A$
defined by 
\begin{eqnarray*}
\underline v\bullet\underline w \ =\  \sum_{k\ge 0} \underline v\bullet_k\underline w
\end{eqnarray*}
makes $\mathcal{T}\,\Vec A$ into an associative algebra. \end{propos}
\noindent {\bf Proof:}\quad We prove that, for all $\underline w,\underline x\in \mathcal{T}\,\Vec A$,
\begin{eqnarray*}
\underline v\bullet (\underline w\bullet \underline x)\ =\ (\underline v\bullet \underline w)\bullet \underline x \ ,
\end{eqnarray*}
for $\underline v\in \Vec^{\tens n}A$ by induction on $n$. The $n=0$ case is just the left $A$ linearity property which has already been noted in Lemma \ref{hjsvhjkvjhk}. Now suppose that the result is true for some $n\ge 0$. From Lemma \ref{hjsvhjkvjhk}, for $\underline v\in \Vec^{\tens n}A$ and $u\in \Vec A$,
\begin{eqnarray*}
(u\tens\underline v)\bullet_k (\underline w\bullet_s \underline x) &=&
u\tens(\underline v\bullet_{k-1} (\underline w\bullet_s \underline x))+u\bullet_k
(\underline v\bullet_{k} (\underline w\bullet_s \underline x))
- (u\bullet_n\underline v)\bullet_k (\underline w\bullet_s \underline x)\ ,
\end{eqnarray*}
and on summing over $k$ and $s$, we find
\begin{eqnarray}  \label{auigcuagvk}
(u\tens\underline v)\bullet (\underline w\bullet \underline x) \,=\,
u\tens(\underline v\bullet (\underline w\bullet \underline x))
- (u\bullet_n\underline v)\bullet (\underline w\bullet \underline x)
+\sum_k u\bullet_k(\underline v\bullet_{k} (\underline w\bullet \underline x))\ .
\end{eqnarray}
Applying Lemma \ref{hjsvhjkvjhk} twice,
\begin{eqnarray*}
((u\tens\underline v)\bullet_k \underline w)\bullet_s \underline x &=&
\Big( u\tens(\underline v\bullet_{k-1} \underline w)+u\bullet_k(\underline v\bullet_{k} \underline w)
-(u\bullet_n\underline v)\bullet_{k} \underline w
\Big)\bullet_s \underline x   \cr
&=&  u\tens((\underline v\bullet_{k-1} \underline w)\bullet_{s-1} \underline x)
+ u\bullet_s((\underline v\bullet_{k-1} \underline w)\bullet_{s} \underline x) \cr
&& -\  (u\bullet_{k-1}(\underline v\bullet_{k-1} \underline w))\bullet_{s} \underline x
\cr
&& +\, \Big( u\bullet_k(\underline v\bullet_{k} \underline w)
-(u\bullet_n\underline v)\bullet_{k} \underline w
\Big)\bullet_s \underline x \ ,
\end{eqnarray*}
and on summing over $k$ and $s$, we find
\begin{eqnarray*}
((u\tens\underline v)\bullet \underline w)\bullet \underline x 
&=&  u\tens((\underline v\bullet \underline w)\bullet \underline x)
+ \sum_s u\bullet_s((\underline v\bullet \underline w)\bullet_{s} \underline x) \cr
&& -\   \sum_k(u\bullet_{k-1}(\underline v\bullet_{k-1} \underline w))\bullet \underline x
\cr
&& +\, \sum_k\big( u\bullet_k(\underline v\bullet_{k} \underline w)
\big)\bullet \underline x 
-  \big( (u\bullet_n\underline v)\bullet \underline w
\big)\bullet \underline x \cr
&=&  u\tens((\underline v\bullet \underline w)\bullet \underline x)
+ \sum_s u\bullet_s((\underline v\bullet \underline w)\bullet_{s} \underline x)
-  \big( (u\bullet_n\underline v)\bullet \underline w
\big)\bullet \underline x \ .
\end{eqnarray*}
By using the inductive assumption, we can rewrite this as
\begin{eqnarray*} \label{auigcuagvk2}
((u\tens\underline v)\bullet \underline w)\bullet \underline x 
&=&  u\tens(\underline v\bullet (\underline w\bullet \underline x))
+ \sum_s u\bullet_s((\underline v\bullet \underline w)\bullet_{s} \underline x)
-   (u\bullet_n\underline v)\bullet (\underline w \bullet \underline x)\ ,
\end{eqnarray*}
and combining this with (\ref{auigcuagvk}) gives (relabeling one index from $k$ to $s$)
\begin{eqnarray}
&& ((u\tens\underline v)\bullet \underline w)\bullet \underline x - (u\tens\underline v)\bullet (\underline w\bullet \underline x) \cr
&=&   \sum_s u\bullet_s((\underline v\bullet \underline w)\bullet_{s} \underline x)
- \sum_s u\bullet_s(\underline v\bullet_{s} (\underline w\bullet \underline x))\ .
\end{eqnarray}
The degree $s$ part of $\underline v\bullet (\underline w\bullet \underline x)=
(\underline v\bullet \underline w)\bullet \underline x$ gives
the following, completing the proof of associativity:
\begin{eqnarray*}
\underline v\bullet_s (\underline w\bullet \underline x)\ =\ 
(\underline v\bullet \underline w)\bullet_{s} \underline x\ .\quad\square
\end{eqnarray*}

\medskip We should stress that
the associative multiplication $\bullet$ in Proposition
\ref{bdhscvkvhj} is not defined on $\mathcal{T}\,\Vec A\tens_A \mathcal{T}\,\Vec A$, but
just on $\mathcal{T}\,\Vec A\tens_{\mathbb{C}} \mathcal{T}\,\Vec A$. This is not surprising, as
in repeated application of vector fields, the vector fields will themselves become differentiated. However the reader should remember that this problem arises with the usual bimodule structure on the tensor product, which is given, for example, by
$(u\tens v).a=u\tens (v.a)$ and $a.(u\tens v)=(a.u)\tens v$. The fact that $\bullet$ is associative gives alternative $A$-bimodule actions. We define $\mathcal{T}\,\Vec A_\bullet$ to be the same as $\mathcal{T}\,\Vec A$ as a left $A$-module, but with right module structure given by $\bullet$. We have no need to modify the 
left action, as the usual action is the same as $\bullet$ in this case. Now we state the main result in this section. 

\begin{theorem}  \label{kuykjvcvu}
The $A$-bimodule $\mathcal{T}\,\Vec A_\bullet$ with
product $\bullet:\mathcal{T}\,\Vec A_\bullet\tens_A \mathcal{T}\,\Vec A_\bullet\to \mathcal{T}\,\Vec A_\bullet$
defined by 
\begin{eqnarray*}
\underline v\bullet\underline w \ =\  \sum_{k\ge 0} \underline v\bullet_k\underline w
\end{eqnarray*}
is an associative algebra, with unit $1\in \Vec^{\tens 0}A=A$. Further, for a left $A$-module $E$ with left covariant derivative $\nabla$, the map 
in (\ref{bcahkjlclkv}) gives $\la:\mathcal{T}\,\Vec A_\bullet\tens_A E\to E$
which is an action of this algebra. 
\end{theorem}
\noindent {\bf Proof:}\quad Combining the results in this section.\quad$\square$

\begin{remark}  \label{vbsdhkvik}
 We need to say something about the category ${}_A\mathcal{E}$ introduced in Definition \ref{cvgjfxtzhc}. For a morphism $T:(E,\nabla_E)\to (F,\nabla_F)$
it is easy to show by induction that $(\id\tens T)\nabla_E^{(n)}=\nabla_F^{(n)}\,T$. From this we have $\underline{v}\,\la\,T(e)=T(\underline{v}\,\la\,e)$ for all $e\in E$ and all $\underline{v}\in \mathcal{T}\,\Vec A$. 
\end{remark}

\section{Noncommutative Sobolev spaces} \label{lkjsbvkcyxth}
In the classical theory of elliptic differential operators, a vital part is played by Sobolev spaces
\cite{Sobolev1,Sobolev2,Sobolev3}. 
They are also used in the analytic theory of complex manifolds
\cite{GrifHar}. 
The space $W^{k,p}$ is used to denote the functions which have $L^p$ norm of the derivatives up to order $k$ (defined as a completion of smooth functions of compact support). This is normally defined locally, but in noncommutative geometry we are forced to use a global definition. To do this we use a Hermitian inner product on $\Omega^1 A$. Given the machinery that we have in place already, it is convenient to define the Sobolev space for left modules $E$ with covariant derivative, and then specialise to the case $E=A$ if required. We shall only cover the $p=2$ case, but this is the most useful case, as it gives a Hilbert space. 

We must first define inner products on $A$-bimodules, where we assume that $A$ is a star algebra. In fact, we shall assume that $A$ is a subalgebra of a $C^*$ algebra, so that we have the usual ideas of positivity. (The reader should think of the algebra of smooth functions on a compact manifold being a subset of the continuous functions.) We shall also assume some functional calculus for $M_n(A)$, the $n$ by $n$ matrices over $A$. 
We write $a\ge 0$ to indicate that $a\in A$ is positive. 
We use \cite{Lance} as a reference for Hilbert $C^*$-modules.

If $E$ is an $A$-bimodule then $\overline{E}$ is identified with $E$ as a set but has the conjugate actions $a.\bar e=\overline{e.a^*}$ and $\bar e.a=\overline{a^*.e}$. Here $\bar e$ denotes $e\in E$ viewed in $\overline{E}$.  

\begin{defin}
An inner product on a bimodule $E$ is a bimodule map $\<,\>:E\tens_A\overline{E}\to A$
which satisfies the symmetry condition for all $x,y\in E$:
\begin{eqnarray*}
\<x,\overline{y}\>\ =\ \<y,\overline{x}\>^*\ .
\end{eqnarray*}
The inner product is positive if $\<x,\overline{x}\>\ge 0$ for all $x\in E$. 
\end{defin}

We now have the following well known result:

\begin{propos}\label{vcdsuiacvi}
Suppose we have inner products $\<,\>_E$, $\<,\>_F$ on $A$-bimodules $E,F$ respectively.
Then there is an inner product $\<,\>_{E\tens F}$ on $E\tens_A F$, given by the formula
\begin{eqnarray*}
\<x\tens y,\overline{x'\tens y'}\>_{E\tens F}\ =\ \<x.\<y,\overline{y'}\>_F,\overline{x'}\>_E\ ,
\end{eqnarray*}
for all $x,x'\in E$ and $y,y'\in F$. If $\<,\>_E$ and $\<,\>_F$ are both positive, and we can take the square root of positive matrices with enteries in $A$, then $\<,\>_{E\tens F}$
is also positive. 
\end{propos}

If we have positive Hermitian inner products on $E$ and on $\Omega^1 A$, 
using Proposition \ref{vcdsuiacvi} we can form the tensor product inner product
\begin{eqnarray*}
\<,\>_n:\big(\Omega^{\tens n}A\tens_A E\big) \tens_A \big(\overline{\Omega^{\tens n}A\tens_A E}\big)\to A\ .
\end{eqnarray*}
From this we can form the $A$-valued inner products of the $n$th derivatives;
\begin{eqnarray}
\<\!\< e,\overline{f}\>\!\>_n &=& \<\nabla^{(n)}e,\overline{\nabla^{(n)}f}\>_n\ .
\end{eqnarray}
Note that we have set $\nabla^{(0)}$ to be the identity.
Now the only missing part of following the classical definition is integration. Take a state $\phi:A\to\mathbb{C}$, i.e.\ a linear map preserving positivity and having
$\phi(1_A)=1$. We can 
 define a number valued semi-inner product on $E$ by
 \begin{eqnarray}\label{vcdisuyc}
\phi(\<\!\< e,\overline{f}\>\!\>_0)+\phi(\<\!\< e,\overline{f}\>\!\>_1)+\dots +\phi(\<\!\< e,\overline{f}\>\!\>_n)\ .
\end{eqnarray}
If the original inner product $\<,\>_E$ is strictly positive, and $\phi$ is a faithful state
(i.e.\ $\phi(x)=0$ for $x\ge 0$ implies $x=0$), then this is an inner product on $E$. 
The Sobolev space $W^{n,2}(E,\phi,\nabla)$ is the completion of $E$ under the
inner product in (\ref{vcdisuyc}).

\section{Another crossing map}
For $E\in {}_A\mathcal{E}_A$ and $F\in {}_A\mathcal{E}$, we shall define and study the map $\vartheta_E$ given by
Fig.\ \ref{newsig56svvuyd}. The significance of this will be explained in Section \ref{kdhjsvjkhcxr}. 

\unitlength 0.5 mm
\begin{picture}(135,90)(-20,15)
\linethickness{0.3mm}
\put(10,80){\line(0,1){10}}
\linethickness{0.3mm}
\put(30,80){\line(0,1){10}}
\linethickness{0.3mm}
%\multiput(10,60)(0.12,0.12){167}{\line(1,0){0.12}}
\multiput(10,60)(0.12,0.12){42}{\line(1,0){0.12}}
\multiput(30,80)(-0.12,-0.12){42}{\line(1,0){0.12}}
\linethickness{0.3mm}
%\multiput(10,80)(0.12,-0.12){42}{\line(1,0){0.12}}
\multiput(10,80)(0.12,-0.12){167}{\line(1,0){0.12}}
\linethickness{0.3mm}
%\multiput(22.5,65)(0.18,-0.12){42}{\line(1,0){0.18}}
\linethickness{0.3mm}
\put(50,60){\line(0,1){30}}
\linethickness{0.3mm}
\multiput(30,60)(0.01,-0.5){1}{\line(0,-1){0.5}}
\multiput(30.01,59.5)(0.04,-0.5){1}{\line(0,-1){0.5}}
\multiput(30.05,59)(0.06,-0.49){1}{\line(0,-1){0.49}}
\multiput(30.11,58.51)(0.09,-0.49){1}{\line(0,-1){0.49}}
\multiput(30.2,58.02)(0.11,-0.49){1}{\line(0,-1){0.49}}
\multiput(30.31,57.53)(0.14,-0.48){1}{\line(0,-1){0.48}}
\multiput(30.44,57.05)(0.16,-0.47){1}{\line(0,-1){0.47}}
\multiput(30.6,56.58)(0.09,-0.23){2}{\line(0,-1){0.23}}
\multiput(30.79,56.12)(0.1,-0.23){2}{\line(0,-1){0.23}}
\multiput(30.99,55.66)(0.11,-0.22){2}{\line(0,-1){0.22}}
\multiput(31.22,55.22)(0.12,-0.22){2}{\line(0,-1){0.22}}
\multiput(31.47,54.79)(0.14,-0.21){2}{\line(0,-1){0.21}}
\multiput(31.74,54.37)(0.15,-0.2){2}{\line(0,-1){0.2}}
\multiput(32.03,53.96)(0.1,-0.13){3}{\line(0,-1){0.13}}
\multiput(32.34,53.57)(0.11,-0.12){3}{\line(0,-1){0.12}}
\multiput(32.67,53.2)(0.12,-0.12){3}{\line(0,-1){0.12}}
\multiput(33.02,52.84)(0.12,-0.11){3}{\line(1,0){0.12}}
\multiput(33.38,52.5)(0.13,-0.11){3}{\line(1,0){0.13}}
\multiput(33.77,52.18)(0.13,-0.1){3}{\line(1,0){0.13}}
\multiput(34.16,51.88)(0.21,-0.14){2}{\line(1,0){0.21}}
\multiput(34.57,51.6)(0.21,-0.13){2}{\line(1,0){0.21}}
\multiput(35,51.34)(0.22,-0.12){2}{\line(1,0){0.22}}
\multiput(35.44,51.1)(0.22,-0.11){2}{\line(1,0){0.22}}
\multiput(35.89,50.88)(0.23,-0.1){2}{\line(1,0){0.23}}
\multiput(36.35,50.69)(0.47,-0.17){1}{\line(1,0){0.47}}
\multiput(36.82,50.52)(0.48,-0.15){1}{\line(1,0){0.48}}
\multiput(37.29,50.37)(0.48,-0.12){1}{\line(1,0){0.48}}
\multiput(37.77,50.25)(0.49,-0.1){1}{\line(1,0){0.49}}
\multiput(38.26,50.15)(0.49,-0.07){1}{\line(1,0){0.49}}
\multiput(38.76,50.08)(0.5,-0.05){1}{\line(1,0){0.5}}
\multiput(39.25,50.03)(0.5,-0.02){1}{\line(1,0){0.5}}
\put(39.75,50){\line(1,0){0.5}}
\multiput(40.25,50)(0.5,0.02){1}{\line(1,0){0.5}}
\multiput(40.75,50.03)(0.5,0.05){1}{\line(1,0){0.5}}
\multiput(41.24,50.08)(0.49,0.07){1}{\line(1,0){0.49}}
\multiput(41.74,50.15)(0.49,0.1){1}{\line(1,0){0.49}}
\multiput(42.23,50.25)(0.48,0.12){1}{\line(1,0){0.48}}
\multiput(42.71,50.37)(0.48,0.15){1}{\line(1,0){0.48}}
\multiput(43.18,50.52)(0.47,0.17){1}{\line(1,0){0.47}}
\multiput(43.65,50.69)(0.23,0.1){2}{\line(1,0){0.23}}
\multiput(44.11,50.88)(0.22,0.11){2}{\line(1,0){0.22}}
\multiput(44.56,51.1)(0.22,0.12){2}{\line(1,0){0.22}}
\multiput(45,51.34)(0.21,0.13){2}{\line(1,0){0.21}}
\multiput(45.43,51.6)(0.21,0.14){2}{\line(1,0){0.21}}
\multiput(45.84,51.88)(0.13,0.1){3}{\line(1,0){0.13}}
\multiput(46.23,52.18)(0.13,0.11){3}{\line(1,0){0.13}}
\multiput(46.62,52.5)(0.12,0.11){3}{\line(1,0){0.12}}
\multiput(46.98,52.84)(0.12,0.12){3}{\line(0,1){0.12}}
\multiput(47.33,53.2)(0.11,0.12){3}{\line(0,1){0.12}}
\multiput(47.66,53.57)(0.1,0.13){3}{\line(0,1){0.13}}
\multiput(47.97,53.96)(0.15,0.2){2}{\line(0,1){0.2}}
\multiput(48.26,54.37)(0.14,0.21){2}{\line(0,1){0.21}}
\multiput(48.53,54.79)(0.12,0.22){2}{\line(0,1){0.22}}
\multiput(48.78,55.22)(0.11,0.22){2}{\line(0,1){0.22}}
\multiput(49.01,55.66)(0.1,0.23){2}{\line(0,1){0.23}}
\multiput(49.21,56.12)(0.09,0.23){2}{\line(0,1){0.23}}
\multiput(49.4,56.58)(0.16,0.47){1}{\line(0,1){0.47}}
\multiput(49.56,57.05)(0.14,0.48){1}{\line(0,1){0.48}}
\multiput(49.69,57.53)(0.11,0.49){1}{\line(0,1){0.49}}
\multiput(49.8,58.02)(0.09,0.49){1}{\line(0,1){0.49}}
\multiput(49.89,58.51)(0.06,0.49){1}{\line(0,1){0.49}}
\multiput(49.95,59)(0.04,0.5){1}{\line(0,1){0.5}}
\multiput(49.99,59.5)(0.01,0.5){1}{\line(0,1){0.5}}

\linethickness{0.3mm}
\put(40,30){\line(0,1){20}}
\linethickness{0.3mm}
\put(10,30){\line(0,1){30}}
\linethickness{0.3mm}
\put(90,70){\line(0,1){20}}
\linethickness{0.3mm}
\put(110,70){\line(0,1){20}}
\linethickness{0.3mm}
\put(130,70){\line(0,1){20}}
\linethickness{0.3mm}
\put(85,70){\line(1,0){50}}
\put(85,60){\line(0,1){10}}
\put(135,60){\line(0,1){10}}
\put(85,60){\line(1,0){50}}
\linethickness{0.3mm}
\put(100,30){\line(0,1){30}}
\linethickness{0.3mm}
\put(120,30){\line(0,1){30}}
\put(99,25){\makebox(0,0)[cc]{$E$}}

\put(119,25){\makebox(0,0)[cc]{$F$}}

\put(39,25){\makebox(0,0)[cc]{$F$}}

\put(9,25){\makebox(0,0)[cc]{$E$}}

\put(29,95){\makebox(0,0)[cc]{$E$}}

\put(49,95){\makebox(0,0)[cc]{$F$}}

\put(109,95){\makebox(0,0)[cc]{$E$}}

\put(129,95){\makebox(0,0)[cc]{$F$}}

\put(5,95){\makebox(0,0)[cc]{$\mathcal{T}\Vec A$}}

\put(85,95){\makebox(0,0)[cc]{$\mathcal{T}\Vec A$}}

\put(108.5,65){\makebox(0,0)[cc]{$\la$}}

\put(45,47.5){\makebox(0,0)[cc]{$\la$}}

\put(28,70){\makebox(0,0)[cc]{$\vartheta_E$}}

\put(70,67){\makebox(0,0)[cc]{$=$}}

%\put(136,50){\makebox(0,0)[cc]{\textbf{Fig.\ 99}}}
%%%%%%%%%%%%%%%
\refstepcounter{piccie} \label{newsig56svvuyd}
\put(160,45){\makebox(0,0)[cc]{Fig.\ \arabic{piccie}}}
%%%%%%%%%%%%%%%

\end{picture}

\noindent As the reader should be used to by now, we shall construct the map
$\vartheta_E$ by recursion on $n$ where $\vartheta_E:\Vec^{\tens n} A\tens E\to
E\tens_A \mathcal{T}\Vec A$, where we will worry about just what sort of tensor product on the domain or whether it is a module map later. 
To do this we start with $n=0$ and $\vartheta_E:A\tens_A E\to E\tens_A A$ being the identity.
For $n=1$ using the formula for the action of $\Vec A$ on a tensor product
(given by $\nabla_{E\tens F}$) in Fig.\ \ref{newsig56svvuyd} gives
\begin{eqnarray}  \label{cnjdibvb}
\vartheta_E\ =\ \la+ \sigma_E^{-1}:\Vec A\tens E\to
E\tens_A \mathcal{T}\Vec A\ .
\end{eqnarray}
The first term on the right hand side in (\ref{cnjdibvb}) is in $E\tens_A A$ and the second in $E\tens_A \Vec A$. 
Now we calculate, for $v\in \Vec A$, $a\in A$ and $e\in E$
\begin{eqnarray}
\vartheta_E(v\tens a.e) &=& v(\extd a).e+\vartheta_E(v.a\tens e)\ ,
\end{eqnarray}
and conclude that we do not get a map from $\Vec A\tens_A E$. However the reader should recall that this is not the right $A$-module structure
for $\mathcal{T}\Vec A$ used in Theorem \ref{kuykjvcvu}. The structure used there was $\mathcal{T}\Vec A_\bullet$, where we have $v\bullet a=v.a+v(\extd a)$. This means that $\vartheta_E$ does give a well defined map from $(A\oplus \Vec A)_\bullet \tens_A E$ to $E\tens_A \mathcal{T}\Vec A$.  Next note that 
(\ref{cnjdibvb}) gives a left $A$-module map. Finally, for the right module structure:
\begin{eqnarray}
\vartheta_E(v\tens e.a) &=& (\ev\tens\id_E)(v\tens\nabla_E(e.a))+\sigma_E^{-1}(v\tens e.a) \cr
&=& \vartheta_E(v\tens e).a + (\ev\tens\id_E)(\id\tens\sigma_E)(v\tens e\tens\extd a)\cr
&=&  \vartheta_E(v\tens e).a + (\id_E\tens\ev)
(\sigma_E^{-1}(v\tens e)\tens\extd a) \cr
&=& \vartheta_E(v\tens e)\bullet a\ .
\end{eqnarray}
We conclude that $\vartheta_E: (A\oplus \Vec A)_\bullet \tens_A E \to E\tens_A \mathcal{T}\Vec A_\bullet$ is an $A$-bimodule map, and this shall be the basis for our recursive construction. We have used $(A\oplus \Vec A)_\bullet$ as $\Vec A$ itself is not a right $A$ module under the $\bullet$ product. For the same reason, we set $\Vec A^{\tens\le n}=A\oplus\Vec A\oplus\Vec A^{\tens 2}\oplus\dots \oplus\Vec A^{\tens n}$ in the following Proposition.

\begin{propos}  \label{jkytfryjhzzdt}
Suppose that $\vartheta_E$ is defined recursively by, for $w\in \Vec A$ and $\underline{v}\in \Vec^{\tens n} A$,
\begin{eqnarray}  \label{fcauyxaert}
\vartheta_E((w\tens \underline{v})\tens e) &=& (\la\tens\id)(\id\tens\vartheta_E)(w\tens \underline{v}\tens e) \cr
&&+\ (\sigma_E^{-1}\tens\id)(\id\tens\vartheta_E)(w\tens \underline{v}\tens e) \cr
&&+\ (\id_E\tens\hat\bullet)(\sigma_E^{-1}\tens\id)(\id\tens\vartheta_E)(w\tens \underline{v}\tens e) \cr
&& -\   \vartheta_E((w\bullet_n \underline{v})\tens e)\ .
\end{eqnarray}
Here $\hat\bullet=\bullet_m:\Vec A\tens \Vec^{\tens m}A\to \Vec^{\tens m}A$ -- we do not give it a specific index in (\ref{fcauyxaert}) as $m$ may vary. Then, for all $n\ge 1$,

$1_n$)\quad $\vartheta_E:\Vec^{\tens\le n}A\tens E\to E\tens_A \mathcal{T}\Vec A$ is well defined;

$2_n$)\quad $\vartheta_E:\Vec^{\tens\le n}A_\bullet\tens_A E\to E\tens_A \mathcal{T}\Vec A$ is well defined;

$3_n$)\quad $\vartheta_E:\Vec^{\tens\le n}A_\bullet\tens_A E\to E\tens_A \mathcal{T}\Vec A$ is a left $A$ module map;

$4_n$)\quad $\vartheta_E:\Vec^{\tens\le n}A_\bullet\tens_A E\to E\tens_A \mathcal{T}\Vec A_\bullet$ is a right $A$ module map;

$5_n$)\quad $
\underline{v}\,\la\,(e\tens f) \,=\, (\id_E\tens\la)(\vartheta_E\tens\id_F)(\underline{v}\tens e\tens f)$, for all $v\in \Vec^{\tens n}A$.

\end{propos}
\noindent {\bf Proof:}\quad Proof by induction. Assume that $\vartheta_E:\Vec^{\tens\le n}A_\bullet\tens_A E\to E\tens_A \mathcal{T}\Vec A_\bullet$ is defined and satisfies
($1_n$,\dots,$5_n$) -- noting that the $n=1$ case is done already -- see (\ref{cnjdibvb}) and the discussion following it. Now we use (\ref{fcauyxaert}) to give the $n+1$ case, and verify the corresponding statements. 

The first thing is to check that the right hand side of (\ref{fcauyxaert}) is actually well defined, given that $\vartheta_E$ maps into $E\tens_A \mathcal{T}\Vec A_\bullet$ (emphasising the $\tens_A$). For this, we need to check that, for all $a\in A$, $\underline{u}\in \Vec^{\tens m} A$,
$w\in \Vec A$ and $e\in E$,
\begin{eqnarray}  \label{vcasjychgx}
\big(\la\tens\id+(\id_E\tens\bullet_m)(\sigma_E^{-1}\tens\id)\big)(w\tens e.a\tens \underline{u}-w\tens e\tens a.\underline{u})\ =\ 0\ .
\end{eqnarray}
Look at
\begin{eqnarray}  \label{vcasjybvsdchgx}
\big(\la\tens\id\big)(w\tens e.a\tens \underline{u}-w\tens e\tens a.\underline{u}) &=& (\ev\tens\id_E\tens\id)(w\tens \sigma_E(e\tens\extd a)\tens \underline{u})\cr
&=&(\id_E\tens\ev\tens\id)(\sigma_E^{-1}(w\tens e)\tens \extd a\tens \underline{u})\ .
\end{eqnarray}
Verifying (\ref{vcasjychgx}) reduces to showing, for $w'\in \Vec A$,
\begin{eqnarray}
w'(\extd a).\underline{u} +
\bullet_m(w'.a\tens \underline{u}-w\tens a.\underline{u})\ =\ 0\ ,
\end{eqnarray}
which comes from the definition of $\bullet_m$ in Lemma \ref{hjsvhjkvjhk} and the left Leibniz rule for $\square$. 

To prove ($1_{n+1}$), we use ($3_{n}$) to show that the middle two
terms of (\ref{fcauyxaert}) evaluated on $w\tens a.\underline{v}\tens e$ are the same as on $w.a\tens \underline{v}\tens e$. 
For the first and fourth terms, we have
\begin{eqnarray*}
(\la\tens\id)(\id\tens\vartheta_E)(w\tens a.\underline{v}\tens e-w.a\tens \underline{v}\tens e) &=&
(\la\tens\id)(w\tens a.\vartheta_E(\underline{v}\tens e)-w.a\tens \vartheta_E(\underline{v}\tens e))\cr
&=& w(\extd a).\vartheta_E(\underline{v}\tens e)\ ,\cr
\vartheta_E((w\bullet_n (a.\underline{v}))\tens e)-\vartheta_E(((w.a)\bullet_n \underline{v})\tens e) &=& 
\vartheta_E(w(\extd a).\underline{v}\tens e)\ =\ w(\extd a).\vartheta_E(\underline{v}\tens e)\ .
\end{eqnarray*}
These cancel in (\ref{fcauyxaert}), verifying ($1_{n+1}$). 

To prove ($2_{n+1}$), it is convenient to rewrite (\ref{fcauyxaert}) as
\begin{eqnarray}  \label{fcauyxcdsaert}
\vartheta_E((w\bullet \underline{v})\tens e) &=& (\la\tens\id)(\id\tens\vartheta_E)(w\tens \underline{v}\tens e) \cr
&&+\ (\sigma_E^{-1}\tens\id)(\id\tens\vartheta_E)(w\tens \underline{v}\tens e) \cr
&&+\ (\id_E\tens\hat\bullet)(\sigma_E^{-1}\tens\id)(\id\tens\vartheta_E)(w\tens \underline{v}\tens e) \ .
\end{eqnarray}
We need to show from (\ref{fcauyxcdsaert}) that $\vartheta_E((w\bullet \underline{v})\bullet a\tens e)=\vartheta_E((w\bullet \underline{v})\tens a.e)$. By associativity of $\bullet$, it is sufficient to verify that 
$\vartheta_E(w\bullet (\underline{v}\bullet a)\tens e)=\vartheta_E((w\bullet \underline{v})\tens a.e)$. To do this, it is enough to show that the right hand side of (\ref{fcauyxcdsaert}) is the same when applied to $w\tens \underline{v}\bullet a\tens e$ and to $w\tens \underline{v}\tens a.e$, but this is true from ($2_{n}$). 

Proving ($3_{n+1}$) is quite simple, as every term on the right hand side of (\ref{fcauyxaert}) is left $A$-linear in $w$.

To prove ($4_{n+1}$) we use (\ref{fcauyxcdsaert}) again. Set $\vartheta_E(\underline{v}\tens e)=f\tens u$. Then
 ($4_{n}$) implies that $\vartheta_E(\underline{v}\tens e.a)=f\tens u\bullet a$. The right hand side of (\ref{fcauyxcdsaert}) applied to $w\tens \underline{v}\tens e.a$ instead of $w\tens \underline{v}\tens e$ is
\begin{eqnarray} \label{xrtyhgfxzsx}
w\,\la\,f\tens u\bullet a +(\id_E\tens\bullet)(\sigma_E^{-1}(w\tens f)\tens u\bullet a)\ ,
\end{eqnarray}
where we have combined the last two terms of (\ref{fcauyxcdsaert}) to give the last term
of (\ref{xrtyhgfxzsx}). Now associativity of $\bullet$ gives the answer.

Finally we consider ($5_{n+1}$). From Lemma \ref{vcjhckvjvhjvj} we can write
\begin{eqnarray}
(w\tens \underline{v})\,\la\,(e\tens f) &=& w\,\la\,(\underline{v}\,\la\,(e\tens f))-(w\bullet_n \underline{v})\,\la\,(e\tens f)\ ,
\end{eqnarray}
and using the inductive hypothesis ($5_{n}$) gives
\begin{eqnarray}  \label{vcgdjsyvxd}
(w\tens \underline{v})\,\la\,(e\tens f) &=& w\,\la\,((\id_E\tens\la)(\vartheta_E\tens\id_F)(\underline{v}\tens e\tens f)) \cr
&& -\   (\id_E\tens\la)(\vartheta_E\tens\id_F)((w\bullet_n \underline{v})\tens e\tens f)\cr
&=& (\la\tens\la)(\id\tens\vartheta_E\tens\id_F)(w\tens \underline{v}\tens e\tens f) \cr
&&+\ (\id_E\tens\la)(\sigma_E^{-1}\tens\id_F)(\id\tens\id_E\tens\la)(\id\tens\vartheta_E\tens\id_F)(w\tens \underline{v}\tens e\tens f) \cr
&& -\   (\id_E\tens\la)(\vartheta_E\tens\id_F)((w\bullet_n \underline{v})\tens e\tens f)\cr
&=& (\la\tens\la)(\id\tens\vartheta_E\tens\id_F)(w\tens \underline{v}\tens e\tens f) \cr
&&+\ (\id_E\tens\la)(\id_E\tens\id\tens\la)(\sigma_E^{-1}\tens\id\tens\id_F)(\id\tens\vartheta_E\tens\id_F)(w\tens \underline{v}\tens e\tens f) \cr
&& -\   (\id_E\tens\la)(\vartheta_E\tens\id_F)((w\bullet_n \underline{v})\tens e\tens f)\ .
\end{eqnarray}
By using Lemma \ref{vcjhckvjvhjvj} again, we can rewrite the middle term of the result of 
(\ref{vcgdjsyvxd}) as follows, 
\begin{eqnarray}  \label{vcgdjsyvxdff}
&& (w\tens \underline{v})\,\la\,(e\tens f) \cr
&=& (\la\tens\la)(\id\tens\vartheta_E\tens\id_F)(w\tens \underline{v}\tens e\tens f) \cr
&&+\ (\id_E\tens\la)(\sigma_E^{-1}\tens\id\tens\id_F)(\id\tens\vartheta_E\tens\id_F)(w\tens \underline{v}\tens e\tens f) \cr
&&+\ (\id_E\tens\la)(\id_E\tens\hat\bullet\tens\id_F)(\sigma_E^{-1}\tens\id\tens\id_F)(\id\tens\vartheta_E\tens\id_F)(w\tens \underline{v}\tens e\tens f) \cr
&& -\   (\id_E\tens\la)(\vartheta_E\tens\id_F)((w\bullet_n \underline{v})\tens e\tens f)\ .
\end{eqnarray}
Now (\ref{vcgdjsyvxdff})  implies ($5_{n+1}$) by the recursive definition (\ref{fcauyxaert}). \quad
$\square$

\medskip Note that we can rewrite the recursive definition (\ref{fcauyxaert}) in a shorter form as
\begin{eqnarray}  \label{fcauyxaeffrt}
\vartheta_E(w\bullet \underline{v}\tens e) &=& (\la\tens\id)(\id\tens\vartheta_E)(w\tens \underline{v}\tens e) \cr
&&+\ (\sigma_E^{-1}\bullet\id)(\id\tens\vartheta_E)(w\tens \underline{v}\tens e) \ .
\end{eqnarray}

\begin{propos}  \label{cvadhjskxsz1}
\begin{eqnarray*}
\vartheta_{E\tens F}\ =\ (\id_E\tens \vartheta_{F})(\vartheta_{E}\tens\id_F):
\mathcal{T}\Vec A_\bullet\tens_A E\tens_A F\to E\tens_A F\tens_A \mathcal{T}\Vec A_\bullet\ .
\end{eqnarray*}
\end{propos}
\noindent {\bf Proof:}\quad This is proved by induction on $n$, where 
$\vartheta_{E\tens F}:\Vec^{\tens\le n}A_\bullet\tens_A E\tens_A F\to E\tens_A F\tens_A \mathcal{T}\Vec A_\bullet$. 

The $n=1$ case is given by (\ref{cnjdibvb}) as
\begin{eqnarray}  \label{ccacnjdibvb}
\vartheta_{E\tens F}\ =\ \la_{E\tens F}+ \sigma_{E\tens F}^{-1}:\Vec A_\bullet\tens_A E\tens_A F\to
E\tens_A F\tens_A \mathcal{T}\Vec A_\bullet\ .
\end{eqnarray}
Next write
\begin{eqnarray}  \label{ccacnjvddibvb}
(\id_E\tens\vartheta_F)(\vartheta_{E}\tens\id_F)\ =\ (\id_E\tens\vartheta_F)(\la_E\tens\id_F+
\sigma_{E}^{-1}\tens\id_F)\ .
\end{eqnarray}
Remember that $\vartheta_F$ is essentially the identity on $\Vec^{\tens 0}A_\bullet\tens_A F$, so we obtain
\begin{eqnarray}  \label{ccacnvsjdibvb}
(\id_E\tens\vartheta_F)(\vartheta_{E}\tens\id_F)\ =\ \la_E\tens\id_F + (\id_E\tens\vartheta_F)(
\sigma_{E}^{-1}\tens\id_F)\ ,
\end{eqnarray}
and, using  the formula for the action of vector fields on a tensor product,
this is the same as (\ref{ccacnjdibvb}).

Now suppose that the hypothesis works for $n$.
For all $\underline{v}\in \Vec^{\tens n} A$,
$w\in \Vec A$, $e\in E$ and $f\in F$, the recursive definition (\ref{fcauyxaeffrt}) gives
\begin{eqnarray}  \label{fcadsvfrt}
\vartheta_{E\tens F}(w\bullet \underline{v}\tens e\tens f) &=& (\la_{E\tens F}\tens\id)(\id\tens\vartheta_{E\tens F})(w\tens \underline{v}\tens e\tens f) \cr
&&+\ (\sigma_{E\tens F}^{-1}\bullet\id)(\id\tens\vartheta_{E\tens F})(w\tens \underline{v}\tens e\tens f) \ .
\end{eqnarray}
By using the inductive hypothesis, we write
\begin{eqnarray}
&& (\sigma_{E\tens F}^{-1}\tens\id)(\id\tens\vartheta_{E\tens F}) \cr
&=& (\id_E\tens \sigma_{F}^{-1}\tens\id)(\sigma_{E}^{-1}\tens\id_F\tens\id)
(\id\tens\id_E\tens \vartheta_{F})(\id\tens \vartheta_{E}\tens\id_F)  \cr
&=& (\id_E\tens \sigma_{F}^{-1}\tens\id)(\id_E\tens\id\tens \vartheta_{F})
(\sigma_{E}^{-1}\tens\id\tens\id_F)
(\id\tens \vartheta_{E}\tens\id_F)\ ,
\end{eqnarray}
and from this we get the second term of (\ref{fcadsvfrt}),
\begin{eqnarray}
&& (\sigma_{E\tens F}^{-1}\bullet\id)(\id\tens\vartheta_{E\tens F}) \cr
&=& (\id_E\tens \sigma_{F}^{-1}\tens\id)(\sigma_{E}^{-1}\tens\id_F\tens\id)
(\id\tens\id_E\tens \vartheta_{F})(\id\tens \vartheta_{E}\tens\id_F)  \cr
&=& \big(\id_E\tens (\sigma_{F}^{-1}\bullet\id)(\id\tens \vartheta_{F})\big)\ \big(
(\sigma_{E}^{-1}\tens\id)
(\id\tens \vartheta_{E})\tens\id_F\big)\ .
\end{eqnarray}
Using (\ref{fcauyxaeffrt}) twice gives
\begin{eqnarray}  \label{vagyuiuryxz}
&& (\sigma_{E\tens F}^{-1}\bullet\id)(\id\tens\vartheta_{E\tens F}) \cr
&=& \big(\id_E\tens \vartheta_{F}(\bullet\tens\id_F)\big)\ \big(
(\sigma_{E}^{-1}\tens\id)
(\id\tens \vartheta_{E})\tens\id_F\big)\cr
&& -\  \big(\id_E\tens  (\la_F\tens\id)(\id\tens\vartheta_F)  \big)\ \big(
(\sigma_{E}^{-1}\tens\id)
(\id\tens \vartheta_{E})\tens\id_F\big)\cr
&=& \big(\id_E\tens \vartheta_{F}\big)\ \big(
(\sigma_{E}^{-1}\bullet\id)
(\id\tens \vartheta_{E})\tens\id_F\big)\cr
&& -\  \big(\id_E\tens  (\la_F\tens\id)(\id\tens\vartheta_F)  \big)\ \big(
(\sigma_{E}^{-1}\tens\id)
(\id\tens \vartheta_{E})\tens\id_F\big)\cr
&=& \big(\id_E\tens \vartheta_{F}\big)\ \big(
 \vartheta_{E}(\bullet\tens\id_E)\tens\id_F\big)\cr
&&  -\  \big(\id_E\tens \vartheta_{F}\big)\ \big(
(\la_E\tens\id)
(\id\tens \vartheta_{E})\tens\id_F\big)\cr
&& -\  \big(\id_E\tens  (\la_F\tens\id)(\id\tens\vartheta_F)  \big)\ \big(
(\sigma_{E}^{-1}\tens\id)
(\id\tens \vartheta_{E})\tens\id_F\big)\ .
\end{eqnarray}
At this point we consider the last two terms of (\ref{vagyuiuryxz}) separately: The second is
\begin{eqnarray}
&&  (\id_E\tens\vartheta_F)(\la_E\tens\id\tens\id_F)(\id\tens\vartheta_E\tens\id_F) \cr
&=& (\la_E\tens\id_F\tens\id)(\id\tens\id_E\tens\vartheta_F)(\id\tens\vartheta_E\tens\id_F)\cr
&=& (\la_E\tens\id_F\tens\id)(\id\tens\vartheta_{E\tens F})\ ,
\end{eqnarray}
and the third is
\begin{eqnarray}
&& (\id_E\tens  \la_F\tens\id)   (\id_E\tens \id\tens\vartheta_F) 
(\sigma_{E}^{-1}\tens\id\tens\id_F)     (\id\tens \vartheta_{E}\tens\id_F) \cr
&=& (\id_E\tens  \la_F\tens\id)   (\sigma_{E}^{-1}\tens\id_F\tens\id)  
  (\id\tens \id_E\tens\vartheta_F)   (\id\tens \vartheta_{E}\tens\id_F) \cr
  &=& (\id_E\tens  \la_F\tens\id)   (\sigma_{E}^{-1}\tens\id_F\tens\id)  
  (\id\tens \vartheta_{E\tens F})  \ .
\end{eqnarray}
Now we put these results back into (\ref{fcadsvfrt}) to get
\begin{eqnarray}  \label{fcadsvfrsdbvsvt}
&& \vartheta_{E\tens F}(w\bullet \underline{v}\tens e\tens f) \cr
 &=& \big(\id_E\tens \vartheta_{F}\big)\ \big(
 \vartheta_{E}\tens\id_F\big)(w\bullet \underline{v}\tens e\tens f) \cr
  && -\ (\la_E\tens\id_F\tens\id)(\id\tens\vartheta_{E\tens F})(w\tens \underline{v}\tens e\tens f) \cr
 && -\ (\id_E\tens  \la_F\tens\id)   (\sigma_{E}^{-1}\tens\id_F\tens\id)  
  (\id\tens \vartheta_{E\tens F}) (w\tens \underline{v}\tens e\tens f) \cr
&&+\ (\la_{E\tens F}\tens\id)(\id\tens\vartheta_{E\tens F})(w\tens \underline{v}\tens e\tens f) \ .
\end{eqnarray}
The last three terms of (\ref{fcadsvfrsdbvsvt}) cancel (remember that $w\in \Vec A$), proving the hypothesis for the case $n+1$.\quad$\square$

\begin{propos} \label{cvadhjskxsz2} For all $\underline{u},\underline{v}\in \mathcal{T}\Vec A$,
\begin{eqnarray*}
\vartheta_E(\underline{u}\bullet \underline{v}\tens e) &=& (\id_E\tens\bullet)(\vartheta_E\tens\id)(\id\tens \vartheta_E)(\underline{u}\tens \underline{v}\tens e)\ .
\end{eqnarray*}
\end{propos}
\noindent {\bf Proof:}\quad This is proved by induction on $m$, where $\underline{u}\in \Vec^{\tens m}A$. First the $m=1$ case is given by combining  (\ref{cnjdibvb})  and (\ref{fcauyxaeffrt}). Now suppose that it is true for some $m\ge 1$, and consider, for all $w\in \Vec A$,
\begin{eqnarray*}
\vartheta_E((w\bullet \underline{u})\bullet \underline{v}\tens e) &=& \vartheta_E(w\bullet (\underline{u}\bullet \underline{v})\tens e)\cr
&=& (\id_E\tens\bullet)(\vartheta_E\tens\id)(\id\tens \vartheta_E)(w\tens \underline{u}\bullet \underline{v}\tens e)\cr
&=& (\id_E\tens\bullet)(\vartheta_E\tens\id)\big(w\tens (\id_E\tens\bullet)(\vartheta_E\tens\id)(\id\tens \vartheta_E)(\underline{u}\tens \underline{v}\tens e)\big)
\end{eqnarray*}
where we have used, in order, associativity of $\bullet$, the $m=1$ case, and the inductive hypothesis. This can be rearranged to give
\begin{eqnarray*}
&& \vartheta_E((w\bullet \underline{u})\bullet \underline{v}\tens e) \cr
&=& (\id_E\tens\bullet(\id\tens\bullet)(\vartheta_E\tens\id\tens\id)
(\id\tens\vartheta_E\tens\id)(\id\tens\id\tens\vartheta_E)(w\tens \underline{u}\tens \underline{v}\tens e)  \cr
&=& (\id_E\tens\bullet(\bullet\tens\id)(\vartheta_E\tens\id\tens\id)
(\id\tens\vartheta_E\tens\id)(\id\tens\id\tens\vartheta_E)(w\tens \underline{u}\tens \underline{v}\tens e)  \cr
&=& (\id_E\tens\bullet(\vartheta_E\tens\id)
(\id\tens\vartheta_E)(w\bullet \underline{u}\tens \underline{v}\tens e)  \ ,
\end{eqnarray*}
where we have used the associativity of $\bullet$, and the $m=1$ case. This concludes the inductive proof.\quad$\square$

\begin{propos}  \label{cvadhjskxsz3}
The following map is simply the $\bullet$ product:
\begin{eqnarray*}
\vartheta_A:\mathcal{T}\Vec A_\bullet\tens_A A\to A\tens_A \mathcal{T}\Vec A_\bullet\ \cong \mathcal{T}\Vec A_\bullet\ .
\end{eqnarray*}
\end{propos}
\noindent {\bf Proof:}\quad First note that $\vartheta_A(v\tens a)=\vartheta_A(v\bullet a\tens 1_A)$. 
Now check the formulae for $\vartheta_A(u\tens 1_A)$, given that any differential applied to $1_A$ gives zero as $\extd(1_A)=0$.\quad$\square$

\medskip The map $\vartheta_E$ maps $\Vec^{\tens n} A_\bullet\tens_A E$ to the sum of 
$E\tens_A \Vec A^{\tens m}_\bullet$ for $m\le n$. On the level of $E\tens_A \Vec A^{\tens n}_\bullet$
(i.e.\ ignoring the $m<n$ terms) it is simply $\sigma_E^{-1}$, which is assumed to be invertible. By standard graded arguments, it is likely that $\vartheta_E$ is invertible. Here we shall give a recursive definition of its inverse.

\begin{propos}  \label{cvadhjskxsz4}
The inverse $\vartheta_E^{-1}:E\tens_A \Vec^{\tens n} A_\bullet \to \mathcal{T}\Vec A_\bullet \tens_A E$ is given recursively by the identity for $n=0$, $\vartheta_E^{-1}(e\tens w)=(\id\tens\id_E-\la)\,\sigma_E$ for $n=1$, and then by, for  $f\in E$, $v\in \Vec^{\tens n} A$ and $u\in  \Vec A$,
\begin{eqnarray*}
\vartheta_E^{-1}(f\tens u\tens v) &=& (\id\bullet \vartheta_E^{-1})(\sigma_E\tens\id)(f\tens u\tens v)-
\vartheta_E^{-1}(f\tens u\,\la\, v)\cr && -\ \vartheta_E^{-1}((\la\,\sigma_E\tens\id)(f\tens u\tens v))\ .
\end{eqnarray*}
This can be rewritten as
\begin{eqnarray*}
\vartheta_E^{-1}(f\tens u\bullet v) &=& (  \id\bullet \vartheta_E^{-1}  - \vartheta_E^{-1}(\la\tens\id   ) )
(\sigma_E(f\tens u)\tens v)\ .
\end{eqnarray*}
\end{propos}
\noindent {\bf Proof:}\quad The $n=1$ case is by explicit calculation.

Now assume that the formula for $\vartheta_E^{-1}$ on $E\tens_A \Vec^{\tens n} A_\bullet$ works. Then we can use (\ref{fcauyxaeffrt}) to write, for $v\in \Vec^{\tens n} A$ and $w\in  \Vec A$,
\begin{eqnarray}
\vartheta_E(w\bullet \vartheta_E^{-1}(e\tens v)) &=& w\,\la\,e\tens v +\sigma_E^{-1}(w\tens e)\bullet v\ .
\end{eqnarray}
If we write $\sigma_E^{-1}(w\tens e)=f\tens u$ for $u\in  \Vec A$ (summation implicit), then
\begin{eqnarray}  \label{cvtauyiytrs}
\vartheta_E(\id\bullet \vartheta_E^{-1})(\sigma_E\tens\id)(f\tens u\tens v) &=&
(\la\,\sigma_E\tens\id)(f\tens u\tens v) + f\tens u\bullet v  \cr
&=& (\la\,\sigma_E\tens\id)(f\tens u\tens v) + f\tens u\tens v \cr && +\  f\tens u\,\la\, v    \ .
\end{eqnarray}
As two of the terms of the right side of (\ref{cvtauyiytrs}) are in the domain of previously defined
$\vartheta_E^{-1}$, we can rewrite  (\ref{cvtauyiytrs}) as
\begin{eqnarray}
f\tens u\tens v &=& \vartheta_E\Big((\id\bullet \vartheta_E^{-1})(\sigma_E\tens\id)(f\tens u\tens v)-
\vartheta_E^{-1}(f\tens u\,\la\, v)\cr && -\ \vartheta_E^{-1}((\la\,\sigma_E\tens\id)(f\tens u\tens v))
\Big)\ ,
\end{eqnarray}
and this gives the recursive formula in the statement. We leave checking other properties of $\vartheta_E^{-1}$ to the reader.\quad$\square$

\section{A rather unusual covariant derivative}
We give a left bimodule covariant derivative on $\mathcal{T}\Vec A_\bullet$. Given 
the dual basis $\coev(1)=\xi\tens u\in\Omega^1 A\tens_A\Vec A$ (summation implicit), define $\nabla:\mathcal{T}\Vec A\to \Omega^1 A\tens\mathcal{T}\Vec A$ by 
\begin{eqnarray}
\nabla(\underline{v}) &=& \xi\tens (u\bullet\underline{v})\ .
\end{eqnarray}
As $a.u=a\bullet u$ for $a\in A$, we see that $\nabla$ only depends on 
$\xi\tens u\in\Omega^1 A\tens_A\Vec A$ (with the emphasis on $\tens_A$). 
We need to look at the defining characteristic of a left covariant derivative,
\begin{eqnarray}
\nabla(a\bullet \underline{v}) &=& \xi\tens (u\bullet(a\bullet\underline{v}))  \cr
 &=& \xi\tens ((u\bullet a)\bullet\underline{v}) \cr
  &=&  \xi\tens (u(\extd a)\bullet\underline{v}) + \xi\tens ((u. a)\bullet\underline{v}) \ .
\end{eqnarray}
Next we use $\coev(a)=a.\xi\tens u=\xi\tens u.a\in\Omega^1 A\tens_A\Vec A$ to compute
\begin{eqnarray}
\nabla(a\bullet \underline{v}) 
  &=&  \xi\tens u(\extd a).\underline{v} + u.\xi\tens (u\bullet\underline{v})\cr
&=&  \xi.u(\extd a)\tens \underline{v} + u.\xi\tens (u\bullet\underline{v}) \cr
&=& \extd a \tens \underline{v} + u.\nabla(\underline{v})  \ .
\end{eqnarray}
Here we have also used the defining property of a dual basis, $\xi.u(\extd a)=\extd a$. 
Now we look at $\nabla$ as a bimodule covariant derivative on $\mathcal{T}\Vec A_\bullet$. Calculate
\begin{eqnarray} \label{bchdjkshjcgvgxfh}
\nabla(\underline{v}\bullet a)-\nabla(\underline{v})\bullet a &=& 
\xi\tens (u\bullet(\underline{v}\bullet a)) - \xi\tens ((u\bullet\underline{v})\bullet a)\ =\ 0\ .
\end{eqnarray}
We come to the rather surprising conclusion that $(\mathcal{T}\Vec A_\bullet,\nabla,0)$ is a bimodule covariant derivative. 

Suppose that $(E,\nabla_E,\sigma_E)$ is a bimodule covariant derivative. Then consider the tensor product covariant derivatives,
\begin{eqnarray}  \label{cvgdjsahgfx}
\nabla_{E\tens \mathcal{T}\Vec A_\bullet}(e\tens \underline{v}) &=& 
\nabla_E(e)\tens  \underline{v} + \sigma_E(e\tens\xi)\tens (u\bullet  \underline{v})\ ,\cr
\nabla_{ \mathcal{T}\Vec A_\bullet\tens E}(\underline{v}\tens e) &=& 
\xi\tens (u\bullet \underline{v})\tens e\ .
\end{eqnarray}
The last derivative looks unusual, until you remember the zero arrived at in 
(\ref{bchdjkshjcgvgxfh}), meaning that differentiating $e$ gives no contribution to the derivative in the second line of (\ref{cvgdjsahgfx}). 
 Note that from (\ref{fcauyxaeffrt}) we have
\begin{eqnarray}  \label{fcauycdsvxaeffrt}
\vartheta_E(w\bullet \underline{v}\tens e) &=& (\la\tens\id)(\id\tens\vartheta_E)(w\tens \underline{v}\tens e) \cr
&&+\ (\sigma_E^{-1}\bullet\id)(\id\tens\vartheta_E)(w\tens \underline{v}\tens e) \ .
\end{eqnarray}
From this and (\ref{cvgdjsahgfx}),
\begin{eqnarray}\label{housadbcu}
(\id\tens
\vartheta_E)\nabla_{ \mathcal{T}\Vec A_\bullet\tens E}(\underline{v}\tens e) &=& 
\xi\tens \vartheta_E(u\bullet \underline{v}\tens e) \cr
&=& \xi\tens (\la\,\tens\id) (u\tens\vartheta_E(\underline{v}\tens e)) \cr
&&+\ \xi\tens (\sigma_E^{-1}\bullet\id)(u\tens \vartheta_E(\underline{v}\tens e))\ .
\end{eqnarray}
Remembering that $\coev(1)=\xi\tens u\in\Omega^1 A\tens_A\Vec A$, we obtain
\begin{eqnarray}
(\id\tens\sigma_E^{-1})(\xi\tens u\tens f) &=& (\sigma_E\tens\id)(f\tens \xi\tens u)\ ,
\end{eqnarray}
and from (\ref{cvgdjsahgfx}) we have
\begin{eqnarray}  \label{cvgdjsvdahgfx}
\nabla_{E\tens \mathcal{T}\Vec A_\bullet}(e\tens \underline{v}) &=& 
\nabla_E(e)\tens  \underline{v} + 
\xi\tens (\sigma_E^{-1}\bullet\id)(u\tens e\tens \underline{v}) \ .
\end{eqnarray}
By definition of $\xi\tens u$,  $\xi\tens u\,\la\,e=\nabla_E(e)$, and we have proved the following result:

\begin{propos}  \label{bchudsabvuyas}
For the coevaluation $\coev:A\to \Omega^1 A\tens_A \Vec A$, 
the map $\nabla(\underline{v})=\coev(1)\bullet \underline{v}$ gives a 
left covariant derivative on $\mathcal{T}\Vec A_\bullet$. It is also a right module map, so $(\nabla,0)$ is a left bimodule covariant derivative on $\mathcal{T}\Vec A_\bullet$. Given another object $(E,\nabla_E,\sigma_E)$ in ${}_A\mathcal{E}_A$,
the map $\vartheta_E: \mathcal{T}\Vec A_\bullet\tens_A E\to  E\tens_A\mathcal{T}\Vec A_\bullet$ defined in Proposition \ref{jkytfryjhzzdt} is a morphism in ${}_A\mathcal{E}_A$, i.e.\ 
\begin{eqnarray*}
(\id\tens\vartheta_E)\nabla_{ \mathcal{T}\Vec A_\bullet\tens E} &=& 
\nabla_{E\tens \mathcal{T}\Vec A_\bullet} \ \vartheta_E \ .
\end{eqnarray*}
\end{propos}

\section{The centre of a category}    \label{kdhjsvjkhcxr}

Recall that a monoidal category $(\mathcal{C},\tens,\Phi,1_\mathcal{C},l,r)$ means a category, with a functor $\tens:\mathcal{C}\times\mathcal{C}\to \mathcal{C}$, a natural equivalence $\Phi:((\ \tens\ )\tens\ )\to (\  \tens(\ \tens\ ))$ subject to Mac Lane's pentagon coherence identity and an identity object $1_\mathcal{C}$ and associated natural isomorpisms $l:\id\to \id\tens 1_{\mathcal{C}}$ and $r:\id\to 1_\mathcal{C}\tens\id$ compatible with $\Phi$. We refer to \cite{Mac} for details. In the trivially associated case (i.e.\ where $(X\tens Y)\tens Z=X\tens (Y\tens Z)$) we can set $\Phi$ to be the identity. 

For simplicity, rather than necessity, the following definition (\cite[Example~3.4]{Maj:rep} and \cite[Definition~3]{JoyStr:tor}) is given in the trivially associated case.
The {\em centre} $\mathcal{Z}(\mathcal{C})$ of a monoidal category $\mathcal{C}$
(with product $\tens$ and identity object $1_{\mathcal{C}}$) is a category which consists of objects which are pairs
$(X,\varphi)$, where $X$ is an object in $\mathcal{C}$ and $\varphi:X\tens-\Rightarrow -\tens X$ is a natural transformation from the functor $A\mapsto X\tens A$ to the functor $A\mapsto A\tens X$. Given an object $A$ of $\mathcal{C}$, we write this as $\varphi_A:X\tens A \to
A\tens X$. The natural transformation is related to the tensor product by
\begin{eqnarray} \label{bcdsuuhvvvv}
\varphi_{A\tens B}=(\id\tens\varphi_B)(\varphi_A\tens\id)\quad\mathrm{and}\quad 
\varphi_{1_\mathcal{C}}=\id_X\ .
\end{eqnarray}
The centre $\mathcal{Z}(\mathcal{C})$ has morphisms $\alpha:(X,\varphi)\to (Y,\vartheta)$
so that $\alpha:X\to Y$ is a morphism in $\mathcal{C}$ with
\begin{eqnarray}
\vartheta_A\,(\alpha\tens\id_A)=(\id_A\tens \alpha)\,\varphi_A:X\tens A\to A\tens Y\ .
\end{eqnarray}
From this definition we can derive certain facts about $\mathcal{Z}(\mathcal{C})$.

The centre $\mathcal{Z}(\mathcal{C})$ is a monoidal category, with 
\begin{eqnarray}
(X,\varphi)\tens (Y,\vartheta) &=& (X\tens Y,(\varphi\tens\id_Y)(\id_X\tens\vartheta))\ ,
\end{eqnarray}
and identity $(1_{\mathcal{C}}, lr^{-1})$. 
It is also a braided category, with
\begin{eqnarray}
\varphi_Y=\Psi_{(X,\varphi), (Y,\vartheta)} : (X,\varphi)\tens (Y,\vartheta) \to
 (Y,\vartheta) \tens (X,\varphi)\ .
\end{eqnarray}

We now give Example \ref{vujhgcvvjhgcf} for two reasons. Firstly, it may reassure those who are not familiar with the centre construction that it is not too complicated. Secondly, as we shall see it shares relevant features with our main example of the differential operators.

\begin{example}\label{vujhgcvvjhgcf} Let $H$ be a Hopf algebra. 
Consider a category of left modules ${}_H\mathcal{M}$ of a Hopf algebra $H$, with morphisms compatible with the $H$ action.  The category ${}_H\mathcal{M}$
has a tensor product, given by $h\,\la(v\tens w)=h_{(1)}\,\la\,v\tens h_{(2)}\,\la\,w$, where $ h_{(1)}\tens h_{(2)} = \Delta(h)$ is the Sweedler notation for the coproduct $\Delta$ of $H$.

We can consider $(H,\la)$ to be an object in ${}_H\mathcal{M}$, where $\la$ is the left adjoint action $g\,\la\,h=g_{(1)}h\,S(g_{(2)})$, where $S$ is the anipode of $H$. 
{}For all $V\in {}_H\mathcal{M}$, define $\varphi_V:H\tens V\to V\tens H$  by
\begin{eqnarray}\label{bcuifnwjfniw}
\varphi_V(h\tens v)\ =\ h_{(1)}\,\la\,v\tens h_{(2)}\ .
\end{eqnarray}
Then $\varphi_V$ is a morphism as
\begin{eqnarray}
\varphi_V(g\,\la\,(h\tens v)) &=& \varphi_V\big(g_{(1)}h\,S(g_{(2)})\tens g_{(3)}\,\la\,v\big) \cr
&=& g_{(1)}h_{(1)}\,S(g_{(4)})\, g_{(5)}\,\la\,v \tens g_{(2)}h_{(2)}\,S(g_{(3)}) \cr
&=& g_{(1)}h_{(1)}\,\la\,v \tens g_{(2)}h_{(2)}\,S(g_{(3)})  \cr
&=&g\,\la\, \varphi_V(h\tens v)\ .
\end{eqnarray}
The condition that $\varphi$ is compatible with the tensor product, 
\begin{eqnarray*}
\varphi_{V\tens W}\,=\, (\id_V\tens\varphi_W)(\varphi_V\tens\id_W)
: H\tens V\tens W\to V\tens W\tens H\ ,
\end{eqnarray*}
is given by the coassociativity of the coproduct. This means that $(H,\la)$ is in the centre $\mathcal{Z}({}_H\mathcal{M})$. If the antipode $S$ of $H$ is bijective, then $\varphi$ is invertible, as $\varphi_V^{-1}(v\tens h)=h_{(2)}\tens S^{-1}(h_{(1)})\,\la\,v$. 

A little more calculation shows that the product $\mu:H\tens H\to H$ is a morphism in the category, and that this makes $(H,\la)$ into an algebra in the centre $\mathcal{Z}({}_H\mathcal{M})$, as 
the required extra condition for the product to be a morphism in $\mathcal{Z}({}_H\mathcal{M})$ is
\begin{eqnarray*}
\varphi_V(\mu\tens\id_V)\,=\,(\id_V\tens \mu)(\varphi_V\tens\id_H)(\id_H\tens \varphi_V)
:H\tens H\tens V\to V\tens H\ ,
\end{eqnarray*}
and is given by the compatibility between the product and coproduct of $H$. 

\end{example}

\begin{theorem} \label{cbadiosiuvcuycf}
The $A$-bimodule $\mathcal{T}\Vec A_\bullet$ with the bimodule 
covariant derivative given in Proposition \ref{bchudsabvuyas} is in the centre
$\mathcal{Z}({}_A\mathcal{E}_A)$ of ${}_A\mathcal{E}_A$ (see Definition \ref{cvgjfxtzhc}), using the natural transformation 
$\vartheta:\mathcal{T}\Vec A_\bullet\tens_A-\Rightarrow -\tens_A \mathcal{T}\Vec A_\bullet$ defined in Proposition \ref{jkytfryjhzzdt}. In addition $\mathcal{T}\Vec A_\bullet$ is a unital associative algebra in $\mathcal{Z}({}_A\mathcal{E}_A)$, using the product $\bullet$ defined in Proposition \ref{bdhscvkvhj}. 
The natural transformation $\vartheta$ is related to the action of differential fields on modules with connection by Fig.\ \ref{newsig56svvuyd}. 
\end{theorem}
\noindent {\bf Proof:}\quad To see that $\vartheta$ is a natural transformation, see
Proposition \ref{bchudsabvuyas}. 
The equations for the tensor products (\ref{bcdsuuhvvvv}) are given by
Proposition \ref{cvadhjskxsz1} and a brief calculation. 
The compatibility between the product $\bullet$ and $\vartheta$ is given by
Proposition~\ref{cvadhjskxsz2}, and compatibility with the unit is given by
Proposition~\ref{cvadhjskxsz3}. \quad$\square$

One would really like to see the structure of a bialgebroid or a $\times_A$-bialgebra \cite{Tak:gro} or, even better, a $\times_A$-Hopf algebra \cite{Sch:dua} on noncommutative vector fields; see \cite{Boh:Hop}. It is instructive to consider Example \ref{vujhgcvvjhgcf} again. Here we have an algebra
$(H,\la,\varphi,\mu)$ in the centre of the category ${}_H\mathcal{M}$. But, given this, can we get back to the Hopf algebra structure of $H$? The answer is yes, but as far as we can see only indirectly. By indirectly, we mean that $H$ can be Tannaka-Krein reconstructed from the representation category \cite{Ma:Tannaka}, and we start with ${}_H\mathcal{M}$. There is a `direct' construction of the action of $H$ by $h\,\la\,v=(\id\tens\epsilon)\varphi_V(h\tens v)$ (this is assuming that we know the counit $\epsilon$), but it is not obvious to see how to construct $\Delta$  by a similar formula. Likewise, it is not clear how to construct the coproduct on  $\mathcal{T}\Vec A_\bullet$.

\end{document}